\documentclass[10pt,a4paper,reqno]{amsart}%

\usepackage{amsthm,amsmath,amsfonts,amssymb,amsxtra,appendix,bookmark,euscript,delarray,dsfont}
\usepackage[latin1]{inputenc}
\usepackage{xcolor}
\usepackage{graphicx}
\usepackage{tikz}
\usepackage{cite}
\usepackage{systeme}
\usepackage{mathtools}

\makeatletter
\def\@tocline#1#2#3#4#5#6#7{\relax
  \ifnum #1>\c@tocdepth 
  \else
    \par \addpenalty\@secpenalty\addvspace{#2}%
    \begingroup \hyphenpenalty\@M
    \@ifempty{#4}{%
      \@tempdima\csname r@tocindent\number#1\endcsname\relax
    }{%
      \@tempdima#4\relax
    }%
    \parindent\z@ \leftskip#3\relax \advance\leftskip\@tempdima\relax
    \rightskip\@pnumwidth plus4em \parfillskip-\@pnumwidth
    #5\leavevmode\hskip-\@tempdima
      \ifcase #1
       \or\or \hskip 1em \or \hskip 2em \else \hskip 3em \fi%
      #6\nobreak\relax
      \dotfill
      \hbox to\@pnumwidth{\@tocpagenum{#7}}
    \par
    \nobreak
    \endgroup
  \fi}
\makeatother

\theoremstyle{plain}
\newtheorem{theorem}{Theorem}[section]
\newtheorem{thm}[theorem]{Theorem}
\newtheorem{notation}[theorem]{Notation}
\newtheorem{lemma}[theorem]{Lemma}
\newtheorem{corollary}[theorem]{Corollary}
\newtheorem{proposition}[theorem]{Proposition}

\newtheorem{claim}[theorem]{Claim}
\theoremstyle{remark}
\newtheorem{remark}[theorem]{Remark}

\numberwithin{equation}{section}

\newcommand\R{{\ensuremath {\mathbb R} }}
\newcommand\C{{\ensuremath {\mathbb C} }}
\newcommand\N{{\ensuremath {\mathbb N} }}

\renewcommand\phi{\varphi}

\newcommand{\cG}{\mathcal{G}}
\newcommand{\eps}{\epsilon}

\renewcommand{\epsilon}{\varepsilon}

\renewcommand{\ge}{\geqslant}
\renewcommand{\le}{\leqslant}
\renewcommand{\geq}{\geqslant}
\renewcommand{\leq}{\leqslant}
\renewcommand{\hat}{\widehat}
\renewcommand{\tilde}{\widetilde}

\title{Asymptotic $N$-soliton-like solutions of the fractional Korteweg-de Vries equation}
\author{Arnaud EYCHENNE}
\address{Department of Mathematics, University of Bergen, Postbox 7800, 5020 Bergen, Norway }
\email{arnaud.eychenne.waxweiler@gmail.com}


\begin{document}

\maketitle

\begin{abstract}
We construct  $N$-soliton solutions for the  fractional Korteweg-de Vries frfKdV) equation 
$$
\partial_t u - \partial_x\left(|D|^{\alpha}u - u^2 \right)=0, 
$$
in the whole sub-critical range $\alpha \in\left(\frac12,2\right)$. More precisely, if $Q_c$ denotes the ground state solution associated to fKdV evolving with velocity $c$, then given  $0<c_1< \cdots < c_N$, we prove the existence of a solution $U$ of (fKdV)
satisfying  
$$
\lim_{t\to\infty} \| U(t,\cdot) - \sum_{j=1}^NQ_{c_j}(x-\rho_j(t)) \|_{H^{\frac{\alpha}2}}=0,
$$
where $\rho'_j(t) \sim c_j$ as $t \to +\infty$. 

The proof adapts the construction of Martel in the generalized KdV setting [Amer. J. Math. 127 (2005), pp. 1103-1140]) to the fractional case. The main new difficulties are the polynomial decay of the ground state $Q_c$ and the use of local techniques (monotonicity properties for a portion of the mass and the energy) for a non-local equation. To bypass these difficulties, we use symmetric and non-symmetric weighted commutator estimates. The symmetric ones were proved by Kenig, Martel and Robbiano [Annales de l'IHP Analyse Non Lin{\'e}aire 28 (2011), pp. 853-887], while the non-symmetric ones seem to be new. 
\end{abstract}

\tableofcontents

\section{Introduction}


\subsection{The fractional Korteweg-de Vries equation}
We consider the fractional Korteweg-de Vries equation (fKdV), also called the dispersion generalized Benjamin-Ono equation,
\begin{equation}
	\partial_t u -|D|^{\alpha} \partial_{x} u + \partial_{x}(u^2)=0, \quad (t,x)\in \R\times\R,
	\label{GC}
\end{equation}
where $\alpha \in \mathbb R$, $|D|^{\alpha}$ is the Riesz potential of order $-\alpha$, defined by $\mathcal{F}(|D|^{\alpha}u)(\xi)=|\xi|^{\alpha}\mathcal{F}(u)(\xi)$ and $\mathcal{F}$ is the Fourier transform. 

In the cases $\alpha=2$, respectively $\alpha=1$, this equation corresponds to the well-known Korteweg de Vries (KdV), respectively Benjamin-Ono (BO) equations, which are completely integrable (see \cite{lax1968integrals,fokas1983inverse}). In the case $\alpha=0$, one recovers the inviscid Burgers' equation after a suitable change of variable, while the case $\alpha=-1$ corresponds to the Burgers-Hilbert equation. Finally, the cases $\alpha=\frac{1}{2}$ and $\alpha=-\frac12$ are somehow reminiscent of the linear dispersion of the finite depth water waves equation with and without surface tension. In other words, for large frequencies, equation \eqref{GC} corresponds in those cases to the Whitham equations with and without surface tension (see \cite{MR3763731} for more details).

From a mathematical point of view, these equations are also useful to understand the \lq\lq fight\rq\rq \, between  nonlinearity and dispersion. Instead of fixing the dispersion (\textit{e.g.} that of the KdV equation) and increasing the nonlinearity (\textit{e.g.} $u^p\partial_xu$ for the generalized KdV equation), one chooses to fix the nonlinearity $u\partial_xu$ and lower the dispersion, allowing then fractional dispersion of the form $|D|^{\alpha}$, $\alpha<2$. As pointed out by Linares, Pilod and Saut in \cite{linares2014dispersive},  this viewpoint is probably more physical since in many problems arising from physics or continuum mechanics the nonlinearity is quadratic with terms like $(u \cdot \nabla)u$ and the dispersion is in some sense weak. Here will focus on positive values of $\alpha$'s. Note however, that the dynamics for negative $\alpha$'s is quite different with the formation of shocks (see \cite{hur2018wave}, \cite{saut2020wave}, \cite{oh2021gradient}).

Although equation \eqref{GC} is not completely integrable outside of the cases $\alpha=1$ and $2$, it enjoys a hamiltonian structure.  In particular, the mass 
\begin{align}
	M(u)(t):=\int u^2(t,x)dx,
\end{align}
and the energy
\begin{align}
	E(u)(t):=  \frac{1}{2}\int\left(|D|^{\frac{\alpha}{2}}u(t,x)\right)^2dx - \frac{1}{3}\int u^3(t,x) dx. \label{energy}
\end{align}
are  formally preserved by the flow of \eqref{GC}.

Moreover, we have the scaling-translation invariance of \eqref{GC}. Let $u$ be a solution of $\eqref{GC}$ then
\begin{align}
	\forall x_0\in \R, c>0, u_{c}(t,x)=c u(c^{\frac{1+\alpha}{\alpha}} t, c^{\frac{1}{\alpha}}(x-x_0)), \label{scalinv}
\end{align}
is also a solution. A straightforward computation shows that $\|u_{c}\|_{\dot{H}^{s}}=c^{s+\alpha-\frac{1}{2}}\| u \|_{\dot{H}^{s}}$. In particular, equation \eqref{GC} is \emph{mass-critical} for $\alpha=\frac12$ and \emph{energy-critical} for $\alpha=\frac13$. 

In this paper, we focus on the \emph{mass-subcritical} case $\alpha \in \left(\frac12,2\right)$. 
We assume that the initial value problem associated to \eqref{GC} is globally well-posed in the energy space $H^{\frac{\alpha}2}(\mathbb R)$ in the whole subcritical range $\frac12<\alpha<2$, in the sense that for all $u_0 \in H^{\frac{\alpha}2}(\mathbb R)$ and $T>0$, there exists a solution $u \in C([0,T] : H^{\frac{\alpha}2}(\mathbb R))$ of \eqref{GC} satisfying $u(0,\cdot)=u_0$ which is unique in some class $X_T \subset C([0,T] : H^{\frac{\alpha}2}(\mathbb R))$, and that the flow $:u_0 \in  H^{\frac{\alpha}2}(\mathbb R) \mapsto u \in C([0,T] : H^{\frac{\alpha}2}(\mathbb R))$ is continuous. 
Such a result has been proved by Herr, Ionescu, Kenig and Koch in \cite{herr2010differential}  in the range $1 \le \alpha < 2$, extending a previous result of Ionescu and Kenig for the BO equation \cite{ionescu2007global}. For weaker dispersion, the global well-posedness in the energy space has been conjectured through numerical simulations by Klein and Saut \cite{klein2015numerical}  in the whole range $\frac12<\alpha<1$. Progress has been made in this direction: recently Molinet, Pilod and Vento proved in \cite{molinet2018well}  global well-posedness in $H^{\frac{\alpha}2}(\mathbb R)$ for $\frac67<\alpha < 1$ (see also Linares, Pilod, Saut \cite{linares2014dispersive} for former results). Note however, that the problem is still open in the case $\frac12<\alpha<\frac67$.

Finally, we mention some other interesting results concerning the fractional KdV equation with positive dispersion $\alpha$. In \cite{ehrnstrom2019enhanced}, Ehrnstr\"om and Wang proved long time existence for small initial data. Fonseca, Linares and Ponce in \cite{MR3103170} proved some persistence results in weighted Sobolev spaces. Kenig, Ponce and Vega \cite{MR4046209}, Kenig, Ponce, Pilod and Vega in \cite{kenig2020unique} and Riano in \cite{riano2021persistence} proved some unique continuation results, while Mendez in  \cite{mendez2020propagation}, \cite{mendez2020propagation2} proved propagation of regularity results. We also refer to Linares, Pilod and Saut \cite{linares2014dispersive} and Klein and Saut \cite{klein2015numerical} for other results, conjectures and numerical simulations regarding the fractional KdV equation.

\subsection{Solitary wave solutions}

A fundamental property of this equation is the existence of solitary wave solutions of the form 
\begin{equation*}
u(t,x)=Q_c(x-ct) \quad \text{with} \quad Q_c(x) \underset{|x| \to +\infty}{\longrightarrow} 0 ,
\end{equation*}
for $c>0$, where $Q_c(x)=cQ(c^{\frac{1}{\alpha}}x)$ and $Q$ is solution of the non-local ODE
\begin{align}
	|D|^{\alpha}Q+Q-Q^2=0.\label{eqQG}
\end{align}
In other words, $Q_c$ satisfies 
\begin{align}
	|D|^{\alpha}Q_{c} + cQ_{c} - Q_{c}^2=0  \label{eqQG2}.
\end{align}
For some particular values of $\alpha$ the solution of \eqref{eqQG} is explicit and unique (up to translations). For $\alpha=2$,  $Q_{KdV}(x)=\displaystyle\frac{3}{2}\cosh^{-2}\left(\displaystyle\frac{x}{2}\right)$, while for $\alpha=1$, $Q_{BO}(x)=4(1+x^2)^{-1}$. The uniqueness result for BO is non-trivial and was proved by Benjamin \cite{benjamin1967internal} and Amick and Toland \cite{amick1991uniqueness} by combining complex analysis techniques with properties of the harmonic extension of the Hilbert transform.

For the other values of $\alpha$, there does not exist, as far as we know, any explicit formulation of $Q$. However, the existence of solutions of \eqref{eqQG} minimising the functional 
\begin{align}\label{GagliardoNirenberg}
		\displaystyle J^{\alpha}(u)=\frac{\left(\displaystyle\int ||D|^{\frac{\alpha}{2}}u|^2\right)^{\frac{1}{2\alpha}}\left(\displaystyle\int |u|^2 \right)^{\frac{\alpha-1}{2\alpha+1}}}{\displaystyle\int |u|^3}.
		\end{align}
is well-known since the work of Weinstein in \cite{weinstein1987existence} and Albert, Bona and Saut \cite{albert1997model}. Such solutions are called \emph{ground states} solutions of \eqref{eqQG}.  They decay polynomially at infinity (see \cite{kenig2011local}), this property being related to the singularity at the origin of the symbol $|\xi|^{\alpha}$. Moreover, their uniqueness is delicate and was proved by Frank and Lenzmann in \cite{frank2013uniqueness} relying on the non-degeneracy of the kernel of the linearized operator associated to $Q$. Below, we summarize the properties of the ground states of \eqref{eqQG}. 

\begin{thm}[\cite{weinstein1987existence,albert1997model,kenig2011local,frank2013uniqueness}]
	Let $\alpha\in\left(\frac{1}{3},2\right)$. There exists $Q\in H^{\frac{\alpha}{2}}(\R)\cap C^{\infty}(\R)$ such that  
	\begin{enumerate}
		\item (\textit{Existence}) The function $Q$ solves \eqref{eqQG} and  $Q=Q(|x|)>0$ is even, positive and strictly decreasing in $|x|$. Moreover, the function $Q$ is a minimizer of $J^{\alpha}$ in the sense that 
		\begin{equation} 
		J^{\alpha}(Q)=\inf_{u \in H^{\frac{\alpha}2}(\mathbb R)}J^{\alpha}(u).
		\end{equation}
		\item (\textit{Decay}) The function $Q$ verifies the following decay estimate
		\begin{align}
			\frac{1}{C(1+|x|)^{k+1+\alpha}}\leq Q^{(k)}(x)\leq \frac{C}{(1+|x|)^{k+1+\alpha}}, \quad k=0,1,2,\label{decay}
		\end{align} 
	for some $C>0$.
	\item (\textit{Uniqueness}) The even ground state solution $Q=Q(|x|)>0$ of \eqref{eqQG} is unique. Furthermore, every optimizer $v\in H^{\frac{\alpha}{2}}(\R)$ for the Gagliardo-Nirenberg problem \eqref{GagliardoNirenberg} is of the form $v=\beta Q(\gamma(\cdot+y))$ with some $\beta\in\C, \beta\neq0, \gamma>0$ and $y\in \R$.
	\item (\textit{Linearized operator}) Let $\mathcal{L}$ be the unbounded operator defined on $L^{2}(\R)$ by 
	\begin{align}
		\mathcal{L}u=|D|^{\alpha}u+u-2Qu.\label{oplinea}
		\end{align}
	Then, the continuous spectrum of $\mathcal{L}$ is $[1,+\infty)$, $\mathcal{L}$ has one negative eigenvalue $\mu_0$, associated to an even eigenfunction $W_0>0$, and $\text{ker} \, \mathcal{L}= \text{span} \, \{ Q'\}$.
	\end{enumerate}
	\end{thm} 

 \begin{remark} 
 The uniqueness problem for the solutions of \eqref{eqQG} which are not ground states is still an open question when $\alpha \neq 1$. 
 \end{remark}

 These solitary waves are orbitally stable under the flow of \eqref{GC} (see Linares, Pilod and Saut \cite{linares2015remarks} and \cite{arnesen2015existence,pava2018stability} for other proofs) in the mass sub-critical range $ \alpha \in\left(\frac12,2\right)$. They were proven to be linearly unstable in the the mass super-critical range $ \alpha \in\left(\frac13,\frac12\right)$ (see \cite{pava2018stability}).
 
The case $0<\alpha<\frac{1}{3}$ is energy super-critical. It has been proved in Theorem 4.1 \cite{linares2014dispersive} that fKdV does not posses solitary waves moving to the right and belonging to the energy space $H^{\frac{\alpha}{2}(\R)}$.
 
Sometimes, we also call these solutions \emph{solitons} even though they are not known to have elastic interactions outside of the integrable case $\alpha=1$.
 
\subsection{N-soliton solutions} 

An important conjecture for nonlinear dispersive equations is to prove the \emph{soliton resolution property}, which states that arbitrary initial data eventually resolve over time into a finite or infinite sum of solitary waves and an oscillatory remainder of essentially linear type. It has been proved in the KdV case for sufficiently smooth and decaying initial data by using the complete integrable structure (see \cite{eckhaus1983emergence}). Note however that despite some recent progress (see \cite{wu2017jost}, \cite{sun2021complete}), it is still an open problem for the Benjamin-Ono equation on the line. 

For KdV type equations, we are still far from a complete understanding of this phenomenon. In this direction, an important question is to construct solutions behaving like a superposition of $N$ solitary waves at infinity. Indeed, such objects are expected to be universal attractors in the region $x>0$ for any smooth and decaying solutions at infinity. These solutions, also called \emph{$N$-soliton solutions} by abuse of language, were first constructed by Martel in \cite{martel2005asymptotic2} for the sub-critical and critical gKdV equations by adapting the construction by Merle in \cite{merle1990construction} of solutions blowing up at $k$ given points for the critical nonlinear Schr\"odinger equation to the KdV setting, and by relying on the energy methods by Martel, Merle and Tsai \cite{martel2002stability}. This construction was extended to the super-critical gKdV equations by C\^ote, Martel and Merle \cite{cote2011construction} (see also Combet \cite{combet2010multi}).

 For the fractional KdV equations, outside of the case $\alpha=1$, no result concerning  construction of $N$-solitary wave solutions at infinity seems to be known.  Of course, in the case $\alpha=1$, the $N$-soliton solutions of the Benjamin-Ono equation are explicit by using inverse scattering method \cite{lax1968integrals},\cite{eckhaus1983emergence}, \cite{bock1979two}, \cite{nakamura1979backlund}, \cite{matsuno1979exact}. These $N$-solitons were also proved to be orbitally stable by Neves and Lopes \cite{neves2006orbital} and Gustafson, Takaoka and Tsai \cite{MR2492606} and even asymptotically stable by Kenig and Martel \cite{kenig2009asymptotic}. 

The main result of this paper states the existence of such $N$-soliton solutions for any given set of velocities $0<c_1<c_2<\cdots<c_N$.  
\begin{thm}\label{exsol}
	We assume $\alpha\in\left(\frac{1}{2},2\right)$. Let $N\in\mathbb{N}$, $0< c_1 < \cdots < c_N<+\infty $. Then, there exist some constants $ T_0>0, C_0>0$, $N$ functions $\rho_1,\cdots, \rho_N\in C^{1}([T_0,+\infty))$ and $U\in C^{0}([T_0,+\infty):H^{\frac{\alpha}{2}}(\R))$ solution of \eqref{GC} such that, for all $t \ge T_0$,
	\begin{align} \label{theo:est1}
		\Bigg\lVert U(t,\cdot) - \displaystyle\sum_{j=1}^{N}Q_{c_j}(\cdot-\rho_j(t)) \Bigg\rVert_{H^{\frac{\alpha}2}}\leq \displaystyle\frac{C_0}{ t^{\frac{\alpha}{2} }},
	\end{align} 
	\begin{align} \label{theo:est2}
		\left|\rho_j(t)-c_jt \right| \le t^{1-\frac{\alpha}4} \quad  \text{and} \quad |\rho_j'(t) - c_j | \le \frac{C_0}{t^{\frac{\alpha}{2} }} ,
	\end{align}
	for all $j \in \{1,\cdots,N \}$.  
\end{thm}

\begin{remark}
	Due to the polynomial decay of the error in \eqref{theo:est2}, it is not clear whether we are  in the case of \emph{strong interactions} or not, and thus whether the asymptotic of $\rho_j(t)$ in \eqref{theo:est1} may be more complicated than just $c_jt$. The construction introduced to understand the strong interaction phenomenon, see for example \cite{vinh2017strongly}, could provide tools to obtain a better estimate on $\rho_j'$ \eqref{theo:est2}. Recently, in a joint work with Valet in \cite{eychenne2022strongly}, we proved the existence of strongly interacting $2$-soliton for the fractional modified Korteweg de Vries equation. An adaptation of the construction obtained in this paper could potentially  help to improve the result in Theorem \ref{exsol}. 
\end{remark}

\begin{remark}
	The construction in the case $\alpha\in\left(\frac{1}{2},\frac{6}{7}\right)$ is conditional to the well-posedness of the equation in the energy space, which is still an open problem for this range of $\alpha$'s.
\end{remark}

\begin{remark}
	Since the solitons $Q_{c_j}$ are smooth, we expect the $N$-soliton solution $U$ to  be also smooth. To obtain the convergence \eqref{theo:est1} in higher Sobolev norm, one  possibility could be to follow the method in \cite{martel2005asymptotic2} based on a Gronwall argument and modified energies. However, it is not clear how to construct such modified energies in the non-local setting. \\
	 Related to this question, the uniqueness of these $N$-soliton solutions is an interesting open problem. 
\end{remark}

Similar construction of $N$-soliton-like solutions have already been performed for other nonlinear dispersive equations. Outside of the gKdV equations commented above, we refer to Martel and Merle \cite{martel2006multi} and C{\^o}te, Martel and Merle \cite{cote2011construction} for the non-linear Schr{\"o}dinger (NLS) equation, and more recently to Ferriere for the logarithmic-NLS equation  \cite{ferriere2021existence}. We also refer to the works of Martel and Merle  \cite{martel2016construction} and Jendrej \cite{jendrej2020construction} for the wave equation, to the works of  C{\^o}te-Mu\~ noz \cite{cote2014multi}, Bellazzini, Ghimenti, Le Coz \cite{bellazzini2014multi} and C{\^o}te, Martel \cite{cote2018multi} for the Klein-Gordon equation, to the work Rousset-Tzvetkov   \cite{ming2015multi} for the water-waves equation, and to the work of Valet \cite{valet2021asymptotic} for the Zakharov-Kuznetsov equation.  

A different method of construction of multi-solitons, based on the fixed point argument of Merle in \cite{merle1990construction}, has been introduced by Le Coz, Li and Tsai in \cite{le2015fast}  for the NLS equation to construct an infinite sum of solitary waves. This strategy has also been used by Chen for the wave equation \cite{chen2018multisolitons} and by Van Tin for the derivative NLS \cite{vantin2021construction}.

Recently Jendrej, Kowalczyk and Lawrie introduced in \cite{jendrej2019dynamics} a new version of the Liapunov-Schmidt reduction in the setting of dispersive equations to derive a complete classification of all kink-antikink pairs in the strongly interacting regime for the classical nonlinear scalar field models on the real line.

Finally, let us observe that the result of Theorem \ref{exsol} would be the first step to study the collision of multi-soliton solutions in the cases $\alpha \in \left(\frac12,2\right)$, $\alpha \neq 1$. We refer to the works of Martel and Merle \cite{martel2011description} and \cite{martel2011inelastic}  for the study of the inelastic collision of two solitons of the quartic KdV equation.

\subsection{Outline of proof of Theorem \ref{exsol}}
The proof of Theorem \ref{exsol} follows the strategy of \cite{merle1990construction}, \cite{martel2005asymptotic2}, \cite{martel2002stability}.
After fixing a sequence of time $(S_n)\nearrow+\infty$, one considers the sequence $(u_n)$ of solutions to \eqref{GC} evolving from the initial data $\displaystyle\sum_{j=1}^{N}Q_{c_j}(\cdot-c_jS_n)$ at time $S_n$. As long as the solution remains sufficiently close to the sum of $N$ solitary waves, one introduces modulated translation parameters $(\rho_{j,n}(t))_{j=1}^N$ allowing to satisfy suitable orthogonality conditions. The goal is to obtain backwards uniform estimates for the difference $u_n(t)-\displaystyle\sum_{j=1}^{N}Q_{c_j}(\cdot-\rho_j(t))$ on some time interval $ [T_0,S_n]$, for some $T_0$ independent of $n$. The $N$-soliton is then obtained by letting $n \to +\infty$ and using a compactness argument. Moreover, it is worth to observe that the uniform estimate relies on monotonicity properties for suitable  portions of the mass and  the energy of the solution.

Compared to the previous constructions, we have to deal here with two major new difficulties. Firstly, due to the singularity at the origin of the symbol $|\xi|^{\alpha}$ related to the non-local operator $|D|^{\alpha}$, the solitary waves have only polynomial decay\footnote{The decay of the solitary waves of gKdV is always exponential.} of order $(1+|x|)^{-(1+\alpha)}$. As a consequence the uniform estimates on the parameters $\rho_{j,n}(t)$ are only polynomial and thus cannot be integrated directly.  Relying on a topological argument introduced in \cite{cote2011construction}, we need then to adapt carefully the initial data of the translation parameters $\rho_{j,n}$ at time $S_n$ to be able to close the bootstrap estimates.

Secondly, observe that the monotonicity techniques introduced by Martel and Merle for gKdV are local in space, and are therefore tailored for differential but not integral (nonlocal) equations. To adapt these techniques to the fKdV equations, one need to use suitable weighted commutator estimates (see Lemma \ref{commG}). Those estimates were introduced in the symmetric case by Kenig and Martel \cite{kenig2009asymptotic} in the case $\alpha=1$ and Kenig, Martel and Robbiano \cite{kenig2011local} for the general case $0<\alpha<2$ (see also \cite{martel2017construction} for an application to the critical modified Benjamin-Ono equation). Note however that to derive the monotonicity property of the energy, one needs a non-symmetric version of these estimates (see estimates \eqref{nsc1G}-\eqref{nsc2G}), whose proof is based on pseudo-differential calculus and follows the one of Kenig, Martel and Robbiano for the symmetric case.

\bigskip

The paper is organised as follows: in Section \ref{sec:cons},  we modulate the geometrical translation parameters for a solution close to $N$ solitary waves, set up the bootstrap setting and close the construction of the $N$-soliton solution after assuming the main bootstrap estimate. In Section \ref{sec:weight}, we state several weighted estimates whose proofs are given in the appendices. These weighted estimates are useful to derive the monotonicity properties and to prove the bootstrap estimate in Section \ref{sec:mono}.

\subsection{Notation}

	\begin{enumerate}
		\item From now on, $C$ will denote a positive constant changing from line to line and  independent of the different parameters. We also denote by $C_{*}$ a positive constant changing from line to line and depending only on the parameters $\{c_1,\cdots,c_N\}$.
		\item Unless stated otherwise, all the integrals will be over $\R$ with respect to the space variable.
		\item For $x\in \R^{N}$, we recall the definition of the Japanese brackets $\langle x \rangle:= \sqrt{1+|x|^2}$.
		\item We denote by $\| f\|_{L^{p}}:=\displaystyle\left(\int |f|^p\right)^{\frac{1}{p}}$ and $\| f\|_{H^{s}}:=\|\langle\cdot\rangle^{\frac{s}{2}}\mathcal{F}(f)(\cdot) \|_{L^2}$, where $\mathcal{F}(f)$ is the Fourier transform of $f$. Finally, $\mathcal{S}(\R)$ denotes the  Schwartz space of real-valued functions.
		\item We fix $0<c_1<\cdots<c_N$ and we set $\beta=\frac{1}{2}\min(c_1,c_2-c_1,\cdots, c_N-c_{N-1})$.
	\end{enumerate}


\section{Construction of the asymptotic $N$-soliton} \label{sec:cons}

\begin{notation}
	\begin{enumerate}
		\item For $L>0$ and $N\in\N$, we define 
		\begin{align}
			\R^{N}_L=\{(y_j)_{j=1}^{N}\in \R^N : y_j-y_{j-1}>L, \forall j\in \{2,\cdots,N\}\}\label{espace}.
		\end{align}
		\item For $Y=(Y_j)_{j=1}^N\in \R^N_L$, we denote 
		\begin{align}
			R_Y(x)=\sum_{j=1}^{N}R_{Y,j}(x):=\sum_{j=1}^NQ_{c_j}(x-Y_j).\label{RY}
		\end{align}
		\item Let $M=(m_{i,j})_{i,j=1}^{N}\in M_N(\R)$, be a  $N\times N $ matrix. 
	\end{enumerate}
	
\end{notation}


\subsection{Modulation of the geometrical parameters}

\begin{proposition}[Modulation]\label{modg}
	There exist $L_1,\gamma_1,T_1>0$ such that the following is true. Assume that $u$ is a solution of \eqref{GC} satisfying that for $L>L_1$, $0<\gamma<\gamma_1$, $S>t^*>T_1$,
	\begin{align}
		\displaystyle\sup_{t^*\leq t\leq S}\left(\inf_{(Y_j)_{j=1}^{N}\in\R^N_L } \bigg\|u(t,\cdot) - \displaystyle\sum_{j=1}^{N} Q_{c_j}(\cdot-Y_j) \bigg\|_{H^\frac{\alpha}{2}}    \right)<\gamma \label{alpha}.
	\end{align}
	Then, there exist $N$ unique $C^1$ functions $\rho_j: [t^*,S]\longrightarrow \R$, $j\in\{1,...,N\}$, such that 
	\begin{align}
		\eta(t,x)=u(t,x) - R(t,x) ,\label{defeta}
	\end{align}
	where
	\begin{align}
		R(t,x) =\displaystyle\sum_{j=1}^N R_j(t,x):= \displaystyle\sum_{j=1}^N Q_{c_j}(x-\rho_j(t) ), \label{Rg}
	\end{align}
	satisfies the following orthogonality conditions 
	\begin{align}
		\displaystyle\int(\partial_x R_j) \eta=0, \quad \forall j\in \{1,...,N \} , \forall t\in [t^*,S]. \label{Ortg}
	\end{align}
	Moreover, for all $ t \in [t^*,S] $ 
	\begin{align}
		\| \eta(t,\cdot) \|_{H^{\frac{\alpha}{2}}} &\leq C\gamma,\label{size}\\
				\inf_{j\in\{1,...,N-1\}} \left(\rho_{j+1}(t)-\rho_j(t)\right) &\geq \frac{L}{2} .\label{vitinf2}
	\end{align}
\end{proposition}


\begin{remark}
	A solution $u$ satisfying \eqref{alpha} lives for all time $t\in[t^*,S]$ in the tube  
	\begin{align}
		\mathcal{T}_{\gamma,L}:=\bigg\{v\in H^{\frac{\alpha}{2}}(\R) :\displaystyle\inf_{Y_j\geq Y_{j-1} + L }\bigg\|v - R_Y \bigg\|_{H^\frac{\alpha}{2}} <\gamma\bigg\}.\label{tube}
	\end{align}
	The proof of Proposition \ref{modg} is an application of the implicit function theorem to the functional 
	\begin{align}
		\Phi : \mathcal{T}_{\gamma,L}\times \R^N_L &\to \R^N \label{fonct} \\
		(v,Y)&\mapsto\left(\displaystyle\int\left(v-R_Y \right)  Q'_{c_j}(\cdot-Y_j)\right)_{j=1}^{N}.\notag
	\end{align}
	Note that a direct application of the implicit function theorem at the point $(R_Y,Y)$ for $Y\in \R_L^N$ would imply
	\begin{align}
		\forall Y\in\mathbb{R}_{L}^N, \exists\varepsilon_Y>0,\exists!\left(\rho_j\right)_{j=1}^N\in C^1(\mathcal{T}_{\gamma,L}\cap B(R_Y,\varepsilon_Y)):\R), \label{IFT1}
	\end{align} 
	such that $\Phi(v,(\rho_j(v))_{j=1}^N)=0$ for $v\in B(R_Y,\eps_Y)$. This would not be enough to conclude the proof of \eqref{Ortg} due to the lack of control of $\varepsilon_{u(t)}$ uniformly in $[t^*,S]$. Indeed, an application of \eqref{IFT1} to a solution $u$ satisfying $u(S,\cdot)=\displaystyle\sum_{j=1}^NQ_{c_j}(\cdot-\rho_{j,n}^{\text{in}})$ for $(\rho_{j,n}^{\text{in}})\in\R_{L}^{N}$, and a continuity argument would provide the existence of $\eps>0$ such that $u(t,\cdot)\in B:=B(\displaystyle\sum_{j=1}^NQ_{c_j}(\cdot-\rho_{j,n}^{\text{in}}),\eps)$ for all $t\in (t_1,S]$, where $t_1$ is the first time before $S$ with $u(t_1,\cdot)\notin B$. Nevertheless, nothing would guarantee that $u(t_1,\cdot)$ belongs to a ball $B(R_Y,\eps_Y)$ for some $Y\in R_L^N$.   \\
	
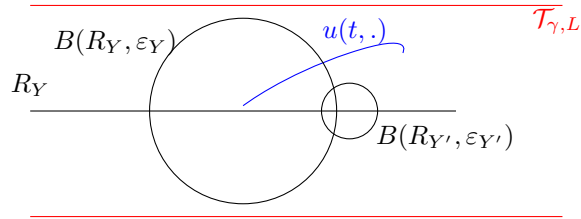
	\begin{figure}[h!]
		\centering
		\begin{tikzpicture}[scale=0.7]
			\draw [red]   (-1,2) -- (9,2);
			\draw [red]  (-1,-2) -- (9,-2);
			\draw [blue, xshift=3cm, yshift=0.1cm] plot [smooth, tension=1.5] coordinates { (0,0) (2,1) (3,1)};
			\draw [black] plot [smooth, tension=1] coordinates { (-1,0)  (7,0)};
			\draw [xshift=3cm] (0,0) circle (50pt);
			\draw [xshift=3cm] (2,0) circle (15pt);
			\draw (5.1,1.5) node[blue]{$u(t,.)$};
			\draw (8.9,1.7) node[red]{$\mathcal{T}_{\gamma,L}$};
			\draw (-1,0.5) node[black]{$R_Y$};
			\draw (0.6,1.3) node[black]{$B(R_Y,\varepsilon_Y)$};
			\draw (6.8,-0.5) node[black]{$B(R_{Y'},\varepsilon_{Y'})$};
		\end{tikzpicture}
	\end{figure}
	To bypass this difficulty, we will use the following quantitative version of the implicit function theorem (see  section 2.2 in \cite{chow2012methods}). We refer to \cite{MR2322685} Lemma 3,  \cite{jendrejnonexistence} Lemme 3.3, \cite{gustafson2006effective} Proposition 3 ,\cite{zaag2012isolatedness} Proposition 3.1  for applications of this theorem in a similar context. 
\end{remark}

\begin{thm}\label{tfi}
	Let $X,Y$ and $Z$ be Banach spaces, $x_0\in X, y_0\in Y$, $\gamma,\delta>0$ and $\Phi:B(x_0,\gamma)\times B(y_0,\delta)\longrightarrow Z$ be continuous in $x$,   continuously differentiable in $y$,  satisfy $\Phi(x_0,y_0)=0$, $M_0:=d_y\Phi(x_0,y_0)$ has a bounded inverse in $\mathcal{L}(Z,Y)$. Assume moreover that
	\begin{align}
		\lVert M_0 - d_y\Phi(x,y) \rVert_{\mathcal{L}(Y,Z)}\leq \frac{1}{3} \lVert M_0^{-1} \rVert^{-1}_{\mathcal{L}(Z,Y)}, \quad \forall x\in B(x_0,\gamma), y\in B(y_0,\delta),\label{FI1}
	\end{align}
	\begin{align}
		\lVert \Phi(x,y_0) \rVert_{Z}\leq \frac{\delta}{3} \lVert M_0^{-1} \rVert^{-1}_{\mathcal{L}(Z,Y)}, \quad \forall  x\in B(x_0,\gamma).\label{FI2}
	\end{align}
	Then there exists $y\in C^1(B(x_0,\gamma):B(y_0,\delta))$ such that for $x\in B(x_0,\gamma)$, $y(x)$ is the unique  solution of the equation $\Phi(x,y(x))=0$ in $B(x_0,\gamma)$.
\end{thm}

Before giving the proof of Proposition \ref{modg}, we need the following lemma.

\begin{lemma}\label{estnonmod}
	There exist $C>0$, $L_2>0$, such that for all $L>L_2$, and all $Y=(Y_{j})_{j=1}^{N}\in\R^{N}_L$ we have 
	\begin{align}
		\bigg| \int (\partial_xR_{Y,j})(\partial_xR_{Y,k}) \bigg|\leq \frac{C}{1+L^{2+\alpha}},\quad j\neq k, \label{estnonmod1}
	\end{align} 
	with $R_{Y,j}$ defined in \eqref{RY} .
	
	Moreover, let $(\rho_j)_{j=1}^{N}\in C^1([t^{*},S]:\R)$ satisfying $\rho_{j+1}-\rho_{j}\geq \frac{L}{2}$ for all $j,k\in\{1,\cdots,N-1\}$, with $j\neq k$, then 
	\begin{align}
		\displaystyle\bigg|\int (\partial_xR_j) (\partial_xR_k)\bigg|&\leq \frac{C}{1+ L^{2+\alpha}},  \label{estnonmod2}\\
		\displaystyle\bigg|\int R_j (\partial_x^2R_k)\bigg|&\leq \frac{C}{1+ L^{1+\alpha}},  \label{estnonmod3}\\
		\displaystyle\bigg|\int R_jR_k (\partial_x^2R_l)\bigg|&\leq \frac{C}{1+ L^{1+\alpha}},  \label{estnonmod4}
	\end{align}
	with $R_j$ defined in \eqref{Rg}.
	
	Furthermore, if the functions $(\rho_j)_{j=1}^{N}$ satisfy $|\rho_{j+1}(t)-\rho_j(t)|\geq \beta t$, with $\beta>0$, then  
		\begin{align}
			\displaystyle\bigg|\int(\partial_xR_j) (\partial_xR_k)\bigg|&\leq \frac{C}{(\beta t)^{2+\alpha}},  \label{estnonmod5}\\
				\displaystyle\bigg|\int R_j (\partial_x^2R_k)\bigg|&\leq \frac{C}{( \beta t)^{1+\alpha}} . \label{estnonmod6}\\
				\displaystyle\bigg|\int R_j R_k (\partial_x^2R_l)\bigg|&\leq \frac{C}{( \beta t)^{1+\alpha}} . \label{estnonmod7}
	\end{align}
\end{lemma} 

\begin{proof}[Proof of Lemma \ref{estnonmod} ]
	By symmetry, we can suppose $j<k$. Let $\Omega:=\left\{x\in \R: x<\frac{Y_j+Y_k}{2}\right\}$. By \eqref{decay} and $Y_k-Y_j>L$, we deduce that
	\begin{align*}
		\bigg|\int_{\Omega} (\partial_xR_{Y,j})(\partial_xR_{Y,k})\bigg|\leq \frac{C}{1+\left(Y_k-\frac{Y_j+Y_k}{2}\right)^{2+\alpha}}\int_{\Omega} |\partial_xR_{Y,j}| \leq \frac{C}{1+L^{2+\alpha}}.
	\end{align*}
	On the other hand, by \eqref{decay} and $Y_j-Y_k<-L$, we get on $\Omega^{c}$
	\begin{align*}
		\bigg|\int_{\Omega^c} (\partial_xR_{Y,j})(\partial_xR_{Y,k})\bigg|\leq \frac{C}{1+\left(\frac{Y_j+Y_k}{2}-Y_j\right)^{2+\alpha}}\int_{\Omega^c} |\partial_xR_{Y,k}|\leq \frac{C}{1+L^{2+\alpha}},
	\end{align*}
	which concludes \eqref{estnonmod1}. To prove the other estimates, we use the same argument on $\Omega:=\{x\in \R: x<\frac{\rho_j+\rho_k}{2} \}$ with the estimates \eqref{decay}, $\rho_{j+1}-\rho_j\geq \frac{L}{2}$ for \eqref{estnonmod2}, \eqref{estnonmod3}, \eqref{estnonmod4} and  \eqref{decay},$|\rho_{j+1}(t)-\rho_j(t)|\geq \beta t$ for \eqref{estnonmod5}, \eqref{estnonmod6}, \eqref{estnonmod7}.
\end{proof}

\begin{proof}[Proof of Proposition \ref{modg}]We decompose the proof in two steps. First, by using Theorem \ref{tfi}, we show that we can find $N$ unique functions $\rho_j$ continuous on $\mathcal{T}_{\gamma,L}$, defined in \eqref{tube}, satisfying \eqref{Ortg} - \eqref{vitinf2}. To obtain the regularity of the functions, we use  the Cauchy-Lipschitz theorem.

	\noindent\emph{First step : existence of the functions $\rho_j$.}
	We recall the definition of $\R_L^N$ and $R_Y$  given respectively in \eqref{espace} and \eqref{RY}.  First, we check that the functional $\Phi$ defined in \eqref{fonct} satisfies
	 the hypotheses of Theorem \ref{tfi}. It is clear that, for all $Y\in\R^{N}_L$, 
	\begin{align*}
		\Phi\left(R_Y,Y\right)=0.
	\end{align*}
	Let us define
	\begin{align*}
		M_0=M_0(R_Y,Y) := d_Y\Phi\left(R_Y,Y\right)=A+B, 
	\end{align*}
	where
	$$
	A:=\begin{pmatrix}
		\int(Q'_{c_1})^2 &0& \cdots& & 0\\
		0& \int (Q'_{c_2})^2&0& \cdots& 0\\
		\vdots & & \ddots & &\vdots\\
		\vdots & &  & \ddots &  \vdots\\
		0&\dots & &0 &\int (Q'_{c_N})^2
	\end{pmatrix},
	$$
	and 
	$$
	B=B(Y):=\begin{pmatrix}
		0&\mathcal{Q}_{1,2}& \cdots& &\mathcal{Q}_{1,N} \\
		\mathcal{Q}_{2,1}& 0&\mathcal{Q}_{2,3}& \cdots& \mathcal{Q}_{2,N}\\
		\vdots & & \ddots & &\vdots\\
		\vdots & & & \ddots &\vdots\\
		\mathcal{Q}_{N,1}&\dots &  &\mathcal{Q}_{N,N-1} &0
	\end{pmatrix},
	$$
	with $\mathcal{Q}_{j,k}:=\displaystyle \int \partial_xR_{Y,j}\partial_xR_{Y,k}$.
	The matrix $A$ is invertible, and by \eqref{estnonmod1}, we get for all $L>L_2$
	$$
	|\mathcal{Q}_{j,k}|\leq \frac{C}{1+ L^{2+\alpha}}.
	$$
	Then for $L>L_3$, with $L_3$ big enough, $M_0$ is invertible. Moreover the matrix $A$ is independent of $Y\in\R_L^N$, and $\displaystyle\lim_{L\to\infty}\|B\|_{\infty} \rightarrow 0$. Then, there exists $\kappa$ independent of $L>1$ such that for all $Y\in\R_L^N$
	$$
	\| M_0\left(R_Y,Y \right)^{-1}\|_{{\infty}}=\| A^{-1}(Id + B(Y)A^{-1})^{-1} \|_{{\infty}}\leq \| A^{-1}\|_{\infty} \sum_{n=0}^{\infty}\bigg[C\frac{\|A^{-1}\|_{\infty}}{1+L^{2+\alpha}}\bigg]^{n} \leq  \kappa.
	$$
	Thus, to verify the conditions \eqref{FI1},\eqref{FI2}, since $\|\cdot\|_{\mathcal{L}(\R^N,\R^N)}\leq \|\cdot\|_{\infty}$, it suffices to prove that
	\begin{align}
		\| M_0 - d_Y\Phi(v,Z) \|_{\infty}\leq \frac{1}{3} \kappa^{-1}, \text{ for } v\in B(R_Y,\gamma), Z\in B(Y,C_1\gamma),\label{FI11}
	\end{align}
	\begin{align}
		\| \Phi(v,Y) \|_{\infty}\leq \frac{C_1\gamma}{3} \kappa^{-1}, \text{ for } v\in B(R_Y,\gamma),\label{FI22}
	\end{align}
	for a positive constant $C_1$ to be chosen later.
	First, we show \eqref{FI22}. Let $j\in\{1,...,N\}$, by Cauchy-Schwarz and since $v\in B(R_Y,\gamma)$ 
	\begin{align*}
		|\Phi_j(v,Y)|\leq \displaystyle\int|v - R_Y | |\partial_xR_{Y,j}|\leq \|v - R_Y \|_{L^2} \| \partial_xR_{Y,j} \|_{L^2}\leq \gamma \| \partial_xQ_{c_j}\|_{L^2},
	\end{align*}
	which implies $\eqref{FI22}$ by choosing $C_1=\displaystyle3\kappa\sup_{j}\|\partial_xQ_{c_j}\|_{L^2}$. Since the constant $C_1$ does  not play any role in the rest of the paper, we write $C$ instead of $C_1$.

	Now, let us verify \eqref{FI11}. First we define 
	$$
	\mathcal{Q}^{*}_{j,k}(Z,Y):= \int \left((\partial_xR_{Z,j})(\partial_xR_{Z,k}) - (\partial_xR_{Y,j})(\partial_xR_{Y,k})\right),
	$$
	with $Z=(Z_j)_{j=1}^N,Y=(Y_j)_{j=1}^N\in\R^{N}_{L}$, and 
	$$
	\mathcal{P}_k(v,Z)=\displaystyle\int\left(v(x) - R_{Z}(x) \right) \partial^2_xR_{Y,k}.
	$$
	with $Z=(z_j)\in\R^{N}_{L}$, $v\in B(R_Y,\gamma)$. We have
	\begin{align*}
		\| M_0 - d_Y\Phi(v,Z) \|_{\infty} = \Bigg\| 
		\begin{pmatrix}
			\mathcal{P}_1(v,Z)&\mathcal{Q}^{*}_{2,1}(Z,Y)  & \cdots& &\mathcal{Q}^{*}_{N,1}(Z,Y) \\
			\vdots & & \ddots & &\vdots\\
			\vdots & & & \ddots &\vdots\\
			\mathcal{Q}^{*}_{1,N}(Z,Y)&\dots &  &\mathcal{Q}^{*}_{N-1,N}(Z,Y)   &\mathcal{P}_N(v,Z)
		\end{pmatrix}
		\Bigg\|_{\infty}.
	\end{align*}
	For $j\in\{1,\cdots,N\}$, by Cauchy-Schwarz since $v\in B(R_Y,\gamma)$
	\begin{align}
		|\mathcal{P}_j(v,Z)|\leq \displaystyle\int|v - R_Z | | \partial^2_xR_{Y,j} |\leq \|v - R_Z \|_{L^2} \| \partial^2_xR_{Y,j} \|_{L^2}\leq C\gamma \label{Pk}.
	\end{align}
	For $j,k\in\{1,\cdots,N\}, j\neq k$, by \eqref{estnonmod1}
	\begin{align}
		|\mathcal{Q}^{*}_{j,k}(Z,Y)|\leq \frac{C}{1+L^{2+\alpha}}. \label{Qkj}
	\end{align}
	Gathering \eqref{Pk} and \eqref{Qkj}, we get
	\begin{align*}
		\| M_0 - D_Y\Phi(v,Y) \|_{\infty} \leq C\gamma + \frac{C}{1+L^{2+\alpha}} \leq \frac{1}{3}\kappa^{-1}.
	\end{align*}
	for $\gamma<\gamma_2$ small enough and $L>L_4$ big enough. 
	
	Then for $L>\max(L_2,L_3,L_4)$ and $\gamma<\gamma_2$ we deduce from Theorem \ref{tfi} the existence and uniqueness  of  $(\rho_j)_{j\in\{1,...,N\}}$  in $C^1(B(R_Y,\gamma) : B(Y,C\gamma))$ satisfying \eqref{Ortg}. Moreover, since $\gamma$ can be chosen independently of $Y\in\R^{N}_{L}$, we can extend by uniqueness $(\rho_j)_{j\in\{1,...,N\}}$ to the whole tube $\mathcal{T}_{\gamma,L}$ defined in \eqref{tube}. Furthermore, for all $v\in \mathcal{T}_{\gamma,L}$, there exists $Y=(Y_j)\in \R^{N}_{L}$ such that $(\rho_j(v))_{j\in\{1,...,N\}}\in B(Y,C\gamma)$. Therefore 
	\begin{align}
		|\rho_{j+1}(v)-\rho_{j}(v)|\geq |Y_{j+1}-Y_j| - |\rho_{j+1}(v)-Y_{j+1}| - |\rho_{j}(v)-Y_{j}|\geq L-2C\gamma\geq \frac{L}{2} .  \label{tr1}
	\end{align}
	Now, by abuse of notation, we define $\rho_{j}(t):=\rho_j(u(t,\cdot))$. Then it is clear that $\rho_j$ is $C^{0}([t^{*},S]:\R)$ since $u(t,\cdot)\in C^{0}([t^*,S]: H^{\frac{\alpha}{2}}(\R))$.

	Let us prove the estimate \eqref{size}. By construction of $(\rho_j(t))_{j=1}^{N}$, we have that for all $t\in [t^*,S]$, there exists $(Y_j(t))_{j=1}^{N}\in \R_{L}^{N}$ such that 
	\begin{align*}
		|\rho_j(t)-Y_j(t)|\leq C\gamma,\\
		\| u(t,\cdot) - R_{Y(t)} \|_{H^{\frac{\alpha}{2}}} \leq \gamma.
	\end{align*}
	By the triangle inequality and  mean value theorem, we deduce 
	\begin{align*}
		\| \eta(t,\cdot) \|_{H^{\frac{\alpha}{2}}}&\leq \| u(t,\cdot) - R_{Y(t)} \|_{H^{\frac{\alpha}{2}}} + \| R_{Y(t)} - R \|_{H^{\frac{\alpha}{2}}}\\
		&\leq \gamma +\sum_{j=1}^{N} |\rho_j(t)-Y_j(t)| \|\partial_xQ \|_{H^{\frac{\alpha}{2}}}\leq C\gamma.
	\end{align*} 
	This finishes the proof of \eqref{size}. Note also that \eqref{vitinf2} is a direct consequence of \eqref{tr1}.

	\noindent\emph{Second step : regularity of the functions $\rho_j$.}
	Assume that the $N$ functions $\rho_j$ are $C^1([t^{*},S]:\R)$. First, we compute the equation for $\eta$ using \eqref{GC} and \eqref{defeta}
	\begin{align*}
		\partial_t\eta &=\partial_x\left( \cG(\eta) +|D|^{\alpha}R  - R^2   \right) + \displaystyle\sum_{0\leq k \leq N} \rho'_k \partial_x R_k,\\  
	\end{align*}
	where
	$$
	\cG(\eta):=|D|^{\alpha}  \eta - 2R\eta - \eta^2.
	$$
	Moreover, since
	$$
	R^2=\displaystyle\sum_{1\leq k\leq N} R^2_k + 2\displaystyle\sum_{1 \leq l<m\leq N} R_lR_m ,
	$$
	this implies by using \eqref{eqQG2} that 
	\begin{align}
		\partial_t\eta &=\partial_x\left( \cG(\eta)  - \displaystyle\sum_{1\leq k\leq N} c_k R_k- 2\displaystyle\sum_{1 \leq l<m\leq N} R_lR_m   \right) + \displaystyle\sum_{1\leq k \leq N} \rho'_k \partial_x R_k. \label{ah}
	\end{align}
	Furthermore, we obtain  differentiating in time the relation $\displaystyle\int\left(\partial_x R_j \right)\eta=0$ 
\begin{align}
	0=\frac{d}{dt} \displaystyle\int\left(\partial_x R_j \right)\eta&										    =-\rho'_j\displaystyle\int\left(\partial^2_xR_j \right)\eta
	+ \displaystyle\int\left(\partial_x R_j \right)\partial_t\eta.\label{inteta}
\end{align}
	Replacing \eqref{ah} in \eqref{inteta}, and integrating by parts, we obtain that
	\begin{align}
		0&=-\displaystyle\int\left(\partial^2_x R_j \right)\left( \cG(\eta)  - \displaystyle\sum_{1\leq k\leq N} c_kR_k - 2\displaystyle\sum_{1\leq l<m\leq N} R_lR_m   \right) + \displaystyle\sum_{1\leq k \leq N} \rho'_k\int (\partial_x R_k)(\partial_xR_j) \notag\\
		& - \rho'_j\displaystyle\int\left(\partial^2_xR_j \right)\eta. \label{2.30}
	\end{align}

	Finally, we deduce that, for all $ j\in\{ 1,...,N\}$,
	\begin{align*}
		\displaystyle\sum_{1\leq k \leq N} \rho'_k\int (\partial_x R_k)(\partial_xR_j)
		-&\rho'_j\displaystyle\int\left(\partial^2_xR_j \right)\eta
		\\
		&=\displaystyle\int\left(\partial^2_x R_j \right)\left( \cG(\eta)  - \displaystyle\sum_{1\leq k\leq N} c_kR_k- 2\displaystyle\sum_{1\leq l<m\leq N} R_lR_m   \right).
	\end{align*}
	We can rewrite this ODE system in the matrix form
	\begin{align}
		AY'=B, \label{matg}
	\end{align}
where $Y:=(\rho_j)_{j=1}^N$ and $A:=A_{0}+A_{\eta}$ where 
	\begin{align*}
		A_{\eta}:=
		\begin{pmatrix}
			-\displaystyle\int\left(\partial^2_xR_1 \right)\eta & \displaystyle\int(\partial_x R_1)(\partial_xR_2) & \ldots & \displaystyle\int(\partial_x R_1)(\partial_xR_N)\\
			\displaystyle\int(\partial_x R_2)(\partial_xR_1) & -\displaystyle\int\left(\partial^2_xR_2 \right)\eta & \ldots & \displaystyle\int(\partial_x R_2)(\partial_xR_N)\\
			\vdots & & & \vdots\\
			\displaystyle\int(\partial_x R_N)(\partial_xR_1) &\ldots & \ldots & -\displaystyle\int\left(\partial^2_xR_N \right)\eta
		\end{pmatrix},
	\end{align*}
\begin{align*}
		&A_0:=\begin{pmatrix}
			\displaystyle\int\left(\partial_xQ_{c_1} \right)^2 &0 &\ldots & 0\\
			0& \ddots & & 0\\
			\vdots & & & \vdots\\
			0&\ldots & \ldots & \displaystyle\int\left(\partial_xQ_{c_N} \right)^2
		\end{pmatrix},\\
		&B:=\left(\displaystyle\int\left(\partial^2_x R_j \right)\left( \cG(\eta)  - \displaystyle\sum_{1\leq k\leq N} c_k R_k- 2 \displaystyle\sum_{1\leq l < m\leq N} R_lR_m   \right) \right)_{1\leq j\leq N}.
	\end{align*}
	In order to prove that $A$ is invertible, it suffices to prove that $A_0$ is invertible and $\lVert A_{\eta} \rVert_{\infty}$ can be taken small enough.
	By using the Cauchy-Schwarz inequality and \eqref{size}, we have 
	$$
	\bigg|\int\left(\partial^2_xR_j \right)\eta\bigg|\leq C\gamma.
	$$  
	By  \eqref{estnonmod2}, we obtain 
	$$
	\bigg|\int(\partial_x R_k)(\partial_xR_j)\bigg|\leq \frac{C}{1+L^{2+\alpha}}, \quad k\neq j.
	$$
	Taking $L>L_{5}$ big enough, and $\gamma<\gamma_{3}$ small enough, the matrix $A$ is invertible and we can rewrite  \eqref{matg} as 
	$$
	Y'=A^{-1}B.
	$$
	Now, we have to prove that $A^{-1}B$ is globally Lipschitz. Let us begin with the term  $\displaystyle\int\left(\partial_x^2R_1\right)\eta $. Let $(\rho_j)_{j=1}^N,(\widetilde{\rho}_{j})_{j=1}^N\in\mathbb{R}^N$, by the Plancherel identity and the Cauchy-Schwarz inequality
		\begin{align*}
		&\bigg| \int \partial_x^2 \left(Q_{c_1}(x- \rho_1)\right)\left(u(t,x) - \sum_{j=1}^{N} Q_{c_j}\left(x-\rho_j \right) \right) - \partial_x^2 \left(Q_{c_1}(x-\widetilde{\rho}_1)\right)\left(u(t,x) - \sum_{j=1}^{N} Q_{c_j}\left(x-\widetilde{\rho}_j \right) \right)   dx\bigg|\\
		&\leq \int |\xi|^2|\widehat{Q_{c_1}}|  |u|  \bigg|e^{i\xi \rho_1}-e^{i\xi\widetilde{\rho}_1} \bigg|d\xi + \sum_{j=1}^{N}\int |\xi|^2|\widehat{Q_{c_1}}| |\widehat{Q_{c_j}}| \bigg|e^{i\xi (\rho_1+\rho_j)}-e^{i\xi(\widetilde{\rho}_1+\widehat{\rho}_j)} \bigg|d\xi\\
		&\leq|\rho_1-\widetilde{\rho}_1|\int |\xi|^3|\widehat{Q_{c_1}}|  |\widehat{u}|d\xi  + \left( |\rho_1-\widetilde{\rho}_1| + \sum_{j=1}^{N} |\rho_j-\widetilde{\rho}_j|\right)\sum_{j=1}^{N} \int |\xi|^3 |\widehat{Q_{c_1}}| |\widehat{Q_{c_j}}| d\xi \\
		&\leq C \sum_{j=1}^{N}|\rho_j-\widetilde{\rho}_j| \lVert \partial_x^3Q_{c_1}\rVert_{L^2}\left( \lVert u_0 \rVert_{L^2} + \| Q_{c_j}\|_{L^2}\right),
\end{align*}
where we have used for the last inequality that $\lVert u(t,\cdot) \rVert_{L^2} = \lVert u_0 \rVert_{L^2} $.

	Using the same argument for the other term in $A$ and $B$, we get $A^{-1}B$ is globally Lipschitz. Therefore, we obtain $N$ unique $C^1$ functions $\tilde{\rho}_j: [t^*,S]\longrightarrow \R$ satisfying \eqref{inteta} with $\tilde{\rho}_j(S)=\rho_j(S)$ as initial condition, where $(\rho_j)_{j=1}^N$ is given by the first step.
	Since \eqref{Ortg} is verified at time $S$ with $\rho_j(S)$, we deduce that for all $t\in [t^*,S]$,
	$$
	\displaystyle\int(\partial_x Q_{c_j}(x-\tilde{\rho}_j(t)))(u-Q_{c_j}(x-\tilde{\rho}_j(t)))=0.
	$$
	By the uniqueness statement of the first step, we conclude that the $N$ functions $\rho_j$, constructed in the first step, are $C^1$ functions. 
	This concludes the proof of Proposition \ref{modg} by taking $\gamma<\gamma_1=\min(\gamma_2,\gamma_3)$ and  $L>L_1=\max(L_2,L_3,L_4,L_5)$ .
\end{proof}



\subsection{Bootstrap setting} \label{sec:boot}
Let $\left(S_n\right)_{n=0}^{+\infty}$ be a non-decreasing sequence of time going to infinity, with $S_n>T_0$, for $T_0>1$ large enough to be chosen later. We define by $u_n$ the solution of \eqref{GC} satisfying  
\begin{align}
	u_n(S_n,\cdot)=\displaystyle\sum_{j=1}^{N} Q_{c_j}(\cdot-\rho_{j,n}^{\text{in}}), \label{initdata}
\end{align}
with \begin{align}
	\rho_{j,n}^{\text{in}}\in I_{j,n}:=[c_jS_n-S_n^{1-\frac{\alpha}{4}},c_jS_n+S_n^{1-\frac{\alpha}{4}}] ,\quad  \text{ for all } j\in\{1,\cdots,N\},\label{initintervalle}
	\end{align}
 to be fixed later.

For $t\leq S_n$, as long as the solution $u_n$ exists and satisfies \eqref{alpha} for suitable $0<\gamma_0<\gamma_1$ and $L_0>L_1$ (which will also be fixed later), we consider the $C^{1}$ functions $(\rho_{j,n})_{j=1}^{N}$ provided by Proposition \ref{modg} and satisfying \eqref{defeta}-\eqref{vitinf2}. At $S_n$, the decomposition satisfies 
\begin{align}
	\eta(S_n)=0,\quad \rho_{j,n}(S_n)=\rho_{j,n}^{in}, \quad j=1,\cdots,N.\label{initcond}
\end{align}

We introduce the bootstrap estimates at $t\leq S_n$, assuming that $u_n$ satisfies \eqref{alpha}:
\begin{align}
	\| \eta(t,x) \|_{H^{\frac{\alpha}{2}}}< \gamma_0 \label{bstr5},
	\end{align}
\begin{align}
	\sup_{j\in\{1,\cdots,N\}}|\rho_{j,n}(t)-c_j t| \le t^{1-\frac{\alpha}{4}} \label{bstr6},
	\end{align}
with $\eta$ defined in \eqref{defeta}.

For $T_0>1$, to be chosen later, we  define 
\begin{align*}
	t^*_n=\inf\left\{ T_0<\tilde{t}\leq S_n:   \exists \eps_n>0 \text{ such that } \eqref{bstr5}-\eqref{bstr6} \text{ holds for all } t\in[\tilde{t},S_n+\eps_n]  \right\}.
	\end{align*}

Note by \eqref{initcond} and by continuity that there exists $\eps_n>0$ such that \eqref{bstr5} holds on $[S_n-\eps_n,S_n+\eps_n]$. Moreover, if $\rho_{j,n}\in \mathring{I}_{j,n}$ for all $j\in\{1,\cdots,N\}$, then by possibly taking $\eps_n$ smaller, \eqref{bstr6} holds also on $[S_n-\eps,S_n+\eps]$ so that $t^{*}_n$ is well-defined. In the case where $\rho_{j_0,n}\in \partial I_{j_0,n}$ for some $j_0\in\{1,\cdots,N\}$, it follows from the transversality property (see \eqref{topo1} below) that $t^{*}_n=S_n$.

The  main result of this section states that there exists at least one choice of $(\rho_{j,n}^{in})_{j=1}^{N}\sim(c_j S_n)_{j=1}^{N}$ such that $t^{*}_n=T_0$. In other words, the bootstrap estimates \eqref{bstr5}-\eqref{bstr6} are valid up to a time $T_0$ independent of $n$.  

\begin{proposition}\label{unifG}
	Let $\alpha\in\left(\frac{1}{2},2\right)$. There exist $T_0>1$, $C_0>1$, $\gamma_0>0$ satisfying $\frac{C_0}{ T_0^{\frac{\alpha}{2}}}<\frac{\gamma_0}{2}$ and $L_0:=\frac{\beta T_0}{2}>L_1$ such that the following is true. For all $n\in\N$, there exists $(\rho_{j,n}^{in})_{j=1}^{N}$$\in I_{j,n}$, with $I_{j,n}$ defined in \eqref{initintervalle}, satisfying 
	\begin{align}
		|\rho_{j,n}^{in}-c_j S_n |\leq S_n^{1-\frac{\alpha}{4}}, \quad j\in\{1,\cdots,N\},\label{initbstr}
		\end{align}
	and $t_n^{*}=T_0$
\end{proposition}

Subsections \ref{sec:mod:est} and \ref{sec:prop6} are dedicated to the proof of Proposition \ref{unifG}. In every step of the proof, $T_0$ will be taken large enough and $\gamma_0>0$ small enough independently of $n$.

\subsection{Modulation estimates} \label{sec:mod:est}

\begin{proposition} \label{prop:mod:est}
	For all $\rho_{j,n}^{\text{in}}\in I_{j,n}$, and for all $t\in[t^{*}_{n},S_n]$, we have 
\begin{align}
	\inf_{j\in\{1,\cdots,N\}}|\rho_{j+1,n}(t)-\rho_{j,n}(t)|\geq \beta t, \label{vitinf}
	\end{align}
	\begin{align}
	|\rho'_{j,n}(t)-c_j|\leq C_{*}\left(\frac{1}{ \left(\beta t\right)^{\alpha + 1}}+ \left(\int \frac{1}{(1+|x-\rho_{j,n}(t)|)^{1+\alpha}}\eta^2\right)^{\frac{1}{2}} +  \|\eta \|^2_{L^2}  \right).  \label{VitTransg}
\end{align}
	\end{proposition}

\begin{proof}
By the triangle inequality and \eqref{bstr6}, for $t$ large enough, we deduce that 
\begin{align*}
	|\rho_{j+1,n}(t)-\rho_{j,n}(t)|\geq  (c_{j+1}-c_j)t - |\rho_{j+1,n}-c_{j+1}t| - |\rho_{j,n}-c_{j}t|\geq 2\beta t - 2t^{1-\frac{\alpha}{4}}\geq \beta t.
	\end{align*}

Now, we prove \eqref{VitTransg}. We deduce from \eqref{2.30} that
\begin{align}
	\left(\rho'_{j,n}-c_j\right)\int\left(\partial_xR_j\right) ^2
	=&\displaystyle\int \left(\partial^2_x R_j \right)\left( \cG(\eta)  - \displaystyle\sum_{1\leq k\neq j\leq N} c_kR_k- 2\displaystyle\sum_{1\leq l<m\leq N} R_lR_m   \right)\label{*}\\
	&-\displaystyle\sum_{1\leq k\neq j \leq N} \rho'_{k,n}\int (\partial_x R_k)(\partial_xR_j)
	+\rho'_{j,n}\displaystyle\int \left(\partial^2_xR_j \right)\eta\notag
\end{align}
for all $j\in\{1,...,N \}$. By using the fact the operator $|D|^{\alpha}$ is self adjoint and the Cauchy-Schwarz inequality, we deduce 
\begin{align*}
	\bigg|\displaystyle\int\left(\partial^2_x R_j \right)\cG(\eta)\bigg| + \bigg| \int\left(\partial^2_xR_j \right)\eta \bigg| \leq C \left(\|\eta \|_{L^2} + \| \eta \|_{L^2}^2\right).
\end{align*}
Moreover, by \eqref{estnonmod5}, \eqref{estnonmod6}, \eqref{estnonmod7}, and \eqref{vitinf} we get 
\begin{align*}
	\sum_{k\neq j} \left(c_k	\bigg|\int\left(\partial^2_x R_j \right)R_k\bigg| + \rho_j'\bigg|\int \partial_x R_k\partial_xR_j\bigg|\right) + 2\sum_{1\leq l<m\leq N} 	\bigg|\int\left(\partial^2_x R_j \right) R_lR_m   \bigg|  	\leq  \frac{C}{\left(\beta t\right)^{1 +\alpha}}.  
\end{align*} 
Gathering the two former estimates, we deduce that for all $j\in\{1,\cdots,N\}$,
\begin{align}
	|\rho'_{j,n}-c_j|\int\left(\partial_xR_j\right) ^2
	\leq& C\left(\gamma_0 ( 1 +|\rho'_{j,n}| )+\sum_{1\leq k\neq j \leq N}\frac{1+|\rho'_{k,n}|}{(\beta T_0)^{ 1 + \alpha }}  \right), \label{**}
\end{align}
which implies after by summing over $j$,
$$
\sum_{k=1}^{N}|\rho_{j,n}'|\leq C_{*}.
$$
Finally, by reinjecting the former estimate in \eqref{*}, we conclude that for all $j\in\{1,\cdots,N\}$
\begin{align*}
	|\rho'_{j,n}(t)-c_j|\leq C_{*}\left(\frac{1}{\left( \beta t\right)^{\alpha + 1}}+ \left(\int \frac{1}{(1+|x-\rho_{j,n}(t)|)^{1+\alpha}}\eta^2\right)^{\frac{1}{2}} +  \|\eta \|^2_{L^2}  \right),
\end{align*}
which yields \eqref{VitTransg}.	
\end{proof}

\subsection{Proof of Proposition \ref{unifG}} \label{sec:prop6}
 The proof of Proposition \ref{unifG} relies on the following result which will be proved in Section 3.
\begin{proposition}[Bootstrap estimate] \label{s=1/2}
	Let $\alpha\in\left(\frac{1}{2},2\right)$. There exist $C_0>1, 0<\gamma_2<\gamma_1$ and $T_2>T_1$ such that for all $\rho_{j,n}^{\text{in}}\in I_{j,n}$, for all  $t\in[t_n^{*},S_n]$, for all $0<\gamma_0<\gamma_2$, and $T_0>T_2$  
	\begin{align}
		\| u_n(t,\cdot) - \sum_{j=1}^{N} Q_{c_j}(\cdot -\rho_{j,n}(t)) \|_{H^{\frac{\alpha}{2}}}  \leq \frac{C_0}{ t^{\frac{\alpha}{2}}}\label{bstr}.
	\end{align}
\end{proposition}

\begin{proof}[Proof of Proposition \ref{unifG} assuming Proposition \ref{s=1/2}]
Let $0<\gamma_0<\gamma_2$ and $T_0>T_2$ such that $\frac{C_0}{ T_0^{\frac{\alpha}{2}}}<\frac{\gamma_0}{2}$. First, we show that $u_n$ satisfies \eqref{alpha} with $L_0=\frac{\beta T_0}{2}$ and that \eqref{bstr5} is strictly improved on $[t^{*}_n,S_n]$. Indeed, it follows from \eqref{bstr} that  
$$
\lVert u_n(t,\cdot)- \sum_{j=1}^{N}Q_{c_j}(\cdot-\rho_{j,n}(t))\rVert_{H^{\frac{\alpha}{2}}}\leq \displaystyle\frac{C_0}{ T_0^{\frac{\alpha}2} } <\frac{\gamma_0}{2}.
$$
Moreover, \eqref{vitinf} implies that
$$
\inf_{j\in\{1,\cdots,N-1\}\}}|\rho_{j+1,n}(t^*_n)-\rho_{j,n}(t^*_n)|\geq \beta T_0=2L_0.
$$

Now, we prove that there exists  $\rho_n^{\text{in}}=(\rho_{j,n}^{\text{in}})_{j=1}^{N}\in\R^N$, satisfying \eqref{initbstr}, such that $t^*_n=T_0$. Assume by contradiction that for all choices $\rho_n^{\text{in}}$ satisfying \eqref{initbstr}, the associated maximal time  $t^{*}_n(\rho_n^{\text{in}})>T_0$.

First, we remark that $\rho_{j,n}^{\text{in}}=c_jS_n + \lambda_{j,n} S_n^{1-\frac{\alpha}{4}}$ for a unique $\lambda_{j,n}\in[-1,1]$ and we denote  $t^{*}(\lambda_n):=t^*_n(\rho_n^{\text{in}})$ (which will also be denoted $t^{*}$ when there is no risk of confusion), with $\lambda_n=(\lambda_{j,n})_{j=1}^{N}$. By definition of $t^{*}$ and the fact that \eqref{alpha} and \eqref{bstr5} are strictly improved on $[t^{*},S_n]$, we have that
\begin{align}
|\rho_{j_0,n}(t^{*})-c_{j_0}t^{*}|=(t^{*})^{ 1-\frac{\alpha}{4}},	\label{sature}
\end{align}
for at least one $j_0\in \{1,\cdots,N\}$. Then, we define 
\begin{align*}
\Phi: [-1,1]^{N}&\to \quad \partial[-1,1]^{N}\\	
       \lambda\quad\quad &\mapsto \quad \left((\rho_{j,n}(t^{*}(\lambda))-c_jt^{*}(\lambda))(t^{*})^{ \frac{\alpha}{4}-1}(\lambda)\right)_{j=1}^{N}.
\end{align*}
We set 
\begin{align*}
f:\R&\to \R^{+}\\	
   s&\mapsto \sup_{j\in\{1,\cdots,N\}}\left((\rho_{j,n}(s)-c_js)s^{ \frac{\alpha}{4}-1}\right)^2.
\end{align*}
We claim that  if for $s_0 \in [T_0,S_n]$, \eqref{sature} is verified in $s_0$ for at least one $j\in\{1,\cdots,N\}$, then 
\begin{align}
f \text{ is a decreasing function in a neighborhood of } s_0, \label{topo1}
\end{align}
and 
\begin{align}
\Phi\in C^{0}([-1,1]&^{N},\partial[-1,1]^{N}) \label{topo2}.	
\end{align}
Let us assume \eqref{topo1} and \eqref{topo2} and finish the proof of Proposition \ref{unifG}. 
For any $\lambda\in \partial [-1,1]^{N}$, we have that 
$$
|\rho_{j_0,n}^{\text{in}}-c_{j_0}S_n|=S_n^{1-\frac{\alpha}{4}}, \quad \text{ for at least one } j_0\in\{1,\cdots,N\},
$$
which implies by \eqref{topo1} that $t^{*}=S_n$. Hence, we deduce that $\Phi_{|\partial[-1,1]^N}=\text{Id}$. However, it is a well-known topological result that no such continuous function $\Phi: [-1,1]^N \to \partial [-1,1]^N$ can exist (see Theorem 1.4, Chapter 3 in \cite{hirsch2012differential}). This concludes the proof of Proposition \ref{unifG}.

Now, we prove \eqref{topo1} and \eqref{topo2}.  Let $f_j(s)=\left((\rho_{j,n}(s)-c_js)s^{ \frac{\alpha}{4}-1}\right)^2$ for $j\in\{1,\cdots,N\}$. Let $s\in \R$. Note that for all $j\in \{1,\cdots,N\}$ the functions $f_j$ are continuously derivable. Then, to prove \eqref{topo1}, it is enough to show that for a time $s_0$  verifying \eqref{sature}, for all $j\in\{1,\cdots,N\}$ such that $f_{j}(s_0)= f(s_0)$ , we have that $f_{j}'(s_0)<0$.

By direct computations, we have that
\begin{align*}
f_{j_0}'(s)=&	2\left((\rho_{j_0,n}(s)-c_{j_0}s)s^{ \frac{\alpha}{4}-1}\right)\left((\rho'_{j_0,n}(s)-c_{j_0})s^{ \frac{\alpha}{4}-1} + \left(\frac{\alpha}{4}-1\right) s^{ \frac{\alpha}{4}-2}(\rho_{j_0,n}(s)-c_{j_0}s)  \right)\\
      =&2\left(\frac{\alpha}{4}-1\right) (\rho_{j_0,n}(s)-c_{j_0}s)^2s^{\frac{\alpha}{2}-3} +2 (\rho_{j_0,n}(s)-c_{j_0}s)(\rho'_{j_0,n}(s)-c_{j_0})s^{ \frac{\alpha}{2}-2}.
\end{align*}
Moreover, inserting \eqref{bstr} in \eqref{VitTransg}, we get for all $j\in \{1,\cdots,N\}$
\begin{align}
	|\rho_{j,n}'(s)-c_j|\leq \frac{C_{*}}{ s^{\frac{\alpha}{2}}},\label{rhoprime}
	\end{align}
which implies, combined with \eqref{sature} in $s_0$ that
\begin{align*}
f_{j_0}'(s_0)\leq 2\left(\frac{\alpha}{4} - 1 \right)s_0^{-1} + 2s_0^{\frac{\alpha}{4}-1}|\rho'_{j_0,n}(s_0)-c_{j_0}|	\leq 2\left(\frac{\alpha}{4} - 1 \right) s_0^{-1} + 2C_{*} s_0^{-\frac{\alpha}{4}-1}.
\end{align*}
Since $\alpha<2$, for $T_0$ large enough, we conclude that 
$$
f_{j_0}'(s_0)<0.
$$
The same computations yield 
$$
f_{j_1}'(s_0)<0.
$$
Then, we conclude that $f'(s_0^{+})<0$ and $f'(s_0^{-})<0$, in other words $f$ is a decreasing function at $s_0$.
Note that for $s_0=S_n$ and $\lambda\in \partial[-1,1]^{N}$, we get that $f$ is a decreasing function at $S_n$.

To show \eqref{topo2}, we prove that the map $:\lambda\in[-1,1]^{N}\mapsto t^{*}(\lambda)$ is continuous. The continuity of $t^{*}(\lambda)$ follows from the transversality property \eqref{topo1}. Indeed, by \eqref{topo1}, for all $\eps>0$ there exists $\delta>0$ such that $f(t^{*}(\lambda)-\eps)>f(t^{*}(\lambda))+\delta=1+\delta$ and for all $t\in [t^{*}(\lambda)+\eps,S_n]$ (possibly empty), $f(t)<1-\delta$. 

Note that $f$ is depending on the parameter $\lambda$ since $\rho_{j,n}(t)=\rho_{j,n}(u(t,\cdot))$. Moreover the functions $\rho_{j,n}$ are globally defined. 

Then, by the continuity of the flow , there exists $\eta>0$ such that for all $|\lambda-\bar{\lambda}|<\eta$, with $\bar{\lambda}\in[-1,1]^N$, the corresponding $\bar{f}$ satisfies $|\bar{f}(s)-f(s)|<\frac{\delta}{2}$ for $s\in[t^{*}(\lambda)-\eps,S_n]$. We deduce that 
\begin{align*}
	\bar{f}(t^{*}(\lambda)+\eps)<|\bar{f}(t^{*}(\lambda)+\eps) -f(t^{*}(\lambda)+\eps)| +f(t^{*}(\lambda)+\eps)<1-\frac{\delta}{2}=\bar{f}(t^{*}(\bar{\lambda}))-\frac{\delta}{2}.
	\end{align*}
In other words, $t^{*}(\bar{\lambda})<t^{*}(\lambda)+\eps$. Furthermore,
\begin{align*}
\bar{f}(t^{*}(\lambda)-\eps)>f(t^{*}(\lambda)-\eps)-|\bar{f}(t^{*}(\lambda)-\eps)-f(t^{*}(\lambda)-\eps)|>1+\frac{\delta}{2}	
\end{align*}
Then, $t^{*}(\lambda)-\eps<t^{*}(\bar{\lambda})$. This finishes the proof of \eqref{topo2}.
\end{proof}

\begin{remark}
The choice of the exponent $1-\frac{\alpha}{4}$ used in \eqref{bstr6} is related to the algebraic decay $s^{-\frac{\alpha}{2}}$ obtained on $\rho_{j,n}'$ in \eqref{rhoprime}. If one could integrate  directly the estimate \eqref{rhoprime}, one would expect a bound of the form $|\rho_{j,n}(t)-c_jt|\leq C t^{1-\frac{\alpha}{2}}$. However, since $\frac{\alpha}{2}<1$ a direct integration of the quantity \eqref{rhoprime} is not possible to close the bootstrap estimate. For this reason, we make use of a topological argument which allows us to adjust carefully the initial conditions $\rho_{j,n}^ {in}$ to integrate \eqref{rhoprime}. The key point is the transversality condition \eqref{topo1}. The sign of the derivative of $f$ is obtained by comparing the quantity $\rho_{j,n}'-c_j$ and $\rho_{j,n}-c_jt$. Therefore, one could get an estimate of the form $|\rho_{j,n}-c_jt| \le t^{1-\frac{\alpha}{2}+\eps}$, with $\eps\in\left(0,\frac{\alpha}{2}\right)$. By convenience we choose $\eps=\frac{\alpha}{4}$ which yields the bound in \eqref{bstr6}.  Note however that this choice does not affect the bound in \eqref{theo:est1}.    
\end{remark}

\subsection{Proof Theorem \ref{exsol} assuming Proposition \ref{unifG} }
First, we state the weak continuity property of the flow of \eqref{GC}. Relying on the well-posedness result in \cite{molinet2018well}, this result is proved in the Appendix \ref{Prweakflow} in the case $\alpha>\frac67$. It will  be admitted otherwise. 

	For sake of clarity, we recall the well-posedness result stated in \cite{molinet2018well}. To do so, we need to define $X^{s,b}_T$ the Bourgain space restricted on a time interval $[0,T]$. First, we introduced the Bourgain space $X^{s,b}$ as the completion of the Schwartz space $\mathcal{S}(\R^2)$ under the norm 
	\begin{align*}
		\| u \|_{X^{s,b}} :=\| \langle \xi \rangle^{s} \langle \tau -\xi|\xi|^\alpha \rangle^{b}\mathcal{F}_{t,x}(u)\|_{L^2}, 
	\end{align*} 
	where $\mathcal{F}_{t,x}$ the Fourier transform with respect to the both variable $t$ and $x$. The restriction space  $X^{s,b}_T$ is the space of function $u:\R\times(0,T)\longrightarrow\R$ satisfying:
	\begin{align*}
		\|u\|_{X^{s,b}_T} = \inf \{ \|\tilde{u}\|_{X^{s,b}} : \tilde{u}:\R\times\R\longrightarrow\R, \tilde{u}_{|\R\times(0,T)}=u\}<+\infty.
	\end{align*}
Now, we state the well-posedness theorem obtained\footnote{Note that the uniqueness statement is slightly stronger than the one in Theorem 1.2 in \cite{molinet2018well} since we only need to assume that $u \in L^ {\infty}((0,T) : H^ {\frac{\alpha}2})$ instead of $u \in C^ {0}([0,T] : H^ {\frac{\alpha}2})$.  This statement is however direct from the proof of Theorem 1.2 in \cite{molinet2018well}.} in \cite{molinet2018well}
\begin{thm}[Global well-posedness in the energy space] \label{GWP:theo}
	Let $\alpha\in(\frac{6}{7},2)$, and $T>0$. \\
	\underline{\emph{\textit{Existence:}}} For any $u_0\in H^{\frac{\alpha}{2}}(\R)$, there exists $u$ solution of \eqref{GC}, with the initial condition $u_0$, belonging to the class
	\begin{align*}
		Y_T:= C^{0}([0,T]:H^{\frac{\alpha}2}(\R))\cap X_{T}^{\frac{\alpha}2-1,1}\cap L^{2}((0,T),W^{\frac{\alpha}4-\frac12-,\infty}(\R)).
	\end{align*} 
\underline{\emph{\textit{Uniqueness:}}} The solution $u$ is unique in the following class 
	\begin{align*}
	Z_T:= L^{\infty}((0,T):H^{\frac{\alpha}2}(\R))\cap X_{T}^{\frac{\alpha}2-1,1}\cap L^{2}((0,T),W^{\frac{\alpha}4-\frac12-,\infty}(\R)).
\end{align*} 
Moreover, $u\in C^{0}_{b}\left(\R^{+} : H^{\frac{\alpha}2}(\R) \right)$ and the flow-map data solution $u_0\mapsto u$ is continuous from $H^{\frac{\alpha}2}(\R)$ into $C^{0}([0,T] : H^{\frac{\alpha}{2}}(\R))$.
\end{thm}

\begin{lemma}[Weak continuity of flow in the energy space]\label{weakflow}
Let $\alpha \in \left(\frac67,2\right)$. Suppose that $z_{0,n}\rightharpoonup z_0$ in $H^{\frac{\alpha}{2}}(\R)$. We consider solutions $z_n$ of \eqref{GC} corresponding to initial data $z_{0,n}$ and satisfying $z_n \in C^{0}(\mathbb R^+:H^{\frac{\alpha}2}(\mathbb R))$. Then,  $z_n(t)\rightharpoonup z(t)$ in $H^{\frac{\alpha}{2}}(\R)$, for all $t \ge 0$, where $z \in C^{0}(\mathbb R^+ : H^{\frac{\alpha}2}(\mathbb R))$ is the solution of \eqref{GC} emanating from $z_0$ obtained in Theorem \ref{GWP:theo}.
	\end{lemma}

By Proposition \ref{unifG}, there exist $C_{*}>0,T_0>0$ independent of $n$, $\rho_{1,n},\cdots, \rho_{N,n}\in C^{1}([T_0,S_n])$ satisfying  \eqref{Ortg}, \eqref{bstr6} and \eqref{bstr} for all $T_0\leq t\leq S_n$. Then, for all $ t\in[T_0, S_n]$,
$$
\|u_n(t,\cdot)   \|_{H^{\frac{\alpha}{2}}}\leq \|u_n(t,\cdot) - \sum_{j=1}^{N}Q_{c_j}(\cdot-\rho_{j,n}(t)) \|_{H^{\frac{\alpha}{2}}} + \|\sum_{j=1}^{N}Q_{c_j}(\cdot-\rho_{j,n}(t)) \|_{H^{\frac{\alpha}{2}}} \leq C_{*}.
$$
Thus, up to a subsequence, there exists $U_0\in H^{\frac{\alpha}{2}}(\R)$ such that 
\begin{align}
	u_{n}(T_0)\rightharpoonup U_0 \quad \text{in} \ H^{\frac{\alpha}{2}}(\mathbb R). \label{weak}
\end{align}
Now, we prove the convergence of the modulation parameters. Let $t\in [T_0,+\infty)$ and set $T$ such that $T_0<t<T<+\infty$. By \eqref{bstr6}, we find  that for all $j\in\{1,\cdots,N\}$ and $n\in\mathbb{N}$
$$
|\rho_{j,n}(t)|\leq T^{1-\frac{\alpha}{4}} + c_jT.
$$
Moreover from \eqref{rhoprime}, we see that $\rho_{j,n}'$ is uniformly bounded independently of time. Thus, by the Arzela-Ascoli theorem,  there exists $r_j(t)\in C^{0}([T_0,T])$ such that, after extracting a subsequence if necessary, we have 
\begin{align}
\rho_{j,n}(t)\rightarrow r_j(t). \label{cv:rj}
\end{align}

Let $U \in C^0([T_0,+\infty) : H^{\frac{\alpha}2}(\mathbb R))$ be the solution of \eqref{GC} satisfying $U(T_0,\cdot)=U_0$ obtained in Lemma \ref{weakflow}.
We set $R^{*}:= \displaystyle\sum_{j=1}^{N}Q_{c_j}(x-r_j(t))$ and let $t\in [T_0,\infty)$. By Lemma \ref{weakflow},  we know that 
\begin{equation} \label{weak:cv:un}
u_n(t) \rightharpoonup U(t) \quad \text{in} \ H^{\frac\alpha2}(\mathbb R), 
\end{equation} 
for all $t \ge T_0$. We deduce then from \eqref{bstr} and \eqref{cv:rj} that
\begin{align*}
	\|U(t,\cdot) - R^{*}(t,\cdot)  \|_{H^{\frac{\alpha}{2}}} &\leq \liminf_n \|u_n(t,\cdot) - \sum_{j=1}^{N}Q_{c_j}(\cdot-\rho_{j,n}(t)) \|_{H^{\frac{\alpha}{2}}}\\ &\quad +\liminf_n
	\sum_{j=1}^N\|Q_{c_j}(\cdot-\rho_{j,n}(t))-Q_{c_j}(\cdot-r_j(t)) \|_{H^{\frac{\alpha}{2}}} \leq \frac{C_{0}}{ t^{\frac\alpha2}}.
\end{align*}

By Proposition \ref{unifG}, we have $\frac{C_{0}}{ T_0^{\alpha}}\leq \frac{\gamma_0}{2}$ and $ \beta T_0>{2L_1}$. Moreover since $|\rho_{j+1,n}(t)-\rho_{j,n}(t)|\geq \beta T_0$, then $r_{j+1}(t)-r_j(t)\geq \beta T_0>2L_1$. Therefore, $U(t,.)\in\mathcal{T}_{\gamma, 2L_1}$ for all $t\in[T_0,\infty)$. By Proposition \ref{modg}, there exist N unique functions $\rho_1,\cdots,\rho_n \in $ $C^1([T_0,+\infty):\R)$ such that $(\rho_j)_{j=1}^N$ verify \eqref{Ortg}. On the other hand, since the solution $u_n$ satisfies also \eqref{Ortg} with $(\rho_{j,n})_{j=1}^{N}$, we deduce passing to the limit and using \eqref{cv:rj}-\eqref{weak:cv:un} that $r_j$ satisfies also \eqref{Ortg}. Hence, by the uniqueness statement in Proposition \ref{modg}, we see that $r_j(t)=\rho_j(t)$ for all $t\in\R$. Therefore, $R^{*}(t,x)= \displaystyle\sum_{j=1}^{N}Q_{c_j}(x-\rho_j(t))$, which concludes the proof of \eqref{theo:est1}. The first estimate in \eqref{theo:est2} follows passing to the limit in \eqref{bstr6}, while the second is derived arguing as Proposition \ref{prop:mod:est} and using \eqref{theo:est1}.

\section{Weighted estimates} \label{sec:weight}
We define $N$ functions to localize the information around each solitary waves. Let 
\begin{align}
	\phi(x)=1 - C_{\phi}\displaystyle\int_{-\infty}^{x}\frac{dy}{\langle y \rangle^{1+\alpha}}, \quad \text{ where }\quad C_{\phi}=\left(\int_{-\infty}^{+\infty}\frac{dy}{\langle y \rangle^{1+\alpha}}\right)^{-1}.\label{defphi}
\end{align}
We have $0\leq \phi \leq 1$. Using the function $\phi$, we set, for $A>1$ to be fixed later,
\begin{align}
	\phi_{j,A}(t,x)=\phi\left(\frac{x - \frac{\rho_j(t)+\rho_{j+1}(t) }{2}}{A}\right)=\phi \left(\frac{x - m_j(t)}{A}\right), \text{ for } j\in\{1,\cdots, N-1\},  \label{mj}
\end{align}
and $ \phi_{N,A}:=1$, where the $\rho_j$'s are defined in Section \ref{sec:cons} (in particular, they satisfy \eqref{vitinf}). 
The function $\phi_{j,A}$ follows the first $j$ solitary waves. Finally, for $j \in \{ 1,\cdots, N\}$, the function $\psi_{j,A}$ is localised around the $j^{th}$ solitary wave. Let
\begin{align}
	\psi_{1,A} = \phi_{1,A},\quad \psi_{j,A} = \phi_{j,A} - \phi_{j-1,A},\quad \psi_{N,A}= 1 - \phi_{N-1,A}.\label{psi}
\end{align}

In this section, we state some important estimates involving to the weight $\phi_{j,A}$ and its derivative $\phi_{j,A}'$. These estimates will be crucial in the proof of the monotonicity of a localised part of the mass and the energy (see Proposition \ref{mono} in Section \ref{sec:mono}).

\subsection{Weighted commutator estimates}

\begin{lemma}\label{commG}
	Let  $\alpha\in(0,2)$. In the symmetric case, there exists $C>0$ such that
	\begin{align}
		\left|\displaystyle\int \left(|D|^{\alpha} u\right) u|\phi'_{j,A}|-\int\left(|D|^{\frac{\alpha}{2}}\left(u\sqrt{|\phi'_{j,A}|} \right) \right)^2 \right|\leq \frac{C}{A^{\alpha}}\int u^2|\phi'_{j,A}|\label{sc1G},
	\end{align}
	and
	\begin{align}
		\left|\displaystyle\int \left(|D|^{\alpha} u\right) \partial_x u \phi_{j,A}  +\frac{\alpha-1}{2}\int\left(|D|^{\frac{\alpha}{2}}\left(u\sqrt{|\phi'_{j,A}|} \right) \right)^2  \right|\leq \frac{C}{A^{\alpha}}\int u^2|\phi'_{j,A}|\label{sc2G},
	\end{align}
	for any $u\in \mathcal{S}(\mathbb{R})$,  $A>1$ and $j \in \{1,\cdots,N\}$. 
	
	In the non-symmetric case, there exists $C>0$ such that
	\begin{align}  \label{nsc1G}
		\left| \displaystyle\int \left( \left( |D|^{\alpha} u \right)v - \left( |D|^{\alpha} v \right)u \right) |\phi'_{j,A}|  \right| \leq& 
		\begin{cases} \displaystyle\frac{C}{A^{\alpha}}\int \left(u^2 + v^2 \right) |\phi'_{j,A}| , & \text{if} \  \alpha\in(0,1], \\
	\displaystyle \frac{C}{A^{\frac{\alpha}{2}}} \int \left( u^2+ v^2 + \left(|D|^{\frac{\alpha}{2}}u\right)^2\right) |\phi'_{j,A}|, & \text{if} \ \alpha\in(1,2) ,
	\end{cases}
		\end{align}
		and
	\begin{align} \label{nsc2G}
		\bigg| \displaystyle\int \left( \left( |D|^{\alpha} u \right)\partial_xv + \left( |D|^{\alpha} v \right)\partial_xu \right)& \phi_{j,A}  + (\alpha-1)\int |D|^{\frac{\alpha}{2}}\left(u\sqrt{|\phi'_{j,A}|} \right) |D|^{\frac{\alpha}{2}}\left(v\sqrt{|\phi'_{j,A}|} \right)   \bigg| \notag\\ 
		&\leq \begin{cases} \displaystyle\frac{C}{A^{\alpha}}\int \left(u^2 + v^2 \right) |\phi'_{j,A}| , & \text{if} \  \alpha\in(0,1], \\
	\displaystyle \frac{C}{A^{\frac{\alpha}{2}}} \int \left( u^2+ v^2 + \left(|D|^{\frac{\alpha}{2}}u\right)^2\right) |\phi'_{j,A}|, & \text{if} \ \alpha\in(1,2) ,
	\end{cases},
	\end{align}
	for any $u,v\in \mathcal{S}(\mathbb{R})$,  $A>1$ and $j \in \{1,\cdots,N\}$. 
\end{lemma}

\begin{remark}

Instead of \eqref{nsc1G}, we can obtain that  for $\alpha_1+\alpha_2=\alpha-1$, with $0\leq\alpha_1,\alpha_2\leq\alpha-1$ and $\alpha\in(1,2)$, there exists $C>0$ such that for all $u,v\in \mathcal{S}(\R)$
	\begin{align*}
		\displaystyle\bigg|\int\left( \left( |D|^{\alpha} u \right)v - \left( |D|^{\alpha} v \right)u \right) |\phi'_{j,A}| \bigg| \leq \frac{C}{A^{\frac{\alpha}{2}}} \int \left( u^2+ v^2 + \left(|D|^{\alpha_1}u\right)^2 + \left(|D|^{\alpha_2}v\right)^2 \right) |\phi'_{j,A}| .
	\end{align*}

Moreover the estimates \eqref{nsc1G} and \eqref{nsc2G} are given with $|D|^{\frac{\alpha}{2}}$ instead of $|D|^{\alpha-1}$ . This is done to simplify the computations, terms with $|D|^{\frac{\alpha}{2}}$ appear naturally in the proof of Proposition \ref{mono}. 
	
Let us explain why we choose to force a dissymmetry on the right hand side of \eqref{nsc1G} and \eqref{nsc2G}. These two estimates will be applied with the function $v=|D|^{\alpha}u$. However, the natural quantities appearing to prove Proposition \ref{mono} are $\displaystyle\int u^2|\phi_{j,A}'|$, $\displaystyle\int (|D|^{\frac{\alpha}{2}}u)^2|\phi_{j,A}'|$ and $\displaystyle\int (|D|^{\alpha}u)^2|\phi_{j,A}'|$. Therefore, to control the remainder terms in \eqref{nsc1G} and \eqref{nsc2G} we need to impose a dissymmetry to avoid an extra derivative on the function $v$.
	 
\end{remark}

The estimates \eqref{sc1G}, \eqref{sc2G} are proved in Lemmas $6$ and $7$ in \cite{kenig2011local} for $\alpha\in[1,2]$. Observe however that their proofs extend easily to the case  $\alpha\in(0,2)$. Note also that while only one side of the inequalities in \eqref{sc1G}-\eqref{sc2G} is stated in Lemmas $6$ and $7$ in \cite{kenig2011local} , both sides are actually proved. 

While the estimates  \eqref{nsc1G} and \eqref{nsc2G} seem to be new, their proofs follow the lines of the ones of Lemma 6 and 7 of \cite{kenig2011local}. For the sake of completeness, we will present them in Appendix \ref{B.2}.

\begin{lemma}\label{eqnorm}
	Let $\alpha\in(0,2]$. There exists $C>0$ such that 
	\begin{align} \label{eqnorm1}
		\bigg|\int \left( |D|^{\alpha}\left( u\sqrt{|\phi'_{j,A}|} \right) \right)^2  &  - \int   \left( |D|^{\alpha} u  \right)^2|\phi'_{j,A}|\bigg|
		\notag\\ & \leq \begin{cases} \displaystyle\frac{C}{A^{\alpha}}\int \left( u^2 +  \left( |D|^{\alpha}u \right)^2 \right) |\phi'_{j,A}| , & \text{if} \  \alpha\in(0,1], \\
	\displaystyle  \frac{C}{A^{\frac{\alpha}{2}}}\int \left( u^2 +   \left( |D|^{\frac{\alpha}{2}}u \right)^2 + \left( |D|^{\alpha}u \right)^2 \right) |\phi'_{j,A}|, & \text{if} \ \alpha\in(1,2) ,
	\end{cases}
	\end{align}
	for all $u\in \mathcal{S}(\R)$,  $A>1$ and $j \in \{1,\cdots,N\}$.  
\end{lemma}

\begin{lemma}\label{esttc}
	Let $\alpha\in(0,2)$. There exists  $C>0$ such that 
	\begin{align} \label{est:esttc}
		\bigg|\int |D|^{\alpha}\left(u\sqrt{|\phi'_{j,A}|} \right)& \left((|D|^{\alpha}u)\sqrt{|\phi'_{j,A}|} \right)   \bigg|\notag\\ \leq&  \int \left(|D|^{\alpha} u\right)^2|\phi'_{j,A}| + \frac{C}{A^{\frac{\alpha}{2}}}  \left(\int \left( u^2 + \left(|D|^{\frac{\alpha}{2}}u\right)^2+\left(|D|^{\alpha} u\right)^2 \right) |\phi'_{j,A}|\right)	,
	\end{align}
	for all $u\in \mathcal{S}(\R)$, $A>1$ and $j \in \{1,\cdots,N\}$. 
\end{lemma}

The proofs of Lemmas \ref{eqnorm} and \ref{esttc} are also given in Appendix \ref{B.2}.

\subsection{Weighted estimates for the solitary waves}


\begin{lemma}\label{est}
	Let $p,q\geq 0$. Then, we have for all $j,k\in \{1,...,N \}$ with $k\neq j$, 
	\begin{align}\displaystyle\int R_j^p R_k^q&\leq \displaystyle\frac{C}{\left(\beta t\right)^{\left( 1+ \alpha \right) \min(p,q)}}, \label{est1} \\
		\displaystyle\int \partial_xR_j^p \partial_xR^q_k&\leq \displaystyle\frac{C}{\left(\beta t\right)^{\left( 2+ \alpha \right) \min(p,q)}},\label{est2}\\
		\displaystyle\int R_k^p \psi_{j,A}^q&\leq \displaystyle\frac{C}{\left(\beta t\right)^{\min(p(1+\alpha),q\alpha)}}. \label{est3}
	\end{align}
Moreover, we have for all $j,k\in\{1,\cdots,N\}$,
	\begin{align}
		\displaystyle\int R_k^p |\phi'_{j,A}|^q&\leq \displaystyle\frac{C}{\left(\beta t\right)^{(1+\alpha)\min(p,q)}}, \label{est4}\\
		\displaystyle\int \partial_xR_k^p |\phi'_{j,A}|^q&\leq \displaystyle\frac{C}{\left(\beta t\right)^{\min(q(1+\alpha),p(2+\alpha))}}\label{est7},\\
		\displaystyle\int R_j^p \left(1-(\psi_{j,A})^{q}\right)&\leq \displaystyle\frac{C}{\left(\beta t\right)^{\min(q\alpha ,p(1+\alpha))}}, \label{est5}\\
		\displaystyle\int \partial_xR_j^p \left(1-(\psi_{j,A})^{q} \right)&\leq \displaystyle\frac{C}{\left(\beta t\right)^{\min(\alpha q,p(2+\alpha)}}. \label{est6}
			\end{align}
\end{lemma}

Lemma \ref{est} is proven arguing exactly as in the proof of Lemma  \ref{estnonmod}.

\subsection{Weighted estimates for the non-linear terms}

\begin{lemma}\label{nonlinterm}
	Let $\alpha \in (0,2)$ and let $\eta\in H^{\frac{\alpha}{2}}(\R)$ be defined in \eqref{defeta} and verify \eqref{alpha}. Then, we have
	\begin{align}
		\displaystyle\int|\eta|^3|\phi'_{j,A}|&\leq C\gamma \left[ \int u^2|\phi'_{j,A}| + \int  \left(|D|^{\frac{\alpha}{2}}\left(u \sqrt{|\phi'_{j,A}|} \right) \right)^2 \right] + \displaystyle\frac{ C}{\left(\beta t\right)^{1+\alpha}}  \label{nlt3},
	\end{align}
	and 
	\begin{align}
		\displaystyle\int|\eta|^4|\phi'_{j,A}|&\leq C\gamma^2\left[ \int u^2|\phi'_{j,A}| + \int  \left(|D|^{\frac{\alpha}{2}}\left(u \sqrt{|\phi'_{j,A}|} \right) \right)^2 \right] + \displaystyle\frac{C}{\left(\beta t\right)^{1+\alpha}}  \label{nlt4}.
	\end{align}
\end{lemma}

\begin{lemma}\label{esttcnl}
	Let $\alpha\in(0,2)$ and let $u$ verify the hypotheses of Theorem \ref{modg}. Then there exists $C>0$ such that 
	\begin{align*}
		&\bigg|\int |D|^{\frac{\alpha}{2}}\left(u\sqrt{|\phi'_{j,A}|} \right) |D|^{\frac{\alpha}{2}}\left(u^2\sqrt{|\phi'_{j,A}|} \right)  \bigg| \\ \leq&  C\left(\gamma^2 + \frac{1}{A^{\frac{\alpha}{2}}} \right) \left(\int \left( u^2 + \left(|D|^{\frac{\alpha}{2}}u\right)^2 \right) |\phi'_{j,A}|\right)
		+ \frac{1}{8} \int \left(|D|^{\alpha} u\right)^2|\phi'_{j,A}| + \frac{C}{\left(\beta t\right)^{1+\alpha} },
	\end{align*}
	for all $A>1$ and $0<\gamma<\gamma_1$.
\end{lemma}

The proofs of these lemmas are also based on pseudo-differential estimates and are given in Appendix \ref{E}.

\section{Proof of the bootstrap estimate} \label{sec:mono}

The goal of this section is to prove Proposition \ref{s=1/2}. We work in the bootstrap setting of Section \ref{sec:boot}. In particular, the solutions $u_n$ admit the decomposition of Proposition \ref{modg} on the time interval $[t_n^*,S_n]$. We also recall the definitions of the weight functions $\phi_{j,A}$ and $\psi_{j,A}$ in \eqref{mj} and \eqref{psi}. 

In every step of the proof, the values of $T_0$ and $A$ will be taken large enough independently of $n$, while the value of $\gamma$ will be chosen small enough independently of $n$. Moreover, for simplicity of notation, we drop the index $n$ of the functions $u_n$ and $(\rho_{j,n})_{j=1}^N$ and of the time $t_n^{*}$, and the index $A$ of the weight functions $\phi_{j,A}$ and $\psi_{j,A}$.

Finally, we define the part of the mass $M_j$ and of the energy $E_j$ localized around the $j^{th}$ solitary wave $R_j$ by
$$
M_j(t):=\displaystyle\int  u(t,x)^2\psi_j(t,x) dx,\quad E_j(t):= \displaystyle\int  \left(\frac{1}{2}u|D|^{\alpha} u - \frac{1}{3}u^3\right)(t,x)\psi_j(t,x) dx,
$$
so that for all $j\in\{1,\cdots,N\}$
$$
\displaystyle\sum_{k=1}^{j}M_k(t)= \int  u(t,x)^2\phi_j(t,x) dx, \quad \displaystyle\sum_{k=1}^{j}E_k(t) =\int  \left(\frac{1}{2}u|D|^{\alpha} u - \frac{1}{3}u^3\right)(t,x)\phi_j(t,x) dx.
$$
We also define 
\begin{equation} \label{def:Etilde}
\tilde{E}_k=E_k +\sigma_0  M_k.
\end{equation}
with 
\begin{equation}
	\sigma_0:=\min_{j\in\{1,\cdots,N-1\}}\left(\displaystyle \frac{c_j}{4}, \frac{c_N}{4} , \frac{c_j c_{j+1}}{4\left(c_j + c_{j+1}\right)}\right). \label{defsigma}
\end{equation}

\subsection{Monotonicity}

\begin{proposition}[Monotonicity]\label{mono}
	Under the bootstrap assumptions \eqref{bstr5}-\eqref{bstr6}, we have 
	\begin{align}
		\sum_{k=1}^{j}\left( M_k(S_n) - M_k(t_0)\right) \geq - \displaystyle\frac{C}{\left(\beta t_0\right)^{\alpha}}, \label{M1}
	\end{align}
	and
	\begin{align}
		\sum_{k=1}^{j}\left( \tilde{E}_k(S_n) - \tilde{E}_k(t_0) \right) \geq -\displaystyle\frac{C}{\left(\beta t_0\right)^{\alpha}}, \label{M2}
	\end{align}
	for all $j\in\{1,\cdots,N \}$, $t_0\in [t^*,S_n]$.
\end{proposition}

\begin{remark}
	From Proposition \ref{mono}, we see that $M_j(S_n)$ is almost larger than $M_j(t_0)$ for $t_0<S_n$.  In other words, when the time decreases, the portion of the mass on the left of the $(j+1)^{th}$ solitary wave also decreases.  A similar phenomenon occurs also for the energy. This can be seen as a manifestation of the dispersive character of \emph{KdV-type} equations: if a wave moves to the right, then the dispersion effect pushes some mass to the left, see Figure \ref{fig1}.Moreover, if $u$ is a solution of fKdV then $u(-t,-x)$ is also a solution. Therefore, if a wave move to the left, then the dispersion effect pushes some mass to the right.
	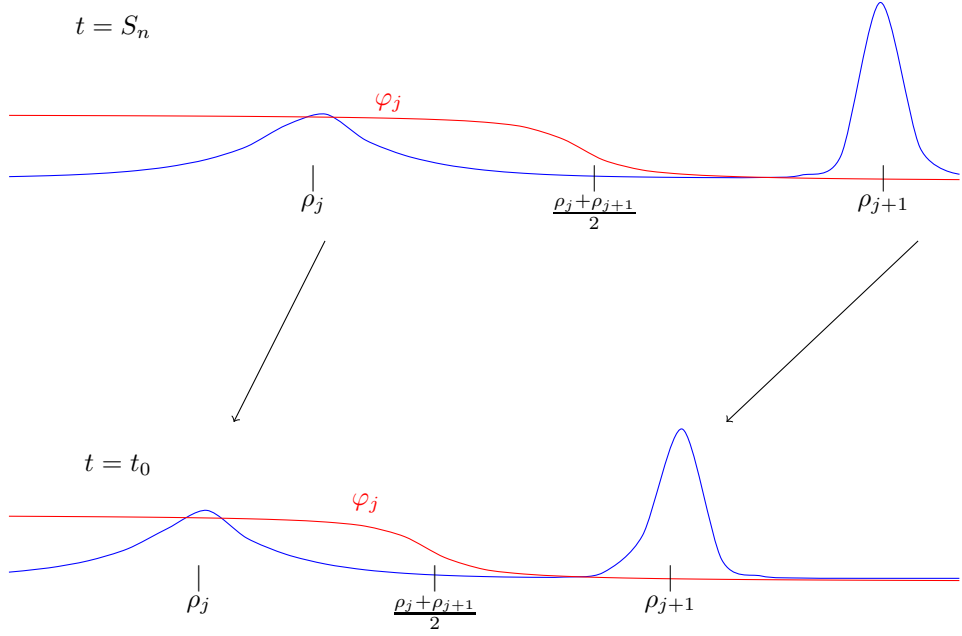
\begin{figure}[h!]
		\centering
		\begin{tikzpicture}
			\draw[scale=0.5,domain=-0:25,smooth,variable=\x,blue] plot ({\x},{2*exp(-1*abs(\x-8)/2) + 5.5*exp(-1*abs(2*(\x-23)))});
			\draw[scale=0.5,domain=-0:25,smooth,variable=\y,red] plot ({\y},{0.8-0.01*atan(\y-15)});
			\draw (4,-0.1) node[below]{$\rho_j$};
			\draw (4,0) node[]{$|$};
			\draw (7.7,-0.1) node[below]{$\frac{\rho_j+\rho_{j+1}}{2}$};
			\draw (7.7,0) node[]{$|$};
			\draw (11.5,-0.1) node[below]{$\rho_{j+1}$};
			\draw (11.5,0) node[]{$|$};
			\draw (5,1) node[red]{$\varphi_j$};
			\draw (2,2) node[left]{$t=S_n$};
		\end{tikzpicture}
		
		\hspace{23mm}
		\begin{tikzpicture}
			\draw[scale=0.6,domain=-5:25,<-] (0,-4) -- (2,0);
			\draw[scale=0.6,domain=-5:25,<-] (10.8,-4) -- (15,0);
		\end{tikzpicture}
		
		\begin{tikzpicture}
			\draw[scale=0.5,domain=-5:20,smooth,variable=\x,blue] plot ({\x},{2*exp(-1*abs(\x)/2) + 6*exp(-1*abs(2*\x-25))});
			\draw[scale=0.5,domain=-5:20,smooth,variable=\y,red] plot ({\y},{0.8-0.01*atan(\y-6)});
			\draw (0,-0.1) node[below]{$\rho_j$};
			\draw (0,0) node[]{$|$};
			\draw (3.1,-0.1) node[below]{$\frac{\rho_j+\rho_{j+1}}{2}$};
			\draw (3.1,0) node[]{$|$};
			\draw (6.2,-0.1) node[below]{$\rho_{j+1}$};
			\draw (6.2,0) node[]{$|$};
			\draw (2.2,1) node[red]{$\varphi_j$};
			\draw (-0.5,1.5) node[left]{$t=t_0$};
		\end{tikzpicture}
		\caption{Monotonicity of the mass} \label{fig1}
	\end{figure}
\end{remark}


\begin{proof}[Proof of Proposition \ref{mono}]
	

	
	We remark that for $j=N$ the inequalities \eqref{M1} and \eqref{M2} are easily verified since $M$ and $E$ are preserved by the flow of \eqref{GC}. Then, we can always assume $1 \le j \le N-1$. 
	
	First, we give the proof of \eqref{M1}. By using \eqref{GC}, integration by parts and $\phi$ is non increasing function, we get
	\begin{align*}
		\frac{1}{2}\frac{d}{dt}\left(\displaystyle\sum_{k=1}^j M_k(t)  \right)= \int \left( |D|^{\alpha}u\left(-\partial_xu \phi_j + u |\phi'_j| \right) -\frac{u^3}{3}|\phi_j'| + \frac{m_j'}{2}u^2|\phi'_j| \right).
	\end{align*}
	Then, we deduce from \eqref{sc1G}, \eqref{sc2G} that
	\begin{align}
		\int \left( |D|^{\alpha}u\left(-\partial_xu \phi_j + u |\phi'_j| \right) \right) \geq-\frac{C}{A^{\alpha}}\int u^2|\phi'_j| + \frac{\alpha+1}{2}\int  \left(|D|^{\frac{\alpha}{2}}\left(u \sqrt{|\phi'_j|} \right) \right)^2 . \label{mass}
	\end{align}
	Observe from \eqref{VitTransg} that $m_j' \ge \frac{c_j+c_{j+1}}4$. Thus,
	\begin{align*}
		\frac{1}{2}\frac{d}{dt}\left(\displaystyle\sum_{k=1}^j M_k(t)  \right)\geq -\int \frac{u^3}{3}|\phi'_j| +\frac{c_j+c_{j+1}}8\int u^2 |\phi_j'|  + \frac{\alpha+1}{2}\int  \left(|D|^{\frac{\alpha}{2}}\left(u \sqrt{|\phi'_j|} \right) \right)^2.
	\end{align*}
	
	\noindent Now, we estimate the nonlinear term. With the notation of Proposition \ref{modg}, we have 
	$$
	|u|^3\leq C \left( \displaystyle\sum_{k=1}^NR_k^3 + |\eta|^3 \right).
	$$
	Therefore, by \eqref{est4} 
	$$
	\displaystyle\sum_{k=1}^N\int |R_k^3||\phi'_j|\leq \displaystyle\frac{C}{\left(\beta t\right)^{\alpha +1}},
	$$
	and by \eqref{nlt3}
	\begin{align*}
		\displaystyle\int|\eta|^3|\phi'_j|&\leq C\gamma\left[ \int u^2|\phi'_j| + \int  \left(|D|^{\frac{\alpha}{2}}\left(u \sqrt{|\phi'_j|} \right) \right)^2 \right] + \displaystyle\frac{C}{\left(\beta t\right)^{\alpha +1}}.
	\end{align*}
	Hence, we can conclude that
	\begin{align*}
		\frac{d}{dt} \left(\displaystyle\sum_{k=1}^j M_k(t)  \right)
		&\geq - \displaystyle\frac{C}{\left(\beta t\right)^{1+\alpha}}.
	\end{align*}
	Thus, we have by integrating between $t_0$ and $S_n$ 
	\begin{align*}
		\displaystyle\sum_{k=1}^j M_k(S_n) - \displaystyle\sum_{k=1}^j M_k(t_0) 
		&\geq - \displaystyle\int_{t_0}^{S_n}\frac{C}{\left(\beta t\right)^{1+\alpha}}dt\geq - \displaystyle\frac{C}{\left(\beta t_0\right)^{\alpha}} ,
	\end{align*}
which proves \eqref{M1}.

	
	 Let us prove \eqref{M2}. We differentiate $E_j$ with respect to time to find that 
	\begin{align*}
		\frac{d}{dt}\left(\displaystyle\sum_{k=1}^j E_k(t)  \right) =&   \displaystyle\int  \left[\left(\frac{1}{2}\partial_tu |D|^{\alpha} u + \frac{1}{2}u|D|^{\alpha}\partial_t u\right) - \partial_tu u^2\right]\phi_j   \\
		&+ m'_j  \displaystyle\int  \left(\frac{1}{2}u|D|^{\alpha} u -   \frac{1}{3}u^3 \right)|\phi'_j|   \\
		=& I_1 + m'_j I_2.							
	\end{align*}
	Using \eqref{GC}, we obtain for $I_1$ that
	\begin{align*}
		I_1=&   \displaystyle \frac{1}{2} \int  \left(|D|^{\alpha} \partial_{x} u  - \partial_{x}(u^2)\right) \left(|D|^{\alpha} u\right)\phi_j + \frac{1}{2}\int u |D|^{\alpha}\left(|D|^{\alpha} \partial_{x} u - \partial_{x}(u^2)\right)\phi_j \\ 
		&- \int \left(|D|^{\alpha} \partial_{x} u - \partial_{x}(u^2)\right) u^2 \phi_j \\
		=& I_{1,1} + I_{1,2} + I_{1,3}.\\
	\end{align*}
	First we compute $I_{1,1}$ by integrating by parts
	\begin{align}
		I_{1,1}=&\frac{1}{4} \displaystyle\int  \left(|D|^{\alpha}u\right)^2|\phi'_j| - \frac{1}{2} \displaystyle\int \partial_{x}(u^2)\left(|D|^{\alpha}u\right)\phi_j ,\label{I11} 
	\end{align}
	since the functions $\phi_j$ are non increasing. Now we decompose $I_{1,2}$ as 
	\begin{align*}
		I_{1,2}=&\frac{1}{2}\displaystyle\int  u\left(|D|^{2\alpha} \partial_{x} u\right) \phi_j - \frac{1}{2} \int u\left(|D|^{\alpha}\partial_{x}(u^2) \right) \phi_j =I_{1,2,1}  + I_{1,2,2}.
	\end{align*}
	First we deal with $I_{1,2,1}$. By using integration by parts we get 
	\begin{align*}
		I_{1,2,1}=&-\frac{1}{2}\displaystyle\int  \partial_{x}u\left(|D|^{2\alpha}  u\right) \phi_j  + \frac{1}{2}\displaystyle\int  u\left(|D|^{2\alpha} u \right)|\phi'_j|  = I_{1,2,1}^1 +I_{1,2,1}^2.
	\end{align*}
	On the one hand, using the estimate \eqref{nsc2G} with $v=|D|^{\alpha}u$ and integration by parts for the last integral, we get
	\begin{align*}
		I_{1,2,1}^1 &=   -\frac{1}{2}\displaystyle\int  \left(\partial_{x}u|D|^{2\alpha}  u  + \left(|D|^{\alpha}u\right) \partial_x|D|^{\alpha}u  \right) \phi_j  - \frac{\alpha-1}{2}\int |D|^{\frac{\alpha}{2}}\left(u\sqrt{|\phi'_j|} \right) |D|^{\frac{\alpha}{2}}\left(\left(|D|^{\alpha}u\right)\sqrt{|\phi'_j|} \right)   \\
		& + \frac{\alpha-1}{2}\int |D|^{\frac{\alpha}{2}}\left(u\sqrt{|\phi'_j|} \right) |D|^{\frac{\alpha}{2}}\left(\left(|D|^{\alpha}u\right)\sqrt{|\phi'_j|} \right)  + \frac{1}{2}\displaystyle\int  \left(\partial_{x}|D|^{\alpha}u\right)\left(|D|^{\alpha}  u\right) \phi_j  \\
		&\geq \frac{1}{4}\displaystyle\int  \left(|D|^{\alpha}u\right)^2|\phi'_j|dx + \frac{\alpha-1}{2}\int |D|^{\frac{\alpha}{2}}\left(u\sqrt{|\phi'_j|} \right) |D|^{\frac{\alpha}{2}}\left(\left(|D|^{\alpha}u\right)\sqrt{|\phi'_j|} \right)  \\
		& -\frac{C}{A^{\frac{\alpha}{2}}}\displaystyle\int  \left( u^2+ \left(|D|^{\frac{\alpha}{2}}u\right)^2+ \left(|D|^{\alpha}u\right)^2 \right) |\phi'_j|  .
	\end{align*}
	On the other hand, we deduce from \eqref{nsc1G} with $v=|D|^{\alpha}u$
	\begin{align*}
		I_{1,2,1}^2 &= \frac{1}{2}\displaystyle\int  \left( u|D|^{2\alpha} u -\left( |D|^{\alpha}u \right)^2 \right) |\phi'_j|  + \frac{1}{2}\displaystyle\int  \left( |D|^{\alpha}u \right)^2  |\phi'_j| \\
		&\geq \frac{1}{2}\displaystyle\int  \left( |D|^{\alpha}u \right)^2  |\phi'_j|  -\frac{C}{A^{\frac{\alpha}{2}}}\displaystyle\int  \left( u^2 + \left(|D|^{\frac{\alpha}{2}}u\right)^2 + \left(|D|^{\alpha}u\right)^2 \right) |\phi'_j|.
	\end{align*}
	Now, we deal with $I_{1,2,2}$. Using integration by parts
	\begin{align*}
		I_{1,2,2}=& \frac{1}{2}\displaystyle\int  \partial_{x}u |D|^{\alpha} (u^2) \phi_j -\frac{1}{2}\displaystyle\int  u|D|^{\alpha}(u^2) |\phi'_j|.
	\end{align*}
	Arguing similarly as for $I_{1,2,1}$, we get from \eqref{nsc1G}, \eqref{nsc2G} with $v=u^2$ that
	\begin{align*}
		I_{1,2,2}\geq & -\frac{\alpha-1}{2}\int |D|^{\frac{\alpha}{2}}\left(u\sqrt{|\phi'_j|} \right) |D|^{\frac{\alpha}{2}}\left(u^2\sqrt{|\phi'_j|} \right)   - \frac{1}{2}\int u^2 \left(|D|^{\alpha}u\right) |\phi'_j|   \\
		&- \frac{1}{2}\int \partial_{x}\left( u^2 \right)\left(|D|^{\alpha}u\right) \phi_j  - \frac{C}{A^{\frac{\alpha}{2}}}\int \left( u^2 + u^4 + \left( |D|^{\frac{\alpha}{2}} u \right)^2 \right) |\phi'_j| .
	\end{align*}
	Hence, we conclude gathering these estimates 
	\begin{align}
		I_{1,2}\geq& \frac{3}{4}\int \left( |D|^{\alpha}u \right)^2 |\phi'_j|  - \frac{1}{2}\int \partial_{x}\left( u^2 \right)|D|^{\alpha}u \phi_j  - \frac{1}{2}\int u^2 |D|^{\alpha}u |\phi'_j|  \notag \\
		&-\frac{\alpha-1}{2}\int |D|^{\frac{\alpha}{2}}\left(u\sqrt{|\phi'_j|} \right) |D|^{\frac{\alpha}{2}}\left(u^2\sqrt{|\phi'_j|} \right)   + \frac{\alpha-1}{2}\int |D|^{\frac{\alpha}{2}}\left(u\sqrt{|\phi'_j|} \right) |D|^{\frac{\alpha}{2}}\left(\left(|D|^{\alpha}u\right)\sqrt{|\phi'_j|} \right)   \notag \\
		& - \frac{C}{A^{\frac{\alpha}{2}}}\int \left( u^2 + u^4 +\left( |D|^{\frac{\alpha}{2}}u \right)^2 + \left( |D|^{\alpha}u \right)^2 \right) |\phi'_j|  .\label{I12}
	\end{align}
	Finally, we compute $I_{1,3}$ by integrating by parts
	\begin{align}
		I_{1,3}=& \displaystyle\int  \left(|D|^{\alpha} u\right) \partial_x\left(u^2\right)\phi_j  -\displaystyle\int  \left(|D|^{\alpha} u\right) u^2|\phi'_j| + \frac{1}{2}\int u^4 |\phi'_j| . \label{I13}
	\end{align}
	Therefore, combining \eqref{I11},\eqref{I12} and \eqref{I13}, we deduce that
	\begin{align*}
		I_1\geq&\int \left( |D|^{\alpha}u \right)^2 |\phi'_j| -\frac{3}{2}\int u^2\left(|D|^{\alpha}u\right) |\phi'_j| +\frac{1}{2}\int u^4 |\phi'_j|  \\
		&-\frac{\alpha-1}{2}\int |D|^{\frac{\alpha}{2}}\left(u\sqrt{|\phi'_j|} \right) |D|^{\frac{\alpha}{2}}\left(u^2\sqrt{|\phi'_j|} \right)   + \frac{\alpha-1}{2}\int |D|^{\frac{\alpha}{2}}\left(u\sqrt{|\phi'_j|} \right) |D|^{\frac{\alpha}{2}}\left(\left(|D|^{\alpha}u\right)\sqrt{|\phi'_j|} \right)  \\
		& - \frac{C}{A^{\frac{\alpha}{2}}}\int \left( u^2 + u^4 + \left( |D|^{\frac{\alpha}{2}}u \right)^2 + \left( |D|^{\alpha}u \right)^2\right) |\phi'_j|.
	\end{align*}
	By using the identity
	\begin{align*}
		&\displaystyle\int  \left(|D|^{\alpha}u\right)^2|\phi'_j| +\frac{1}{2}\int u^4 |\phi'_j|  
		-\displaystyle\frac{3}{2}\int  \left(|D|^{\alpha} u\right) u^2|\phi'_j|\\
		&=\displaystyle\int \left(\frac{1}{4}|D|^{\alpha}u- 3u^2 \right)^2|\phi'_j| +\frac{15}{16}\displaystyle\int  \left(|D|^{\alpha}u\right)^2|\phi'_j| - \displaystyle\frac{17}{2}\int u^4|\phi'_j|,
	\end{align*}
	 $\frac{\alpha-1}{2}\in[-\frac{1}{4},\frac{1}{2}]$, and Lemmas \ref{esttc}, \ref{esttcnl} we conclude that  
	\begin{align*}
		I_1  \geq&\frac{1}{4}\displaystyle\int  \left(|D|^{\alpha}u\right)^2|\phi'_j| - 9\int u^4|\phi'_j| 
		- C\left(\frac{1}{A^{\frac{\alpha}{2}}} + \gamma^2 \right)\int \left( u^2 + \left( |D|^{\frac{\alpha}{2}}u \right)^2 \right) |\phi'_j| - \frac{C}{\left( \beta t \right)^{1+\alpha}}.
	\end{align*}
	Using \eqref{sc1G} to control $I_2$, we obtain that
	\begin{align*}
		\frac{d}{dt}\left(\displaystyle\sum_{k=1}^j E_k(t)  \right)\geq& \frac{1}{4}\int \left(|D|^{\alpha}u \right)^2|\phi'_j|  -C\left(\frac{1}{A^{\frac{\alpha}{2}}} + \gamma^2 \right) \int \left( u^2 + \left( |D|^{\frac{\alpha}{2}}u \right)^2 \right) |\phi'_j| \\& - \frac{m'_j}{3}\int |u|^3|\phi'_j| -9\int u^4|\phi'_j| - \frac{C}{\left( \beta t \right)^{1+\alpha}}.
	\end{align*}
	We need to add the mass to the energy in order to control the remaining terms
	\begin{align*}
		\frac{d}{dt} \left(\displaystyle\sum_{k=1}^j E_k(t) + \sigma_0 M_k(t)  \right)\geq&
		-9\int u^4|\phi'_j| -C\left(\frac{1}{A^{\frac{\alpha}{2}}} + \gamma^2 \right)\int \left( u^2 + \left( |D|^{\frac{\alpha}{2}}u \right)^2 \right) |\phi'_j| 
		- \frac{m'_j}{3}\int |u|^3|\phi'_j| \\
		&-\sigma_0\displaystyle\int |D|^{\alpha}u\left( \partial_x u\phi_j -u|\phi'_j| \right)
		- \sigma_0 \int \frac{1}{3}|u|^3|\phi'_j|
		+\sigma_0  m'_j\int u^2 |\phi'_j|\\
		&-  \frac{C}{\left( \beta t \right)^{1+\alpha}}.
	\end{align*}
	Thus, by using \eqref{mass}, we deduce
	\begin{align*}
		\frac{d}{dt}\left(\displaystyle\sum_{k=1}^j E_k(t) + \sigma_0 M_k(t)  \right)\geq&
		-9\int u^4|\phi'_j| 
		+\left(   \sigma_0m'_j - \frac{C}{A^{\frac{\alpha}{2}}} -\gamma^2 \right)\int u^2|\phi'_j| 
		- \frac{m'_j + \sigma_0  }{3}\int u^3|\phi'_j|\\
		&+\sigma_0\frac{3-\alpha}{2}\int \left(|D|^{\frac{\alpha}{2}}\left(u\sqrt{|\phi'_j|} \right) \right)^2   -\left(\frac{C}{A^{\frac{\alpha}{2}}} + \gamma^2 \right)\int  \left( |D|^{\frac{\alpha}{2}}u \right)^2  |\phi'_j| \\&-  \frac{C}{\left( \beta t \right)^{1+\alpha}}. 
	\end{align*}
	Observe from \eqref{VitTransg} that $\frac{\sigma_0m'_j}{2}-\frac{C}{A^{\frac{\alpha}2}}-\gamma^2>0$. Thus, by \eqref{eqnorm1}, we deduce that
	\begin{align*}
		\frac{d}{dt}\left(\displaystyle\sum_{k=1}^j E_k(t) +  \sigma_0M_k(t)  \right)\geq&
		\frac{\sigma_0m'_j}{2} \int u^2|\phi'_j| 
		-9\int u^4|\phi'_j| 
		- \frac{m'_j + \sigma_0  }{3}\int u^3|\phi'_j|\\
		&+\sigma_0\frac{3-\alpha}{4}\int \left(|D|^{\frac{\alpha}{2}}\left(u\sqrt{|\phi'_j|} \right) \right)^2 - \frac{C}{\left( \beta t \right)^{1+\alpha}}. 
	\end{align*}
	From Lemma \ref{nonlinterm} , we get 
	\begin{align*}
		\frac{d}{dt} \left(\displaystyle\sum_{k=1}^j E_k(t) + \sigma_0 M_k(t)  \right)&\geq
		\left(   \frac{\sigma_0m'_j}{4}- C\gamma \right)\int u^2|\phi'_j| 
		- \displaystyle\frac{C}{\left(\beta t\right)^{1+\alpha}} \\
		&+\left(\sigma_0\frac{3-\alpha}{2} -C\gamma \right)\int \left(|D|^{\frac{\alpha}{2}}\left(u\sqrt{|\phi'_j|} \right) \right)^2 \\
		&\geq - \displaystyle\frac{C}{\left(\beta t\right)^{1+\alpha}}.
	\end{align*}
		Thus we have by integrating between $t_0$ and $S_n$ 
	\begin{align*}
		\displaystyle\sum_{k=1}^j \tilde{E}_k(S_n) - \displaystyle\sum_{k=1}^j \tilde{E}_k(t_0) 
		&\geq - \displaystyle\int_{t_0}^{S_n}\frac{C}{\left(\beta t\right)^{1+\alpha}}dt
		\geq - \displaystyle\frac{C}{\left(\beta t_0\right)^{\alpha}},
	\end{align*}
	which proves \eqref{M2}.
\end{proof}


\subsection{Mass and energy expansion}

\begin{lemma}\label{emest}
	There exist $C>0$ such that the following hold:\\
	\begin{align}
		\bigg| M_j(t) - \left[ \displaystyle\int Q^2_{c_j} + 2\int\eta(t) R_j(t) + \int \eta^2(t) \psi_j(t) \right] \bigg|\leq \displaystyle\frac{C}{\left(\beta t\right)^{\alpha}}, \label{emest1}
	\end{align}
	\begin{align}
		\bigg| E_j(t) -\left[ E(Q_{c_j}) -c_j\displaystyle\int \eta(t) R_j(t) +\right.&\left. \displaystyle\frac{1}{2}\int\left(\eta(t) |D|^{\alpha}\eta(t) - 2 R(t)\eta^2(t)\right)\psi_j(t) \right]  \bigg| \notag \\
		&\leq \displaystyle\frac{C}{\left(\beta t\right)^{\alpha}} + C\gamma \lVert \eta(t) \rVert_{H^{\frac{\alpha}{2}}}^2\label{emest2}
	\end{align}
		and
	\begin{align}
		\bigg| \left(E_j(t) + \displaystyle\frac{c_j}{2} M_j(t)\right) -\left( E(Q_{c_j})+\right.&\left. \displaystyle\frac{c_j}{2} M(Q_{c_j}) \right) -\displaystyle\frac{1}{2}H_j(t) \bigg| \notag \\
		&\leq \displaystyle\frac{C}{\left(\beta t\right)^{\alpha}} + C\gamma\lVert \eta(t) \rVert_{H^{\frac{\alpha}{2}}}^2 ,\label{emest3}
	\end{align}
	where 
	\begin{align}
		H_j(t):=H_j(\eta(t),\eta(t))=\displaystyle\int\left( \eta(t) |D|^{\alpha}\eta(t) + c_j\eta^2(t)- 2R_j(t)\eta^2(t) \right)\psi_j(t) .\label{quadratic}
	\end{align}
\end{lemma}

\begin{proof}
	Using $u=R+\eta$ in the mass, we get
	\begin{align*}
		M_j(t)=\displaystyle\int \left( R^2 + 2R\eta + \eta^2 \right)\psi_j.
	\end{align*}
	Thus by direct computations,
	\begin{align*}
		M_j(t) - \left[ \displaystyle\int Q^2_{c_j} + 2\int\eta R_j + \int \eta^2 \psi_j \right] =  \displaystyle\int \left(R^2\psi_j - Q^2_{c_j} \right)+  2\displaystyle\int \eta\left( R\psi_j - R_j \right) = I_1 + 2 I_2.
	\end{align*}
	We use the translation invariance of the $L^2$ norm of $Q_{c_j}$, and \eqref{est1}, \eqref{est3}, \eqref{est5} to deduce
	\begin{align*}
		|I_1|
		\leq  \displaystyle \sum_{(k,i)\neq (j,j)}\int R_iR_k\psi_j +\int R_j^2\left(1-\psi_j\right)  
		\leq \displaystyle\frac{C}{\left(\beta t\right)^{\alpha}}.
	\end{align*}
	By Cauchy-Schwarz inequality, \eqref{est3}, \eqref{est5}, we obtain for $I_2$
	$$
	|I_2|\leq 2\| \eta\|_{L^2}\left(\sum_{k\neq j} \| R_k \psi_j  \|_{L^2} + \|R_j(1-\psi_j) \|_{L^2}\right) \leq \displaystyle\frac{C}{\left(\beta t\right)^{\alpha}}.
	$$
	Combining these two inequalities we conclude the proof of \eqref{emest1}.
	
	To prove \eqref{emest2},  we expand $E_j$ as
	\begin{align*}
		E_j(t)
		&=\displaystyle\frac{1}{2}  \int R \left(|D|^{\alpha}R\right)\psi_j + \frac{1}{2}  \int \eta \left(|D|^{\alpha}\eta\right)\psi_j +\frac{1}{2}  \int \left(R |D|^{\alpha}\eta + \eta |D|^{\alpha}R \right) \psi_j\\
		&\quad - \displaystyle \frac{1}{3}\int R^3\psi_j -\int R^2\eta\psi_j -\int R\eta^2 - \frac{1}{3}\int \eta^3\psi_j.
	\end{align*}
	Hence,
	\begin{align*}
		\bigg| E_j(t) -&\left[ E(Q_{c_j}) -c_j\displaystyle\int \eta  R_j + \displaystyle\frac{1}{2}\int\left(\eta |D|^{\alpha}\eta - 2 R\eta^2\right)\psi_j \right]  \bigg|\\
		&\leq \bigg|\int \left(\frac{1}{2} R|D|^{\alpha}R-\frac{1}{3} R^3 \right)\psi_j- E(Q_{c_j})\bigg| + \frac{1}{3} \bigg| \int \eta^3\psi_j\bigg|\\
		&\quad +\frac{1}{2}  \bigg|  \int \left(R |D|^{\alpha}\eta + \eta |D|^{\alpha}R \right) \psi_j - 2 \displaystyle \int R^2\eta\psi_j + 2 c_j\int \eta R_j  \bigg| \\
		&= J_1 + J_2 + J_3.
	\end{align*}
	We use the translation invariance of $E(Q_{c_j})$, \eqref{est3}, \eqref{est5} and $\||D|^{\alpha} R_k \|_{L^{\infty}}\leq C$ to bound $|J_1|$ by
	\begin{align*}
		 C \left(\int \left(  |R_j \|D|^{\alpha}R_j |+ R_j^3\right)\left(1-\psi_j  \right)+\sum_{(i,k)\neq (j,j)} \bigg|\int  R_i \left(|D|^{\alpha}R_k\right)\psi_j\bigg|+ \sum_{(i,k,l) \neq (j,j,j)}\int R_iR_kR_l\psi_j\right)
		 	\end{align*}
	Thus, we get 
	\begin{equation}
	|J_1| \leq \displaystyle\frac{C}{\left(\beta t\right)^{\alpha}}.\label{J1}
	\end{equation}
	Replacing $|\phi_j'|$ by $\psi_j$ in \eqref{nlt3est}, we have
	\begin{align*}
		|J_2|&\leq C\gamma \lVert \eta \sqrt{\psi_j}  \rVert_{H^{\frac{\alpha}{2}}}^2.
	\end{align*}
	so that it follows arguing as \eqref{q2} that
	\begin{align}
		|J_2|\leq C\gamma \lVert \eta \rVert_{H^{\frac{\alpha}{2}}}^2.\label{J2}
	\end{align}
	By using \eqref{eqQG2}, we get 
	\begin{align*}
		2|J_3|&\leq  \bigg|  \int \left(R |D|^{\alpha}\eta + \eta |D|^{\alpha}R \right) \psi_j - 2 \displaystyle \int \left(|D|^{\alpha}R_j\right)\eta  \bigg| + 2\bigg| \displaystyle \int \left(R_j^2\eta -R^2\eta\psi_j \right) \bigg|\\
		&\leq \bigg|  \int R \left(|D|^{\alpha}\eta\right)\psi_j -  \displaystyle \int \left(|D|^{\alpha}R_j\right)\eta \bigg|+ \bigg|\int \eta \left(|D|^{\alpha}R \right) \psi_j -  \displaystyle \int \left(|D|^{\alpha}R_j\right)\eta \bigg|\\
		&\quad + 2\bigg| \displaystyle \int R_j^2\eta -R^2\eta\psi_j  \bigg| \\ & = |J_{3,1}| + |J_{3,2}| + |J_{3,3}|.
	\end{align*}
	By using the Cauchy-Schwarz inequality, we obtain 
	\begin{align*}
		|J_{3,2}| + |J_{3,3}|\leq& \|\eta \|_{L^2} \left(  \| \left(1-\psi_j\right) |D|^{\alpha}R_j  \|_{L^2} + \| R_j^2 \left(1-\psi_j \right)\|_{L^2} \right) \\& + \|\eta \|_{L^2} \left( \sum_{k\neq j} \|  \psi_j |D|^{\alpha}R_k \|_{L^2} + \sum_{(k,l)\neq(j,j)} \| R_kR_l \psi_j \|_{L^2} \right)
	\end{align*}
	From \eqref{eqQG2}, we rewrite $|D|^{\alpha}R_j=R_j^2-c_jR_{j}$. Thus, it follows from \eqref{est3} and \eqref{est5}	\begin{align}
		|J_{3,2}| + |J_{3,3}|\leq \frac{C}{(\beta t)^{\alpha}}.\label{J323}
	\end{align}
	 Now, we estimate $|J_{3,1}|$. By using the Cauchy-Schwarz inequality, 
	\begin{align}
		|J_{3,1}|&\leq \sum_{k\neq j}\bigg|  \int R_k (|D|^{\alpha}\eta)\psi_j  \bigg| + \bigg|  \int R_j (|D|^{\alpha}\eta )\left(1-\psi_j\right)  \bigg|\notag\\
		&\leq \lVert |D|^{\frac{\alpha}{2}}\eta\rVert _{L^2} \left( \sum_{k\neq j}\lVert |D|^{\frac{\alpha}{2}}\left( R_k \psi_j  \right) \rVert_{L^2}  +\lVert |D|^{\frac{\alpha}{2}} \left(R_j \left(1-\psi_j\right) \right) \rVert_{L^2} \right)\notag
	\end{align}
By interpolation $\lVert u\rVert_{H^{\frac\alpha2}}\leq \lVert u\rVert_{L^2}^{1-\frac\alpha2}\lVert u\rVert_{H^{1}}^{\frac\alpha2}$. Thus, we deduce from  \eqref{est4}, \eqref{est5}, \eqref{est6},\eqref{est7} that 
	\begin{align}
		|J_{3,1}|
		&\leq \lVert |D|^{\frac{\alpha}{2}}\eta\rVert _{L^2} \left( \sum_{k\neq j}\lVert  R_k \psi_j   \rVert_{L^2}^{1-\frac{\alpha}{2}} \lVert  R_k \psi_j   \rVert_{H^1}^{\frac{\alpha}{2}}  +\lVert R_j \left(1-\psi_j\right)  \rVert_{L^2}^{1-\frac{\alpha}{2}} \lVert R_j \left(1-\psi_j\right)  \rVert_{H^1}^{\frac{\alpha}{2}} \right) \leq \displaystyle\frac{C}{\left(\beta t\right)^{\alpha}}. \label{J31}
	\end{align}
	Gathering \eqref{J1}, \eqref{J2}, \eqref{J323} and \eqref{J31} we conclude the proof of \eqref{emest2}. 
	
	To prove \eqref{emest3}, gathering \eqref{emest1} and \eqref{emest2} we get 
	\begin{align*}
		\bigg| \left(E_j(t) + \displaystyle\frac{c_j}{2} M_j(t)\right) -\left( E(Q_{c_j})+\right.&\left. \displaystyle\frac{c_j}{2} M(Q_{c_j}) \right) -\displaystyle\frac{1}{2}H_j(t) \bigg| \\
		&\leq \displaystyle\frac{C}{\left(\beta t\right)^{\alpha}} + C\gamma\lVert \eta(t) \rVert_{H^{\frac{\alpha}{2}}}^2 + \bigg| \int \eta^2\left(R - R_j\right)\psi_j \bigg|
	\end{align*}
	Then by Cauchy-Schwarz, \eqref{est3} and the Sobolev imbedding $L^4(\R)\xhookrightarrow{}H^{\frac{1}{4}}(\R)$, we obtain that 
	\begin{align*}
		\bigg| \int \eta^2\left(R - R_j\right)\psi_j \bigg| \leq \| \eta\|_{L^4}^2 \sum_{k\neq j}\|R_k\psi_j \|_{L^2} \leq \frac{C}{\left(\beta t\right)^{\alpha}} ,
	\end{align*}
	which concludes the proof of \eqref{emest3}.
\end{proof}


\subsection{Control of the $R_j$ directions}
We recall $C_{*}$ is a positive constant changing from line to line and depending only on the parameters $\{c_1,\cdots,c_N\}$.
\begin{proposition}\label{R_j}
	For all $j\in\{1,\cdots,N\}$, $t_0\in[t^*,S_n]$,
	\begin{equation} \label{eta:Rj}
	\sum_{k=1}^{j}\bigg|\int\eta(t_0)R_k(t_0)\bigg|\leq\displaystyle\frac{C_{*}}{ t_0^{\alpha}} + C_{*} \| \eta(t_0) \|_{H^{\frac{\alpha}{2}}}^2 .
	\end{equation}
\end{proposition}

%

\begin{proof}[Proof of Proposition \ref{R_j}]
	The proof is by induction. For $j=1$, by \eqref{emest1} at time $t=t_0$ and $t=S_n$, we deduce that
	\begin{align*}
		2\displaystyle\int\eta(t_0)R_1(t_0) &\leq \frac{C}{\left(\beta t_0\right)^{\alpha}} -\int\eta^2\psi_1 -\int Q_{c_1}^2 + M_1(t_0)\\
		&\leq \frac{C}{\left(\beta t_0\right)^{\alpha}} -\int\eta^2\psi_1 -M_1(S_n) + M_1(t_0) +M_1(S_n) - \int  Q_{c_1}^2 \\
		&\leq \frac{C}{\left(\beta t_0\right)^{\alpha}} -\int\eta^2\psi_1 -M_1(S_n) + M_1(t_0).\\
	\end{align*}
	Moreover, by using estimate \eqref{M1}, we deduce that
	\begin{align*}
		2\displaystyle\int\eta(t_0)R_1(t_0) &\leq \frac{C}{\left(\beta t_0\right)^{\alpha}} -\int\eta^2\psi_1 \leq \frac{C}{\left(\beta t_0\right)^{\alpha}}.
	\end{align*}
	Now we want to obtain a lower bound of this scalar product. We recall that $\tilde{E}_k$ and $\sigma_0$ are respectively defined in \eqref{def:Etilde} and \eqref{defsigma}. By \eqref{emest1}, \eqref{emest2} at time $t=t_0$ and $t=S_n$,  we observe that 
	\begin{align*}
		\left(c_1 - 2\sigma_0\right)\int\eta R_1\geq& -\frac{C}{\left(\beta t_0\right)^{\alpha}} -C\gamma\lVert \eta \rVert_{H^{\frac{\alpha}{2}}}^2 - E_1(t_0) + E(Q_{c_1}) \\
		&+ \frac{1}{2}\int\eta (|D|^{\alpha}\eta)\psi_1- \int R\eta^2\psi_1\\
		&-\sigma_0 M_1(t_0) + \sigma_0\int Q_{c_1}^2 + \sigma_0\int\eta^2\psi_1\\ 
		\geq& -\frac{C_{*}}{ t_0^{\alpha}} -C_{*} \lVert \eta \rVert_{H^{\frac{\alpha}{2}}}^2- \widetilde{E}_1(t_0) + \widetilde{E}_1(S_n) \\
		&+ \frac{1}{2}\int\eta (|D|^{\alpha}\eta)\psi_1- \int R\eta^2\psi_1 + \sigma_0\int\eta^2\psi_1.\\
	\end{align*}
	Thus, we deduce from \eqref{M2} and the fact $\lVert R\psi_1 \rVert_{\infty}<C$
	\begin{align*}
		\left(c_1 - 2\sigma_0\right)\int\eta R_1\geq&-\frac{C_{*}}{t_0^{\alpha}} -C_{*}\|\eta\|_{H^{\frac\alpha2}}^2 + \frac{1}{2}\int\eta (|D|^{\alpha}\eta)\psi_1. 
	\end{align*}
	Note that replacing $\phi'_{j}$ by $\psi_j$ in \eqref{q1} we deduce 
	\begin{align*}
		\||D|^{\frac{\alpha}{2}}\left( \eta\psi_j \right) \|_{L^2} \leq C \| \eta \|_{H^{\frac{\alpha}{2}}} ,
	\end{align*}
	so that 
	\begin{align*}
		\left(c_1 - 2\sigma_0\right)\int\eta R_1\geq&-\frac{C_{*}}{ t_0^{\alpha}} -C_{*}\lVert \eta \rVert_{H^{\frac{\alpha}{2}}}^2. 
	\end{align*}
	Combining the lower and upper bound, we conclude that
	$$
	\bigg|\int\eta(t_0)R_1(t_0)\bigg|\leq\displaystyle\frac{C_{*}}{t_0^{\alpha}} + C_{*}\lVert \eta \rVert_{H^{\frac{\alpha}{2}}}^2.
	$$
	
	Now, we prove the inductive step. We assume that \eqref{eta:Rj} holds true for some $j \in \{1,\cdots,N-1\}$ and we prove it for $j+1$. Arguing similarly as in the case $j=1$, we deduce by \eqref{emest1} at time $t=t_0$ and $t=S_n$, \eqref{M1} and then the induction hypothesis \eqref{eta:Rj} in $j$, that
	\begin{align*}
		2\int\eta(t_0)R_{j+1}(t_0)\leq& \displaystyle\frac{C}{\left(\beta t_0\right)^{\alpha}}  + \sum_{k=1}^{j+1}\left( M_{k}(t_0) - M_{k}(S_n)\right)- \int\eta^2\psi_{j+1}\\
		& -\left[\sum_{k=0}^{j} \left(M_{k}(t_0) - M_{k}(S_n)\right) - \sum_{k=1}^{j}\int\eta^2\psi_k \right] - \sum_{k=1}^{j}\int\eta^2\psi_k\\
		\leq&\displaystyle\frac{C}{\left(\beta t_0\right)^{\alpha}} +2\sum_{k=1}^{j}\bigg|\int\eta R_k\bigg|\\
		\leq&\frac{C_{*}}{t_0^{\alpha}} + C_{*}\lVert \eta \rVert_{H^{\frac{\alpha}{2}}}^2.
	\end{align*}
	Arguing similarly as $j=1$ for the lower bound, we obtain from  \eqref{emest1}, \eqref{emest2} at time $t=t_0$ and $t=S_n$  that
	\begin{align*}
		\left(c_{j+1}-2\sigma_0 \right)\int\eta R_{j+1}\geq&- \frac{C_{*}}{t_0^{\alpha}} -\tilde{E}_{j+1}(t_0) +\tilde{E}_{j+1}(S_n) + \sigma_0 \int\eta^2\psi_{j+1}\\ 
		&+\frac{1}{2}\int\left(\eta |D|^{\alpha} - 2R\eta^2 \right)\psi_{j+1}- C_{*} \lVert \eta \rVert_{H^{\frac{\alpha}{2}}}^2\\
		\geq& -\displaystyle\frac{C_{*}}{ t_0^{\alpha}} - C_{*}\lVert \eta \rVert_{H^{\frac{\alpha}{2}}}^2 + \sigma_0\sum_{k=1}^{j+1} \int \eta^2\psi_k\\
		& + \sum_{k=1}^{j+1} \left(\widetilde{E}_{k}(S_n) - \widetilde{E}_{k}(t_0) \right)+  \frac{1}{2}\sum_{k=1}^{j+1}\int\left(\eta |D|^{\alpha}\eta - 2R\eta^2 \right)\psi_k\\
		& -\left[ \sum_{k=1}^{j} \left(\widetilde{E}_{k}(S_n) - \widetilde{E}_{k}(t_0) \right)+  \frac{1}{2}\sum_{k=1}^{j}\int\left(\eta |D|^{\alpha}\eta - 2R\eta^2 \right)\psi_k + \sigma_0\sum_{k=1}^{j}\int \eta^2\psi_k\right ]. 
	\end{align*}
	Thus, by using again \eqref{emest1}, \eqref{emest2} at time $t=t_0$ and $t=S_n$, \eqref{M2}, and then the induction hypothesis \eqref{eta:Rj} in $j$, it follows that
	\begin{align*}
		\left(c_{j+1}-2\sigma_0 \right)\int\eta R_{j+1}
		\geq& -\displaystyle\frac{C_{*}}{t_0^{\alpha}} - C_{*}\lVert \eta \rVert_{H^{\frac{\alpha}{2}}}^2 - \sum_{k=1}^{j}(c_{k}-2\sigma_0)\bigg| \int\eta R_k \bigg|\\
		\geq& -\displaystyle\frac{C_{*}}{t_0^{\alpha}} - C_{*}\lVert \eta \rVert_{H^{\frac{\alpha}{2}}}^2. \\
	\end{align*}
	This concludes the proof of \eqref{eta:Rj} in $j+1$, and thus the proof of Proposition \ref{R_j} by induction. 
\end{proof}


\subsection{Proof of Proposition \ref{s=1/2}}
Recalling the notation $\eta=u-R$, it suffices to prove that 
\begin{align}
	\lVert \eta \rVert_{H^{\frac{\alpha}{2}}}^2\leq \displaystyle\frac{C_{*}}{ t_0^{\alpha}}.\label{errorest}
\end{align}
\begin{proof}[Proof of \eqref{errorest}]
	The proof of the estimate \eqref{errorest} relies on the quadratic form $H_j(t)$ defined in \eqref{quadratic} 
	On the one hand, we have from \eqref{emest3},  
	\begin{align}
		\sum_{j=1}^{N}\frac{1}{c_j^{2}}H_j(t_0)\leq& \sum_{j=1}^{N}\frac{1}{c_j^{2}}\left(E_j(t_0) + \frac{c_j}{2}M_j(t_0) \right) - \sum_{j=1}^{N}\frac{1}{c_j^{2}}\left(E(Q_{c_j}) + \frac{c_j}{2}M(Q_{c_j}) \right) \label{sum}\\
		&+ \displaystyle\frac{C_{*}}{t_0 ^{\alpha}} + C_{*}\gamma \lVert \eta \rVert_{H^{\frac{\alpha}{2}}}^2. \notag
	\end{align}
	On the other hand, by a direct resummation argument, we observe that 
	$$
	\sum_{j=1}^{N}\frac{1}{c_j^{2}}\widetilde{E}_j=\sum_{j=1}^{N-1}\left( \frac{1}{c_j^{2}} - \frac{1}{c_{j+1}^{2}}\right)\sum_{k=1}^{j}\widetilde{E}_k + \frac{1}{c_{N}^{2}}\sum_{k=1}^{N}\widetilde{E}_k
	$$
	and
	\begin{align*}
		\sum_{j=1}^{N}\frac{1}{c_j^{2}}\left(\frac{c_j}{2} - \sigma_0\right)M_j=&\sum_{j=1}^{N-1}\left[\frac{1}{2}\left( \frac{1}{c_j} - \frac{1}{c_{j+1}}\right)\left(1- 2\sigma_0 \left( \frac{1}{c_j} + \frac{1}{c_{j+1}}\right) \right) \right]\sum_{k=1}^{j}M_k \\
		&+\frac{1}{2c_{N}}\left(1-2\frac{\sigma_0}{c_N} \right)\sum_{k=1}^{N}M_k.\\
	\end{align*}
	Combining these two identities, since $\widetilde{E}_j=E_j+\sigma_0M_j$, we deduce that
	\begin{align*}
		\sum_{j=1}^{N}\frac{1}{c_j^{2}}\left(E_j+ \frac{c_j}{2} M_j\right)=&\sum_{j=1}^{N-1}\left( \frac{1}{c_j^{2}} - \frac{1}{c_{j+1}^{2}}\right)\sum_{k=1}^{j}\widetilde{E}_k + \frac{1}{c_{N}^{2}}\sum_{k=1}^{N}\widetilde{E}_k\\
		&+\sum_{j=1}^{N-1}\left[\frac{1}{2}\left( \frac{1}{c_j} - \frac{1}{c_{j+1}}\right)\left(1- 2\sigma_0 \left( \frac{1}{c_j} + \frac{1}{c_{j+1}}\right) \right) \right]\sum_{k=1}^{j}M_k \\
		&+\frac{1}{2c_{N}}\left(1-2\frac{\sigma_0}{c_N} \right)\sum_{k=1}^{N}M_k.\\
	\end{align*}
	Note that all the coefficients in front of the partial sums on the right hand side of the above estimate are positive by definition of $\sigma_0$ in \eqref{defsigma}. Therefore, we deduce from \eqref{sum}, \eqref{emest1}, \eqref{emest2} at time $t=t_0$ and $t=S_n$ and the monotonicity results \eqref{M1} and \eqref{M2} in Proposition \ref{mono},  that 
	\begin{align}
		\sum_{j=0}^{N}\frac{1}{c_j^{2}}H_j(t_0)\leq \displaystyle\frac{C_{*}}{ t_0^{\alpha}} + C_{*}\gamma\lVert \eta \rVert_{H^{\frac{\alpha}{2}}}^2.\label{Hj1}
	\end{align}
	On the other hand, by Corollary \ref{almostcoer} and \eqref{Ortg}, there exists $\lambda_0>0$ such that 
	\begin{align*}
		\sum_{j=0}^{N}\frac{1}{c_j^{2}}H_j(t_0)\geq& \lambda_0\lVert \eta \rVert_{H^{\frac{\alpha}{2}}}^2 -\frac{C_{*}}{  t_0 ^{\alpha}}
		- \frac{1}{\lambda_0}\sum_{j=0}^{N} \left(\int \eta(t_0)R_j(t_0) \right)^2. 
	\end{align*}
	The control of the $R_j$ directions derived in  Proposition \ref{R_j} yields 
	\begin{align}
		\sum_{j=0}^{N}\frac{1}{c_j^{2}}H_j(t_0)\geq& \lambda_0\lVert \eta \rVert_{H^{\frac{\alpha}{2}}}^2-\frac{C_{*}}{ t_0 ^{\alpha}}   - \frac{1}{\lambda_0}\displaystyle\frac{C_{*}}{ t_0^{2\alpha}} - \frac{1}{\lambda_0}C_{*}\lVert \eta \rVert_{H^{\frac{\alpha}{2}}}^4 . \label{Hj2}
	\end{align}								     
	Therefore, we conclude the proof of \eqref{errorest} by combining \eqref{Hj1} and \eqref{Hj2}, which finishes the proof of Proposition \ref{s=1/2}.

\end{proof}


\appendix

\section{Weak continuity of the flow}\label{Prweakflow}

In this appendix, we give the proof of Lemma \ref{weakflow}, where the IVP associated to \eqref{GC} is globally well-posed in the energy space (see \cite{molinet2018well}). We follow a general argument given by L. Molinet \cite{Mol} (see also \cite{goubet2009global}). 

\begin{proof}[Proof of Lemma \ref{weakflow}]
Let $\alpha>\frac67$. For $T>0$, we recall the definition of $Y_t$ the resolution space  $$Y_T:= C^{0}([0,T]:H^{\frac{\alpha}2}(\R))\cap X_{T}^{\frac{\alpha}2-1,1}\cap L^{2}((0,T),W^{\frac{\alpha}4-\frac12-,\infty}(\R)),$$ and let $\|\cdot\|_{Y_T}$ be the norm associated to $Y_T$.

Assume $z_{0,n}\rightharpoonup z_0$ in $H^{\frac{\alpha}2}(\mathbb R)$. By the Banach-Steinhaus theorem, we deduce that there exists $C>0$ such that $\| z_{0,n} \|_{H^{\frac{\alpha}2}}\leq C$. Moreover,  by the global well-posedness result Theorem \ref{GWP:theo},  there exists $C>0$ such that the solution $z_n$ of \eqref{GC}, associated to $z_{0,n}$, verifies
\begin{equation}
	\|z_n(t) \|_{Y_T}\leq C, \quad \forall t\in[0,T].\label{bornunif}
	\end{equation} 
Thus, by Banach-Alaoglu's theorem, there exists 
$$ z\in Z_T:= L^{\infty}((0,T):H^{\frac{\alpha}2}(\R))\cap X_{T}^{\frac{\alpha}2-1,1}\cap L^{2}((0,T),W^{\frac{\alpha}4-\frac12-,\infty}(\R)), $$ such that $z_n\stackrel{\ast}{\rightharpoonup} z$ in $ L^{\infty}([0,T]:H^{\frac{\alpha}2}(\R))$, up to extracting a subsequence.
By \eqref{bornunif}, we get 
$$\| D^{\alpha}\partial_x z_n \|_{L^{\infty}_{T}H^{-\frac{\alpha}{2}-1}}\leq C,$$ 
  and since $z_n^2\in L^1(\R)\xhookrightarrow{}H^{-\frac{1}{2}^{-}}(\R)$, we have that 
 $$\| \partial_x (z_n^2) \|_{L^{\infty}_{T}H^{-\frac{3}{2}^{-}}}\leq C\| z_n\|_{L^2}^2\leq C.$$
  Then, we obtain, by \eqref{GC} that 
 \begin{equation} \label{bound:dzndt}
 \| \partial_t z_n\|_{L^{\infty}_{T}H^{\min(-\frac{3}{2}^{-},-\frac{\alpha}{2}-1 )}}\leq C.
 \end{equation}
 Therefore, by the Aubin-Lions theorem (Theorem 1.71 in \cite{novotny2004introduction}), we deduce that $z_n\rightarrow z$ in $L^{2}([0,T]:L^{2}([-k,k])$, for all $k\in\N$. In particular, this implies $z_n^2\rightarrow z^2$ in $L^{1}([0,T]:L^{1}([-k,k]))$. 

Now, since $z_n$ is a weak solution of \eqref{GC} in the distributional sense satisfying $z_n(0,\cdot)=z_{0,n}$, we know that for all $\phi\in C^{\infty}_{c}([-T,T]\times\R)$, 
\begin{align*}
\int_{0}^{T}\int_{\R}\left(\partial_t\phi -\partial_x|D|^{\alpha}\phi\right) z_n dx dt + \int_{0}^{T}\int_{\R} \left(\partial_x \phi\right) z_n^2dxdt -\int_{\R}\phi(0,x) z_{0,n}(x)dx   =0.	
\end{align*}    
Thus passing to the limit, we conclude that
\begin{align*}
	\int_{0}^{T}\int_{\R}\left(\partial_t\phi -\partial_x|D|^{\alpha}\phi\right) z dx dt + \int_{0}^{T}\int_{\R} \left(\partial_x \phi\right) z^2dxdt -\int_{\R}\phi(0,x) z_{0}(x)dx   =0	,
\end{align*}    
which proves that $z$ is a weak solution of \eqref{GC} corresponding to the initial datum $z_0$. Since $z \in Z_T$, we deduce from the uniqueness part in Theorem \ref{GWP:theo} that $z \in Y_T$ and that $z_n\stackrel{\ast}{\rightharpoonup} z$ in $ L^{\infty}([0,T]:H^{\frac{\alpha}2}(\R))$ for the whole sequence.
 
Finally, let $\psi \in C_c^{\infty}(\mathbb R)$. It follows from the Arzela-Ascoli theorem and the bounds \eqref{bornunif}-\eqref{bound:dzndt} that the function $\displaystyle v_n:t\in[0,T]\mapsto \int_{\R}\psi(x)z_n(t,x)dx$ converges up to a subsequence in $C^{0}([0,T]:\R)$. Moreover, by uniqueness, this limit holds for the whole sequence and is equal to  $\displaystyle\int_{\mathbb R}z(t,x)\psi(x)dx$, which implies that $z_n(t) \rightharpoonup z(t)$ in $H^{\frac{\alpha}2}(\mathbb R)$ for all $t \in [0,T]$.  

\end{proof}

\section{Pseudo-differential toolbox}
In this section $u$ will denoted a function in $\mathcal{S}(\R)$.
First, we recall some well-known results on pseudo-differential operators (see \cite{alinhac2012operateurs}, or \cite{hormander2007analysis} chapter $18$).
Let $D=-i\partial_x$. We define the symbolic class $\mathcal{S}^{m,q}$ by 
\begin{align*}
	a(x,\xi)\in \mathcal{S}^{m,q}\iff \left\{ 
	\begin{array}{ll}
		a\in C^{\infty}(\R^2)\\
		\forall k,\beta\in \N,\exists C_{k,\beta}>0 \text{ such that } |\partial_x^k\partial_{\xi}^{\beta}a(x,\xi)|\leq C_{k,\beta}\langle x \rangle^{q-k} <\xi>^{m-\beta} 
	\end{array} 
	\right.
\end{align*}
For all $u\in \mathcal{S}(\R)$, we set the operator associated to the symbol $a(x,\xi)\in \mathcal{S}^{m,q}$ by 
$$
a(x,D)u:=\displaystyle\frac{1}{2\pi}\int e^{ix\xi} a(x,\xi) \mathcal{F}(u)(\xi) d\xi.
$$ 
We state the three following results 
\begin{enumerate}
	\item  Let $a\in \mathcal{S}^{m,q}$, there exists $C>0$, such that for all $u\in \mathcal{S}(\R)$ 
	\begin{align}
		\| a(x,D)u \|_{L^2}\leq C\|\langle x \rangle^q\langle D \rangle^m u \|_{L^2}. \label{nps}
	\end{align}
	\item Let $a\in \mathcal{S}^{m,q}$ and $b\in \mathcal{S}^{m',q'}$, then there exists $c\in \mathcal{S}^{m+m',q+q'}$ such that 
	\begin{align}
		a(x,D)b(x,D)=c(x,D).\label{pps}
	\end{align}
	\item If $a\in \mathcal{S}^{m,q}$ and $b\in \mathcal{S}^{m',q'}$ are two operators we define the commutator by $[a(xD),b(x,D)]:=a(x,D)b(x,D)-b(x,D)a(x,D)$.
	Moreover there exists $c\in \mathcal{S}^{m+m'-1,q+q'-1}$ such that 
	\begin{align}
		[a(x,D),b(x,D)]=c(x,D).\label{comps}
	\end{align}
\end{enumerate}
As a consequence of \eqref{pps},  $\langle D \rangle^m\langle x \rangle^q\langle D \rangle^{-m}\in \mathcal{S}^{0,q}$. Therefore, by \eqref{nps}, we have 
\begin{align*}
	\|\langle D \rangle^m\langle x \rangle^qu \|_{L^2}&=\|\langle D \rangle^m\langle x \rangle^q\langle D \rangle^{-m}\langle D \rangle^mu \|_{L^2}\\
	&\leq C_2\|\langle x \rangle^q\langle D \rangle^m u \|_{L^2},
\end{align*}
for $C_2>0$.
By the same computations with $\langle x \rangle^q$ instead of $\langle D \rangle^{m}$, there exists $C_1>0$ such that
\begin{align*}
	C_1\|\langle x \rangle^q\langle D \rangle^m u \|_{L^2}&\leq  \|\langle D \rangle^m\langle x \rangle^qu \|_{L^2}.
\end{align*}
Gathering these two estimates, we conclude that
\begin{align}
	C_1\|\langle x \rangle^q\langle D \rangle^m u \|_{L^2}&\leq  \|\langle D \rangle^m\langle x \rangle^qu \|_{L^2}\leq C_2\|\langle x \rangle^q\langle D \rangle^m u \|_{L^2}.  \label{oh}
\end{align}

\section{Proof of the weighted commutator estimates}\label{B.2}
This section is devoted to the proofs of Lemmas \ref{commG}, \ref{eqnorm} and \ref{esttc}.

In this section $u$ is a function in $H^{\frac{\alpha}{2}}(\R)$. Let $\phi$ be defined as in \eqref{defphi} and for $A>1$, $\phi_{j,A}$ as in \eqref{mj}.
Moreover, we define 
\begin{align}
	\Phi(x)=\sqrt{|\phi'(x)|} \sim \langle x\rangle^{-\frac{1+\alpha}2} \quad \text{and} \quad \Phi_{j,A}=\sqrt{|\phi'_{j,A}|}. 
\end{align}

Finally, let $\chi\in C^{\infty}_c(\R)$ such that $\chi(\xi)=1$ on $[-1,1]$ and $\chi(\xi)=0$ on $[-2,2]^{c}$.

\medskip 

The proof of Lemma  \ref{commG} is an extension of the proof of Lemmas 6 and 7 in \cite{kenig2011local}. 
Note that, while the estimates in Lemmas 6 and 7 in \cite{kenig2011local} are stated for $\alpha \in [1,2]$, their proofs extend directly to the case $\alpha \in (0,2)$. This yields the estimates  \eqref{sc1G} and \eqref{sc2G}. However, since the estimates \eqref{nsc1G} and \eqref{nsc2G} are not symmetric in $u$,  we cannot use the Claim 3 in \cite{kenig2011local}. Instead, we use the following estimates (which are also derived from the techniques in \cite{kenig2011local}).    
\begin{lemma}\label{K1}
	Let $\alpha\in(0,2)$, $T=|D|^{\alpha}\phi D- D\phi |D|^{\alpha}$. Then, there exists $C>0$ such that for all $u,v\in\mathcal{S}(\R)$ we have
	\begin{align}
		i\int (Tu) v  = -(\alpha-1 ) \int |D|^{\frac{\alpha}{2}}\left(u\Phi \right)|D|^{\frac{\alpha}{2}}\left(v\Phi \right)   + R, \label{k1}
	\end{align}
	with 
	\begin{equation} \label{est:R}
	|R|\leq \begin{cases} \displaystyle C\int (u^2+v^2)|\phi'|, & \text{if} \  \alpha\in(0,1], \\
	\displaystyle C\left( \int u^2|\phi'|  +\frac{1}{A^{\frac{\alpha}{2}}} \int  (|D|^{\frac{\alpha}{2}}u)^2|\phi'| + A^{\frac{\alpha}{2}}\int v^2|\phi'| \right), & \text{if} \ \alpha\in(1,2) ,
	\end{cases}
	\end{equation}
	 for all $A>1$.

\end{lemma}

\begin{proof}
The proof of \eqref{k1} is a combination of the proofs of Claim 1 and Claim 2 in \cite{kenig2011local}. Following \cite{kenig2011local}, we split $T=T_1+T_2$, where
\begin{align*}
T_1&=|D|^{\alpha}(1-\chi(D))\phi D- D\phi (1-\chi(D))|D|^{\alpha} ;\\ 
T_2&=|D|^{\alpha}\chi(D)\phi D- D\phi \chi(D)|D|^{\alpha}.
\end{align*}

First, we have arguing exactly as in Claim 2 in \cite{kenig2011local} that, for any $\alpha \in (0,2)$,
\begin{equation} \label{est:T_2}
i \int (T_2u)v=-(\alpha-1) \int |D|^{\frac{\alpha}{2}}\chi(D)\left(u\Phi \right)|D|^{\frac{\alpha}{2}}\left(v\Phi \right)   + R_2, 
\end{equation}
with 
$$
|R_2|\leq \displaystyle C\int (u^2+v^2)|\phi'| . 
$$

Now, we deal with the operator $T_1$. Let us define $a(x,\xi)=\varphi(x)|\xi|^{\alpha}(1-\chi(\xi)) \in \mathcal{S}^{\alpha,0}$. Then, following the computations in the proof of Claim 1 in \cite{kenig2011local}, we have
\begin{equation} \label{est:T_1}
i \int (T_1u)v=-(\alpha-1) \int |D|^{\frac{\alpha}{2}}(1-\chi(D))\left(u\Phi \right)|D|^{\frac{\alpha}{2}}\left(v\Phi \right)
+ \int(\widetilde{T}_1u) v  + R_1, 
\end{equation}
with 
\begin{align*}
\widetilde{T}_1 &=-\frac{i}2(\partial_x^2\partial_{\xi}^2a)(x,D)D-\Phi[\Phi,|D|^{\alpha}(1-\chi(D))]\\
|R_1|&\leq \displaystyle C\int (u^2+v^2)|\phi'| . 
\end{align*}
To estimate $\left|\int(\widetilde{T}_1u) v\right|$, we cannot use Claim 3 in \cite{kenig2011local} due to the lack of symmetry. Instead, we use classical pseudo-differential calculus estimates. Observe that the symbol $t_1(x,\xi)$ of $\widetilde{T}_1$ belongs to the class $\mathcal{S}^{\alpha-1,-\alpha-2}$. In the case, $0<\alpha<1$, $t_1(x,\xi) \in \mathcal{S}^{0,-(\alpha+1)}$. Thus the Cauchy-Schwarz inequality and \eqref{nps} yield
\[ \left|\int(\widetilde{T}_1u) v\right| =\left|\int(\Phi^{-1}\widetilde{T}_1u) \Phi v\right| \leq \|\Phi^{-1}\widetilde{T}_1u\|_{L^2}\|\Phi v\|_{L^2} \leq \displaystyle C\int (u^2+v^2)|\phi'| .
\]
In the case $1<\alpha<2$,  $t_1(x,\xi) \in  \mathcal{S}^{\frac{\alpha}{2},-(1+\alpha)}$.  
By Cauchy-Schwarz' inequality, \eqref{nps}, and then Young's inequality, we get 
$$
\left|\int(\widetilde{T}_1u) v\right| \leq \|\Phi^{-1}\widetilde{T}_1u\|_{L^2}\|\Phi v\|_{L^2} \leq \displaystyle C\left(\frac{1}{A^{\frac{\alpha}{2}}} \int  (\langle D \rangle^{\frac{\alpha}{2}}u)^2|\phi'| + A^{\frac{\alpha}{2}}\int v^2|\phi'| \right).
$$
for any $A>1$.  Moreover by pseudo-differential calculus \eqref{pps}, \eqref{nps}, and since the symbols of the operators $\Phi\chi(D)\langle D \rangle^{\frac{\alpha}{2}}\Phi^{-1}$ and $\Phi(1-\chi(D))\langle D \rangle^{\frac{\alpha}{2}}|D|^{-\frac{\alpha}{2}}\Phi^{-1}$ belong to $\mathcal{S}^{0,0}$,
\begin{align}
	\int \left(\langle D \rangle^{\frac{\alpha}{2}} u\right)^2 |\phi'| &\leq 2\int \left(\chi(D)\langle D \rangle^{\frac{\alpha}{2}} u\right)^2 |\phi'|  + 2\int \left((1-\chi(D))\langle D \rangle^{\frac{\alpha}{2}} u\right)^2 |\phi'|  \notag\\
	&\leq 2\int  \left(\Phi\chi(D)\langle D \rangle^{\frac{\alpha}{2}}\Phi^{-1} \left(\Phi u\right)\right)^2 \notag \\
	&\quad + 2\int \left(\Phi(1-\chi(D))\langle D \rangle^{\frac{\alpha}{2}}|D|^{-\frac{\alpha}{2}}\Phi^{-1} \left(\Phi |D|^{\frac{\alpha}{2}}u\right)\right)^2  \notag \\
	&\leq C\left( \int  u^2 |\phi'|  + \int \left(|D|^{\frac{\alpha}{2}} u\right)^2 |\phi'| \right) \label{dec}.
\end{align} 
Gathering these estimates concludes the proof of Lemma \ref{K1}.
\end{proof}

\begin{lemma}\label{lemma:S}
	Let $\alpha\in(0,2)$, $S=\Phi[\Phi,|D|^{\alpha}]$. Then, there exists $C>0$ such that for all $u,v\in\mathcal{S}(\R)$ we have
	\begin{equation} \label{est:S}
	\left|\int (Su) v  \right|+\left|\int (Sv) u  \right|\leq \begin{cases} \displaystyle C\int (u^2+v^2)|\phi'|, & \text{if} \  \alpha\in(0,1], \\
	\displaystyle C\left( \int u^2|\phi'|  +\frac{1}{A^{\frac{\alpha}{2}}} \int  (|D|^{\frac{\alpha}{2}}u)^2|\phi'| + A^{\frac{\alpha}{2}}\int v^2|\phi'| \right), & \text{if} \ \alpha\in(1,2) ,
	\end{cases}
	\end{equation}
	 for all $A>1$.	 
\end{lemma}

\begin{proof} We split $S=S_1+S_2,$ where
\begin{align*}
S_1&=\Phi[\Phi,|D|^{\alpha}(1-\chi(D))]   ;\\ 
S_2&=\Phi[\Phi,|D|^{\alpha}\chi(D)].
\end{align*}

We first deal with the high frequency terms involving $S_1$. 
Since $1-\chi$ is supported outside $0$, $S_1$ is a pseudo-differential operator of symbol  $s_1(x,\xi)$, which belongs to the class $\mathcal{S}^{\alpha-1,-\alpha-2}$. In the case, $0<\alpha\leq1$, $s_1(x,\xi) \in \mathcal{S}^{0,-(\alpha+1)}$ and in the case, $1<\alpha<2$, $s_1(x,\xi) \in \mathcal{S}^{\frac{\alpha}{2},-(\alpha+1)}$. Thus by arguing as in the proof of Lemma \ref{K1}, we deduce that
\begin{align*}
\left| \int (S_1 u)v \right|\leq 
\begin{cases}
C\displaystyle\int (u^2+v^2)|\phi'|, &  \text{ if } 0<\alpha\leq 1,\\
C\displaystyle \left(\frac{1}{A^{\frac{\alpha}{2}}} \int  (\langle |D| \rangle^{\frac{\alpha}{2}}u)^2|\phi'| + A^{\frac{\alpha}{2}}\int v^2|\phi'| \right), & \text{ if } 1<\alpha<2 .
\end{cases}
\end{align*}
Observe that the same estimate also holds for $\left| \int (S_1 v)u \right|$. Indeed, the proof is exactly the same in the case $0<\alpha \le 1$. In the case $1<\alpha<2$, $\Phi^{-1}\langle |\xi| \rangle^{-\frac{\alpha}2}s_1(x,\xi) \in \mathcal{S}^{0,-\frac{\alpha+1}2}$, so that 
\begin{align*}
\left| \int (S_1 v)u \right| \leq \|\Phi^{-1}\langle |D| \rangle^{-\frac{\alpha}2}S_1 v\|_{L^2} \|\Phi\langle |D| \rangle^{\frac{\alpha}2}u \|_{L^2}\leq \left(\frac{1}{A^{\frac{\alpha}{2}}} \int  (\langle D \rangle^{\frac{\alpha}{2}}u)^2|\phi'| + A^{\frac{\alpha}{2}}\int v^2|\phi'| \right)
\end{align*}
Moreover, by using \eqref{dec}, we deduce that 
\begin{align*}
	\left|\int (S_1u) v  \right|+\left|\int (S_1v) u  \right|\leq \begin{cases} \displaystyle C\int (u^2+v^2)|\phi'|, & \text{if} \  \alpha\in(0,1], \\
	\displaystyle C\left( \int u^2|\phi'|  +\frac{1}{A^{\frac{\alpha}{2}}} \int  (|D|^{\frac{\alpha}{2}}u)^2|\phi'| + A^{\frac{\alpha}{2}}\int v^2|\phi'| \right), & \text{if} \ \alpha\in(1,2) ,
	\end{cases}
\end{align*}

Now we deal with the low frequency term involving $S_2$. We follow the proof given in \cite{kenig2011local} for the same type of operator. We remark that $|D|^{\alpha}\chi(D)u=k\ast u$, with $\hat{k}=|\xi|^{\alpha}\chi(\xi)$. Then, we can rewrite 
\[ [\Phi,|D|^{\alpha}\chi(D)]u= \int k(x-y) \left( \Phi(x) - \Phi(y) \right) u(y)dy. \] 
We want to prove that the operator defined by the kernel
$$
\Lambda(x,y)= k(x-y)\left( \Phi(x) -\Phi(y) \right) \Phi^{-1}(y),
$$
is bounded in $L^2(\R)$. For this, we need the 3 following results.

\begin{thm}[Schur's test \cite{halmos2012bounded}, Theorem 5.2]\label{K3}
	Let $p, q$ be two non-negative measurable functions. If there exists $\alpha,\beta>0$ such that 
	\begin{enumerate}
		\item $\displaystyle\int_{\R}|K(x,y)| q(y) dy \leq \alpha p(x) \text{ a.e } x\in \R $.
		\item $\displaystyle\int_{\R}|K(x,y)| p(x) dx \leq \beta q(y) \text{ a.e } y\in \R$.
	\end{enumerate}
	Then $Tf:=\displaystyle\int_{\R}K(x,y)f(y) dy$ is a bounded operator on $L^2(\R)$.
\end{thm}

\begin{claim}[\cite{kenig2011local} Claim 8]\label{K4}
	There exists $C>0$ such that
	$$
	|\Phi(x) - \Phi(y)| \leq C\frac{|x-y|}{\left( \langle x \rangle + \langle y \rangle \right)^{\frac{\alpha+3}{2}}} \text{ if } |x-y|\leq \frac{1}{2}\left( \langle x \rangle + \langle y \rangle \right),
	$$
	$$
	|\Phi(x) - \Phi(y)| \leq  \frac{1}{\langle x \rangle ^{\frac{1+\alpha}{2}}} + \frac{1}{\langle y \rangle ^{\frac{1+\alpha}{2}}}  \text{ if } |x-y|\geq \frac{1}{2}\left( \langle x \rangle + \langle y \rangle \right).
	$$
\end{claim}

\begin{lemma}[\cite{kenig2011local}, Lemma A.2]\label{K2}
	Let $p$ be a homogeneous function of degree $\beta>-1$. Let $\chi\in C^{\infty}_0(\R)$ such that $0\leq \chi\leq1$, $\chi(\xi)=1$ if $|\xi|<1$ and $\chi(\xi)=0$ if $|\xi|>2$. Let 
	$$
	k(x)= \frac{1}{2\pi}\int e^{ix\xi}p(\xi)\chi(\xi) d\xi.
	$$
	Then for all $q\in\mathbb{N}$, there exists $C_q>0$ such that, for all $x\in\R$,
	\begin{align}
		|\partial_x^q k(x)|\leq \frac{C_q}{\langle x \rangle^{\beta+q+1}}.
	\end{align}
\end{lemma}

Let $\Lambda= \Lambda_1 +\Lambda_2$, where $\Lambda_1$ and $\Lambda_2$ are restricted respectively to the regions $|x-y|\leq \frac{1}{2}\left( \langle x \rangle + \langle y \rangle \right)$ and $|x-y|\geq \frac{1}{2}\left( \langle x \rangle + \langle y \rangle \right)$. By using  Claim \ref{K4} and Lemma \ref{K2},
\begin{align*}
	|\Lambda_1(x,y)|&\leq C\frac{1}{\langle x-y \rangle^{1+\alpha}}\frac{|x-y|}{\left( \langle x \rangle + \langle y \rangle \right)^{\frac{\alpha+3}{2}}}\langle y \rangle^{\frac{1+\alpha}{2}}\\
	&\leq C \frac{1}{\langle x-y \rangle^{1+\alpha}} .
\end{align*}
Then, by Theorem \ref{K3}, with $p=q=1$, $\Lambda_1$ is the kernel of a bounded operator in $L^2$. Now, we deal with $\Lambda_2$. By using  Claim \ref{K4} and Lemma \ref{K2}, 
\begin{align*}
	|\Lambda_2(x,y)|&\leq C\frac{1}{\langle x-y \rangle^{1+\alpha}}\left(\frac{1}{\langle x \rangle^{\frac{1+\alpha}{2}}} + \frac{1}{\langle y \rangle^{\frac{1+\alpha}{2}}} \right)\langle y \rangle^{\frac{1+\alpha}{2}}\\
	&\leq C\frac{1}{\langle x-y \rangle^{1+\alpha}} + C\frac{\langle y \rangle^{\frac{1+\alpha}{2}}}{\langle x-y \rangle^{1+\alpha} \langle x \rangle^{\frac{1+\alpha}{2}}}\\
	&\leq \Lambda_3(x,y) + \Lambda_4(x,y).
\end{align*}
Then, by Theorem \ref{K3}, with $p=q=1$, $\Lambda_3$ is the kernel of a bounded operator in $L^2$. 
We compute
\begin{align*}
	\int \Lambda_4(x,y)\langle x \rangle^{-\frac{1}{2}}dx \leq C \langle y \rangle^{-\frac{1+\alpha}{2}}, \quad \int \Lambda_4(x,y)\langle y \rangle^{-\frac{1+\alpha}{2}}dy \leq C \langle x \rangle^{-\frac{1}{2}}.
\end{align*}
Then by Theorem \ref{K3}, we deduce that $\Lambda_4$ is the kernel of a bounded operator in $L^2$. Gathering these estimates, we conclude that
\begin{align*}
	\|[\Phi,|D|^{\alpha}\chi(D)]u\|_{L^2}\leq C\|u\Phi \|_{L^2}.
\end{align*}
Therefore, by Young's inequality, we get 
\begin{align*}
\bigg|\int (S_2u)v\bigg|=\bigg|\int v\Phi[\Phi,|D|^{\alpha}\chi(D)]u\bigg| \leq C\left( \int v^2|\phi'|  + \int\left([\Phi,|D|^{\alpha}\chi(D)]u\right)^{2} \right) \leq  C \int \left( u^2 + v^2 \right)|\phi'| . 
\end{align*}
The estimate for $\bigg|\int (S_2v)u\bigg|$ is similar. This concludes the proof of Lemma \ref{lemma:S}.
\end{proof}

\begin{proof}[Proof of \eqref{nsc2G}] 
By integration by parts, we get 
\begin{align*}
	\displaystyle\int \left( \left( |D|^{\alpha} u \right)\partial_xv + \left( |D|^{\alpha} v \right)\partial_xu \right) \phi &= \int \left(|D|^{\alpha}(\phi \partial_x u)-\partial_x(\phi  |D|^{\alpha}u)  \right) v  = i\int Tu v,
\end{align*}
with $T=|D|^{\alpha}\phi D- D \phi  |D|^{\alpha}$. Hence, we deduce from \eqref{k1} that
\begin{align*}
	\bigg| \displaystyle\int \left( \left( |D|^{\alpha} u \right)\partial_xv + \left( |D|^{\alpha} v \right)\partial_xu \right) \phi  +& (\alpha-1)\int |D|^{\frac{\alpha}{2}}\left(u\sqrt{|\phi'|} \right) |D|^{\frac{\alpha}{2}}\left(v\sqrt{|\phi'|} \right)   \bigg| =|R| ,
\end{align*}
where $|R|$ satisfies \eqref{est:R}.
Therefore, we conclude the proof of \eqref{nsc2G} by performing the change  of variable $x'=\frac{x-m_j}{A}$. 
\end{proof}

\begin{proof}[Proof of  \eqref{nsc1G}] By direct computation we get 
\begin{align*}
	\displaystyle\int \left( \left( |D|^{\alpha} u \right)v - \left( |D|^{\alpha} v \right)u \right) |\phi'|
	=\int v\Phi[\Phi,|D|^{\alpha}]u - \int u\Phi[\Phi,|D|^{\alpha}]v =\int (Su)v-\int (Sv)u,
\end{align*}
where $S=\Phi[\Phi,|D|^{\alpha}]$. Therefore, we conclude the proof of \eqref{nsc1G} by using Lemma \ref{lemma:S} and performing the change  of variable $x'=\frac{x-m_j}{A}$. 
\end{proof}

This finishes the proof of Lemma \ref{commG}. Now we prove Lemma \ref{eqnorm}.

\begin{proof}[Proof of Lemma \ref{eqnorm}]
By direct computations, we obtain 
\begin{align}
	\int \left(|D|^{\alpha}\left(u\Phi_{j,A} \right)\right)^2 -  \int \left( |D|^{\alpha}u \right)^{2} |\phi'_{j,A}| &= \int (|D|^{2\alpha}u) u |\phi'_{j,A}|  -\int \left( |D|^{\alpha}u \right)^{2} |\phi'_{j,A}| - \int u\Phi_{j,A}[\Phi_{j,A},|D|^{2\alpha}]u \nonumber \\
	&= \int (|D|^{2\alpha}u) u |\phi'_{j,A}|  -\int \left( |D|^{\alpha}u \right)^{2} |\phi'_{j,A}| \nonumber \\
	& \quad - \frac{1}{2} \int u\left[\Phi_{j,A},\left[\Phi_{j,A},|D|^{2\alpha}\right]\right]u . \label{comm:eqnorm}
\end{align} 
By  applying \eqref{nsc1G} to $v=|D|^{\alpha}u$, we get that 
\begin{align}\label{a1}
	\bigg| \int (|D|^{2\alpha}u) u |\phi'_{j,A}|  -\int \left( |D|^{\alpha}u \right)^{2} |\phi'_{j,A}| \bigg| &\leq\begin{cases} \displaystyle\frac{C}{A^{\alpha}}\int \left( u^2 + \left(|D|^{\alpha}u\right)^2 \right)|\phi'_{j,A}|  , &\text{if }\alpha\in(0,1], \\
		\displaystyle\frac{C}{A^{\frac{\alpha}{2}}}\int \left( u^2 +\left(|D|^{\frac{\alpha}{2}}u\right)^2  + \left(|D|^{\alpha}u\right)^2 \right)|\phi'_{j,A}|,   &\text{if }\alpha\in(1,2).
		\end{cases}
\end{align} 

It remains then to estimate the term \eqref{comm:eqnorm}. By using the proof of Claim 7 in \cite{kenig2011local}, we know that for any $\alpha>0$, there exists $C>0$ such that for all $u, \, v \in \mathcal{S}(\R)$,
	\begin{align}
		\bigg|\int (\big[\Phi ,[\Phi ,\chi(D)|D|^{\alpha}]\big]u) v\bigg|\leq C\int(u^2+v^2)|\phi'|. \label{k4}
	\end{align}
 We observe that the symbol of the pseudo-differential operator $\Phi^{-1}\left[\Phi ,\left[\Phi ,\chi(D) |D|^{2\alpha} \right]\right]$ belongs to $\mathcal{S}^{2\alpha-2,-\frac{1+\alpha}{2}-2}\subset \mathcal{S}^{\alpha,-\frac{1+\alpha}{2}}$, since $\alpha \in (0,2)$. Thus it follows from \eqref{nps}, \eqref{k4} and Young's inequality that
\begin{align*}
	\bigg|\int u \left[\Phi ,\left[\Phi , |D|^{2\alpha} \right]\right]u \bigg|&\leq \bigg|\int u \Phi \Phi^{-1} \left[\Phi ,\left[\Phi ,(1-\chi(D)) |D|^{2\alpha} \right]\right]u \bigg| \notag\\
	&\quad +\bigg|\int u \left[\Phi ,\left[\Phi ,\chi(D) |D|^{2\alpha} \right]\right]u  \bigg|\notag\\
	&\leq \frac{C}{A^{\alpha}}\int \left(\langle D \rangle^{\alpha} u\right)^2 |\phi'| + C A^{\alpha} \int u^2 |\phi'|. 	
\end{align*}
Moreover, arguing as in \eqref{dec}, we obtain that
\begin{align*}
	\int \left(\langle D \rangle^{\alpha} u\right)^2 |\phi'| dx
	&\leq C\left( \int  u^2 |\phi'|  + \int \left(|D|^{\alpha} u\right)^2 |\phi'| \right) .
\end{align*} 
By changing variable $x'=\frac{x-m_j}{A}$, we deduce that
\begin{align}
	\bigg|\int u\Phi_{j,A}\left[\Phi_{j,A},|D|^{2\alpha}\right]u  \bigg| \leq \frac{C}{A^{\alpha}}\int \left( u^2 + \left(|D|^{\alpha} u\right)^2 \right) |\phi'_{j,A}| \label{a2}
\end{align}
Therefore, by gathering \eqref{a1}, \eqref{a2}, we conclude the proof of \eqref{eqnorm1}.

\end{proof} 

\begin{proof}[Proof of Lemma \ref{esttc}]
By Young's inequality and  \eqref{eqnorm1}, we obtain that  
\begin{align*}
	  \bigg|\int |D|^{\alpha}\left(u\Phi_{j,A}  \right) (|D|^{\alpha}u)\Phi_{j,A}  \bigg| 
	\leq&  \frac12 \left( \int \left(|D|^{\alpha}\left(u\Phi_{j,A} \right)\right)^2 + \int \left(|D|^{\alpha}u\right)^{2}|\phi'_{j,A}|   \right) \notag \\ \leq & \int \left(|D|^{\alpha} u\right)^2|\phi'_{j,A}|+ \frac12\left|\int \left(|D|^{\alpha}\left(u\Phi_{j,A} \right)\right)^2 -  \int \left( |D|^{\alpha}u \right)^{2} |\phi'_{j,A}| \right| \notag \\
	\leq& \int \left(|D|^{\alpha} u\right)^2|\phi'_{j,A}|+\frac{C}{A^{\frac{\alpha}{2}}}  \left(\int \left( u^2 + \left(|D|^{\frac{\alpha}{2}}u\right)^2 + \left(|D|^{\alpha}u\right)^2 \right) |\phi'_{j,A}|\right)  ,
\end{align*}
which yields \eqref{est:esttc}.
\end{proof}

\section{Proof of the non-linear weighted estimates} \label{E}

\begin{proof}[Proof of Lemma \ref{nonlinterm}]
	Let $j\in\{1,..., N\}$. First we prove estimate \eqref{nlt3}. By the Cauchy-Schwarz inequality, the Sobolev embedding $\dot{H}^{\frac{1}{4}}(\R) \xhookrightarrow{} L^4(\R)$ and $\frac{1}{4}<\frac{\alpha}{2}$, we get that
	\begin{align}
		\displaystyle\int |\eta|^3|\phi_{j,A}'|&\leq \left( \int \eta^2 \right)^{\frac{1}{2}} \left( \int \eta^4|\phi_{j,A}'|^2 \right)^{\frac{1}{2}}
		\leq C \gamma \lVert |D|^{\frac{1}{4}}\left(\eta \Phi_{j,A}\right) \rVert_{L^2}^{2} 
		\leq C \gamma \lVert \eta \Phi_{j,A} \rVert_{H^{\frac{\alpha}{2}}}^{2}. \label{nlt3est}
	\end{align} 
	From the decomposition $\eta=u-R$, we have
	\begin{align}
		\lVert \eta \Phi_{j,A} \rVert_{H^{\frac{\alpha}{2}}} \leq \|u \Phi_{j,A}  \|_{H^{\frac{\alpha}{2}}} + \| R \Phi_{j,A}\lVert_{H^{\frac{\alpha}{2}}} \leq  \|u \Phi_{j,A}  \|_{H^{\frac{\alpha}{2}}} + \| R \Phi_{j,A} \lVert_{H^{\alpha}} .\label{q1}
	\end{align}
	
	To deal with the second term on the right-hand side of \eqref{q1}, we use pseudo-differential calculus. Observe that the symbols of $\Phi\chi(D)\langle D \rangle^{\alpha}\Phi^{-1}$ and  $\Phi(1-\chi(D))\langle D \rangle^{\alpha}|D|^{-\alpha}\Phi^{-1}$ belong to $\mathcal{S}^{0,0}$. It follows then from 
	\eqref{oh}, and then \eqref{pps} that, for all $v\in\mathcal{S}(\R)$,
	\begin{align}
		\||D|^{\alpha}\left(v \Phi \right)  \|_{L^2}&\leq \| v \Phi \|_{H^{\alpha}} \notag\\
		&\leq C \|\left(\langle D \rangle^{\alpha}v\right)\Phi \|_{L^2} \notag\\
		&\leq C \| \chi(D)\left(\langle D \rangle^{\alpha}v\right)\Phi \|_{L^2}
		+\|(1-\chi(D))\langle D \rangle^{\alpha}|D|^{-\alpha}\left(|D|^{\alpha}v\right)\Phi \|_{L^2} \notag \\
		&\leq C\left( \| v\Phi \|_{L^2} +\|\left(|D|^{\alpha}v\right)\Phi \|_{L^2}\right). \label{est:split}
	\end{align}
	Then, we obtain, by changing variable $x'=\frac{x-m_j}{A}$,
	\begin{align*}
		\||D|^{\alpha}\left(v \Phi_{j,A} \right)  \|_{L^2}\leq C\left( \frac{1}{A^{\alpha}}\| v\Phi_{j,A}\|_{L^2} +\|\left(|D|^{\alpha}v\right)\Phi_{j,A} \|_{L^2}\right) ,
	\end{align*}
	so that
	\begin{align}
		\|v \Phi_{j,A}  \|_{H^{\alpha}} \leq C\left( \| v\Phi_{j,A} \|_{L^2} +\|\left(|D|^{\alpha}v\right)\Phi_{j,A} \|_{L^2}\right) ,\label{q2}
	\end{align}
	for all $v\in \mathcal{S}(\R)$.
	
	Therefore, we deduce combining \eqref{nlt3est}, \eqref{q1} and \eqref{q2} with $v=R$ that
	\begin{align*}
		\displaystyle\int |\eta|^3|\phi_{j,A}'| \leq C\gamma \| u\Phi_{j,A}\|_{H^{\frac{\alpha}{2}}}^2 + C \gamma \int \left(   R^2 + \left(|D|^{\alpha}R\right)^2  \right)|\phi_{j,A}'|.
	\end{align*}
	Moreover, by using the equation  \eqref{eqQG2} and \eqref{est4}, we obtain 
	\begin{align*}
		\int \left(|D|^{\alpha}R\right)^2|\phi'_{j,A}| \leq C\left( \int R^2+ R^4 \right) |\phi'_{j,A}|\leq \frac{C}{(\beta t)^{1+\alpha}},
	\end{align*}
	which concludes the proof of \eqref{nlt3}.
	
	Now we prove  \eqref{nlt4}. Using the Cauchy-Schwarz inequality, the Sobolev embedding and the former estimates, we conclude that
	\begin{align*}
		\displaystyle\int \eta^4|\phi'_{j,A}|&\leq \left( \int \eta^4 \right)^{\frac{1}{2}} \left( \int \eta^4|\phi'_{j,A}|^2 \right)^{\frac{1}{2}}
		\leq C\gamma^2 \left[ \int u^2|\phi'_{j,A}| + \int  \left(|D|^{\frac{\alpha}{2}}\left(u \Phi_{j,A} \right) \right)^2 \right] + \displaystyle\frac{C\gamma^2}{\left(\beta t\right)^{1+\alpha}}  ,
	\end{align*} 
	which is exactly estimate \eqref{nlt4}.
\end{proof}

\begin{proof}[Proof of Lemma \ref{esttcnl}]

By using Young's inequality and the decomposition $u= R + \eta $,  we deduce that
\begin{align*}
	\bigg|\int |D|^{\frac{\alpha}{2}}\left(u\Phi_{j,A} \right) |D|^{\frac{\alpha}{2}}\left(u^2\Phi_{j,A} \right)   \bigg|=&\bigg|\int |D|^{\alpha}\left(u\Phi_{j,A} \right) u^2\Phi_{j,A}   \bigg|\\
	\leq& 2\int u^4|\phi'_{j,A}|  + \frac{1}{8}\int \left(|D|^{\alpha}\left(u\Phi_{j,A}\right) \right)^2\\
	\leq& C\left( \int \eta^4 |\phi'_{j,A}| + \int R^4|\phi'_{j,A}| \right) +  \frac{1}{8}\int \left(|D|^{\alpha}\left(u\Phi_{j,A}\right) \right)^2.
\end{align*}
Using \eqref{nlt4}, \eqref{est4}, \eqref{eqnorm1}, we have that 
\begin{align*}
	\bigg|\int |D|^{\frac{\alpha}{2}}\left(u\Phi_{j,A}  \right) |D|^{\frac{\alpha}{2}}\left(u^2\Phi_{j,A}  \right)   \bigg|\leq& C\left(\gamma^2 + \frac{1}{A^{\frac{\alpha}{2}}} \right) \left(\int \left( u^2 + \left(|D|^{\frac{\alpha}{2}}u\right)^2 + \left(|D|^{\alpha}u\right)^2 \right) |\phi'_{j,A}|\right) \notag  \\
	&+ C\gamma^2 \int \left(|D|^{\frac{\alpha}{2}}\left(u \Phi_{j,A}\right)\right)^2 + \frac{1}{8} \int \left(|D|^{\alpha} u\right)^2|\phi'_{j,A}| \notag \\
	&+ \frac{C}{\left(\beta t\right)^{1+\alpha} }. 
\end{align*}
Furthermore, by using again \eqref{eqnorm1} with $\frac{\alpha}{2}<1$, we deduce that 
\begin{align}
	\bigg|\int |D|^{\frac{\alpha}{2}}\left(u\Phi_{j,A}  \right) |D|^{\frac{\alpha}{2}}\left(u^2\Phi_{j,A}  \right)   \bigg|\leq& C\left(\gamma^2 + \frac{1}{A^{\frac{\alpha}{2}}} \right) \left(\int \left( u^2 + \left(|D|^{\frac{\alpha}{2}}u\right)^2 + \left(|D|^{\alpha}u\right)^2 \right) |\phi'_{j,A}|\right) \notag  \\
	&+ \frac{1}{8} \int \left(|D|^{\alpha} u\right)^2|\phi'_{j,A}| + \frac{C}{\left(\beta t\right)^{1+\alpha} }. \label{cl341}
\end{align}
This concludes the proof of Lemma \ref{esttcnl}.
\end{proof}


\section{Coercivity of the localized operator}

To begin, we recall the definition of $H_j$ given in \eqref{quadratic}
$$
H_j(u,u)=\int \left( u|D|^{\alpha}u +c_ju^2 -2R_ju^2  \right)\psi_{j,A},
$$
where $R_j$ is defined in \eqref{Rg} and $\psi_{j,A}$ is defined in \eqref{psi}, and $u\in H^{\frac{\alpha}{2}}(\R)$. Moreover, let  
$$
\mathcal{L}_ju= |D|^{\alpha}u +c_j u -2R_ju \quad \text{and} \quad \mathcal{L}u= |D|^{\alpha}u + u -2Q u .
$$ 

It has been proved in  Theorem 2.3 in \cite{frank2013uniqueness} that the spectrum of $\mathcal{L}$  is composed by one negative eigenvalue, the eigenvalue $0$ and that the rest is the continuous spectrum $[1,+\infty)$. Moreover, the eigenspaces associated with the negative eigenvalue and $0$ are one-dimensional vector spaces and the eigenspace of $0$ is spanned by $Q'$. 

Furthermore, from Lemma E.1 in \cite{weinstein1985modulational}, since we are in the subcritical case, we can replace the eigenfunction associated with the negative eigenvalue by $Q$ to get the coercivity property stated in the following theorem. 
\begin{thm}\label{Lcoercivite} Let $\alpha \in\left(\frac12,2\right)$. 
	Then, there exists $\mu >0$ such that for all $u\in H^{\frac{\alpha}{2}}(\R)$
	$$
	\displaystyle\int u\mathcal{L}u\geq \mu\| u \|_{H^{\frac{\alpha}{2}}}  -\frac{1}{\mu}\left(\int uQ  \right)^2 - \frac{1}{\mu}\left(\int uQ'  \right)^2 
	$$
\end{thm}

\begin{remark}
By using a scaling argument, the result of Theorem \ref{Lcoercivite} still holds if one replaces $\mathcal{L}$ by $\mathcal{L}_j$ and $Q$ by $R_j$, for $j \in \{1,\cdots,N\}$. 
\end{remark}

As a consequence of the former theorem, we deduce a coercivity property for the bilinear form $H_j$.
\begin{corollary}\label{almostcoer}
	There exist $\nu>0$, $C>0$ such that for all $A>1$, $u\in H^{\frac{\alpha}{2}}(\R)$
	\begin{align*}
		\displaystyle\sum_{j=1}^{N}H_j(u,u)\geq&  \left(\nu - \frac{C}{\left( \beta t \right)^{\alpha}} - \frac{C}{A^{\frac{\alpha}{2}}} \right)\| u \|_{H^{\frac{\alpha}{2}}}^2 -\frac{1}{\nu}\sum_{j=1}^N \left(\left(\int uR_j  \right)^2 +\left(\int u\partial_x R_j  \right)^2 \right) .
	\end{align*}
\end{corollary}
\begin{proof}[Proof of Corollary \ref{almostcoer}]
	For all $j \in \{1,\cdots,N\}$, we have from Theorem \ref{Lcoercivite}  that
	\begin{align*}
		H_j(u,u)=&\int u\sqrt{\psi_{j,A}} \mathcal{L}(u\sqrt{\psi_{j,A}}) + \int u \left( |D|^{\alpha}u \right)\psi_{j,A} -\int \left(|D|^{\frac{\alpha}{2}}\left(u\sqrt{\psi_{j,A}}\right)  \right)^2\\
		\geq& \mu \| u\sqrt{\psi_{j,A}} \|_{H^{\frac{\alpha}{2}}}^2 -\frac{1}{\mu}\left(\int u\sqrt{\psi_{j,A}}R_j  \right)^2 - \frac{1}{\mu}\left(\int u\sqrt{\psi_{j,A}}\partial_x R_j  \right)^2\\
		&+ \int u \left(|D|^{\alpha}u\right) \psi_{j,A} -\int \left(|D|^{\frac{\alpha}{2}}\left(u\sqrt{\psi_{j,A}}\right) \right)^2.\\
	\end{align*}
By \eqref{est5}, \eqref{est6}, and the Cauchy-Schwarz inequality we deduce that for all $j\in\{1,\cdots,N\}$ 
\begin{align*}
	\left(\int u\sqrt{\psi_{j,A}} R_j  \right)^2 +\left(\int u \sqrt{\psi_{j,A}}\partial_x R_j  \right)^2&\leq 2\left(\int uR_j  \right)^2 + 2 \left(\int u\partial_x R_j  \right)^2 \\
	&\quad+2\left(\int u(1-\sqrt{\psi_{j,A}}) R_j  \right)^2 +2\left(\int u (1-\sqrt{\psi_{j,A}})\partial_x R_j  \right)^2\\
	&\leq 2\left(\int uR_j  \right)^2 + 2 \left(\int u\partial_x R_j  \right)^2 + \frac{C\|u\|_{L^2}^2}{(\beta t)^{\alpha}}.
\end{align*}
Observe from $\langle D \rangle^{\frac{\alpha}{2}}\sim 1+|D|^{\frac{\alpha}{2}}$ that 
\begin{equation*} 
\|u\sqrt{\psi_{j,A}}\|_{H^{\frac\alpha2}}^2 \ge c \int (u^2+(|D|^{\frac\alpha2}u)^2)\psi_{j,A} -c \int (|D|^{\frac\alpha2}u)^2\psi_{j,A}+c\int \left(|D|^{\frac{\alpha}{2}}\left(u\sqrt{\psi_{j,A}}\right) \right)^2 ,
\end{equation*}
for a small positive constant $0<c<1$. Since $\displaystyle\sum_{j=1}^{N} \psi_{j,A}=1$, we have
	\begin{align} \label{sum_psi}
		\sum_{j=1}^N \int u\left(|D|^{\alpha}u\right)\psi_{j,A} = \int \left(|D|^{\frac{\alpha}{2}}u\right)^2= \sum_{j=1}^N \int \left(\left(|D|^{\frac{\alpha}{2}}u\right)\sqrt{\psi_{j,A}}\right)^2.
	\end{align}
Hence, we deduce by summing over $j$ that 
\begin{align*}
\sum_{j=1}^{N}H_j(u,u)		\geq&  N  \left(c\mu - \frac{C}{(\beta t)^{\alpha}}\right)  \| u \|_{H^{\frac{\alpha}{2}}}^2 -\frac{2}{\mu}\sum_{j=1}^N  \left(\left(\int u R_j  \right)^2 +\left(\int  u\partial_x R_j  \right)^2\right)\\
&+\left(1-c\mu\right) \sum_{j=1}^N \left( \int u \left(|D|^{\alpha}u\right) \psi_{j,A} -\int \left(|D|^{\frac{\alpha}{2}}\left(u\sqrt{\psi_{j,A}}\right) \right)^2 \right).	
\end{align*}
It remains to estimate the last term on the right hand side of the former inequality.
		By using \eqref{sum_psi} and direct computations, 
	\begin{align*}
		&\sum_{j=1}^N\left( \int u\left(|D|^{\alpha}u\right)\psi_{j,A} - \int \left(|D|^{\frac{\alpha}{2}}\left(u\sqrt{\psi_{j,A}}\right) \right)^2\right) \\
		&=\sum_{j=1}^N \int \left(\left(|D|^{\frac{\alpha}{2}}u\right)\sqrt{\psi_{j,A}} + |D|^{\frac{\alpha}{2}}\left(u\sqrt{\psi_{j,A}}\right) \right) \bigg( \left(|D|^{\frac{\alpha}{2}}u\right)\sqrt{\psi_{j,A}} - |D|^{\frac{\alpha}{2}}\left(u\sqrt{\psi_{j,A}}\right) \bigg)\\
		&=-\sum_{j=1}^N \int \left( \left(|D|^{\frac{\alpha}{2}}u \right)\sqrt{\psi_{j,A}} + |D|^{\frac{\alpha}{2}}\left(u\sqrt{\psi_{j,A}}\right) \right) [|D|^{\frac{\alpha}{2}},\sqrt{\psi_{j,A}}] u. \\						
	\end{align*}
	By arguing as in \eqref{est:split} and using $0 \le \psi_{j,A} \le 1$, we deduce that 
	\begin{align}
		\sum_{j=1}^N \| |D|^{\frac{\alpha}{2}}(u\sqrt{\psi_{j,A}})\|_{L^2}\leq  C\sum_{j=1}^N \left(\| u\sqrt{\psi_{j,A}} \|_{L^2} + \| \left(|D|^{\frac{\alpha}{2}}u\right) \sqrt{\psi_{j,A}}\|_{L^2} \right) \leq C\| u \|_{H^{\frac{\alpha}{2}}}. \label{appq1}
	\end{align}
	Then, it follows from the Cauchy-Schwarz inequality and \eqref{appq1} that
	\begin{align*}
		\sum_{j=1}^N\bigg| \int u\left(|D|^{\alpha}u\right)\psi_{j,A} - \int \left(|D|^{\frac{\alpha}{2}}\left(u\sqrt{\psi_{j,A}}\right) \right)^2 \bigg|&\leq C \| u \|_{H^{\frac{\alpha}{2}}} \sum_{j=1}^N\|[|D|^{\frac{\alpha}{2}},\sqrt{\psi_{j,A}}] u \|_{L^2} 
	\end{align*}
	Finally to estimate the commutator on the right-hand side of the former estimate, we will rely on pseudo-differential calculus and argue as in the previous subsection. 
 By \eqref{comps}, we have that the symbol  of  $[|D|^{\frac{\alpha}{2}}(1-\chi(D)),\sqrt{\psi}]$ belongs to $\mathcal{S}^{\frac{\alpha}{2}-1,-1}\subset \mathcal{S}^{0,0}$, since $\alpha<2$. Then, it follows from \eqref{nps} that
	$$
	\|[|D|^{\frac{\alpha}{2}},\sqrt{\psi}] u \|_{L^2} \leq C \left( \| u \|_{L^2} + \|[|D|^{\frac{\alpha}{2}}\chi(D),\sqrt{\psi}] u \|_{L^2}  \right).
	$$
	We recall $|D|^{\frac{\alpha}{2}}\chi(D)u=k\ast u$, with $\hat{k}=|\xi|^{\frac{\alpha}{2}}\chi(\xi)$, so that 
	$$[\chi(D)|D|^{\frac{\alpha}{2}},\sqrt{\psi}]u= \int k(x-y)\left(\sqrt{\psi}(x)-\sqrt{\psi}(y) \right) u(y)dy .$$ 
	We want to prove that the operator $T$ defined by the kernel 
	$$
	R(x,y)=k(x-y)\left(\sqrt{\psi}(x)-\sqrt{\psi}(y) \right),
	$$
	is bounded in $L^2(\R)$. By Lemma \ref{K2}, we obtain that $\displaystyle |k(x)|\leq \frac{C}{\langle x \rangle^{1+\frac{\alpha}{2}}}$. Since  $\sqrt{\psi}\in L^{\infty}(\R)$, we deduce by Lemma \ref{K3} that
	$$
	\|[|D|^{\frac{\alpha}{2}}\chi(D),\sqrt{\psi}] u \|_{L^2} \leq C \| u \|_{L^2} .
	$$
	By changing the variable $x'=\frac{x-m_j}{A}$, we get that 
	$$
	\|[|D|^{\frac{\alpha}{2}},\sqrt{\psi_{j,A}}] u \|_{L^2} \leq \frac{C}{A^{\frac{\alpha}{2}}} \| u \|_{L^2} .
	$$
	
	We finish the proof of Corollary \ref{almostcoer} by combining all these estimates and by choosing $\nu$ small enough.  
\end{proof}

\section*{Acknowledgments} 
This research was supported by a Trond Mohn Foundation grant. The author is grateful for the anonymous referees for their comments which greatly improved this article. 
The author would like to thank Jacek Jendrej for pointing out the qualitative implicit function theorem, Yvan Martel for many important comments, Luc Molinet for the proof of the weak continuity of the flow and  his advisor Didier Pilod for suggesting this problem, for his constant scientific support, and whose pertinent comments greatly improved this manuscript.

\bibliographystyle{amsplain}
\bibliography{biblio}

\end{document}